\documentclass[11pt]{article}
 
\newcommand{\mb}[1]{\mathbb{#1}}

\newcommand{\bs}[1]{\boldsymbol{#1}}
\newcommand{\cbrac}[1]{\left\{#1\right\}}
\newcommand{\sbrac}[1]{\left[#1\right]}

\usepackage{subfigure}
\usepackage{indentfirst}
\usepackage{xcolor}
\usepackage{verbatim}
\usepackage{tikz}
\usepackage{hyperref}
\usepackage{bm}
\usepackage[a4paper, total={6in, 8in}]{geometry}
\usepackage{amsmath}
\usepackage{amsthm}	
\usepackage{amsfonts}	
\usepackage{amssymb, bbm, units}
\usepackage{enumerate}
\usepackage[nobysame, alphabetic, initials]{amsrefs}
\usepackage{ulem} 

\newcommand{\norm}[1]{\left\lVert#1\right\rVert}
\numberwithin{equation}{section}
\numberwithin{figure}{section}
 \usepackage[nodayofweek]{datetime}

\theoremstyle{plain}
\newtheorem{theorem}{Theorem}[section]
\newtheorem{lemma}[theorem]{Lemma}
\newtheorem{proposition}[theorem]{Proposition}
\newtheorem{corollary}[theorem]{Corollary}

\newtheorem{definition}[theorem]{Definition}
\newtheorem{remark}[theorem]{Remark}

\newcommand{\cF}{\mathcal{F}}
\newcommand{\cG}{\mathcal{G}}

\newcommand{\cL}{\mathcal{L}}

\newcommand{\cT}{\mathcal{T}}

\newcommand{\cZ}{\mathcal{Z}}

\newcommand{\DD}{\mathbb{D}}
\newcommand{\EE}{\mathbb{E}}

\newcommand{\MM}{\mathbb{M}}

\newcommand{\PP}{\mathbb{P}}
\newcommand{\RR}{\mathbb{R}}
\newcommand{\SSS}{\mathbb{S}}
\newcommand{\ZZ}{\mathbb{Z}}

\newcommand{\1}{\mathbbm{1}}
\newcommand{\brac}[1]{\left(#1\right)}
\newcommand{\mf}[1]{\mathbf{#1}}
\newcommand{\mft}[1]{\mathbf{\tilde{#1}}}
\newcommand{\define}{\triangleq}
\newcommand{\indic}[1]{\mathbbm{1}{\brac{#1}}}

\newcommand\numeq[1]%
  {\stackrel{\scriptscriptstyle(\mkern-1.5mu#1\mkern-1.5mu)}{=}}

\makeatletter
\newsavebox{\@brx}
\newcommand{\llangle}[1][]{\savebox{\@brx}{\(\m@th{#1\langle}\)}%
  \mathopen{\copy\@brx\mkern2mu\kern-0.9\wd\@brx\usebox{\@brx}}}
\newcommand{\rrangle}[1][]{\savebox{\@brx}{\(\m@th{#1\rangle}\)}%
  \mathclose{\copy\@brx\mkern2mu\kern-0.9\wd\@brx\usebox{\@brx}}}
\makeatother

\newcommand{\toP}{\stackrel{p}{\rightarrow}}

\definecolor{darkgreen}{rgb}{0,0.35,0}

\title{Asymptotic Optimality of the Speed-Aware Join-the-Shortest-Queue in the Halfin-Whitt Regime for Heterogeneous Systems}
\author{Sanidhay Bhambay$^1$ \and  Burak B\"{u}ke$^2$ \and Arpan Mukhophadyay$^3$}
\date{%
    $^1$Durham University Business School\\%
    $^2$School of Mathematics, The University of Edinburgh\\%
    $^3$Department of Computer Science, University of Warwick
}
\begin{document}
\maketitle
\begin{abstract}
The Join-the-Shortest-Queue (JSQ) load balancing scheme is known to minimise the average response time of jobs in homogeneous systems with identical servers.  However, for {\em heterogeneous} systems with servers having different processing speeds, finding an optimal load balancing scheme remains an open problem for finite system sizes. Recently, for systems with heterogeneous servers, a variant of the JSQ scheme, called the {\em Speed-Aware-Join-the-Shortest-Queue (SA-JSQ)} scheme, has been shown to achieve asymptotic optimality in the fluid-scaling regime where the number of servers $n$ tends to infinity but the normalised the arrival rate of jobs remains constant.
{In this paper, we show that the SA-JSQ scheme is also asymptotically optimal for heterogeneous systems in the {\em Halfin-Whitt} traffic regime where the normalised arrival rate scales as $1-O(1/\sqrt{n})$.} Our analysis begins by establishing that an appropriately scaled and centered version of the Markov process describing system dynamics  weakly converges to a two-dimensional reflected {\em Ornstein-Uhlenbeck (OU) process}. We then show using {\em Stein's method} that the stationary distribution of the underlying Markov process converges to that of the OU process as the system size increases by establishing the validity of interchange of limits. {Finally, through coupling with a suitably constructed system, we show that SA-JSQ asymptotically minimises the diffusion-scaled total number of jobs and the diffusion-scaled number of waiting jobs in the steady-state in the Halfin-Whitt regime among all policies which dispatch jobs based on queue lengths and server speeds.}
\end{abstract}

\section{Introduction}
Load balancing is key to optimising resource utilisation and ensuring low delay of jobs in  multi-server systems. In modern data centers, where numerous servers of various capacities coexist, load balancers play a pivotal role in distributing incoming network traffic or requests across these servers. This distribution of the load prevents any single server from being overwhelmed with jobs when there are idle servers in the system. Hence, load balancing ensures high utilisation of the available resources as well as low mean response time for the jobs. 

The standard model to study the performance of load balancing schemes involves $n$ {\em identical} or {\em homogeneous} servers, each equipped with its own queue to store pending job requests. In this setup, there is a continuous stream of jobs arriving at a rate of $n\lambda^{(n)}$. In contrast to the well-known $M/M/n$ system, there is no global queue in this model and the dispatcher routes all incoming jobs to a server queue immediately upon arrival according to a specific load balancing scheme. This approach is especially beneficial to model the limited memory available to dispatchers in data centers for storing incoming jobs. A natural load balancing scheme to consider in this canonical setting is the {\em Join-the-Shortest-Queue (JSQ)} scheme, wherein each incoming job is assigned to the server with minimum queue length with ties broken uniformly at random. 
The JSQ scheme was first analysed in~\cite{Winston1977optimality} and its optimality 
in terms of minimising the mean response time of jobs was proved for Poisson arrivals and exponential service times. This result was later extended to service time distributions having nondecreasing hazard rate by \cite{Weber1978}. The transient dynamics of the JSQ scheme has been studied by \cite{Eschenfeldt2018} in the {\em Halfin-Whitt regime}, introduced originally for $M/M/n$ systems by \cite{halfin1981heavy}, by proving diffusion limits of the scaled system size processes when $\lambda^{(n)}$ is dependent on $n$, and the quantity $\sqrt{n}(1-\lambda^{(n)})$ has a nondegenerate limit $\beta>0$ as $n\to \infty$. \cite{Braverman2020} showed that this diffusion limit can also be used to study the steady-state behaviour of the model by proving an interchange of limits result. The rate at which the steady-state distributions of the JSQ scheme converges to the steady-state distribution of the diffusion limit was derived in~\cite{braverman2023join}

It is crucial to highlight that the effectiveness of the JSQ scheme primarily relies on the  server homogeneity assumption. However, this assumption does not align with practical scenarios since data centers typically host an array of physical devices spanning multiple generations, having varying processing capabilities. Furthermore, the processing speeds of servers can diverge due to the incorporation of diverse acceleration devices like GPUs, FPGAs, and ASICs (see e.g.\cite{GPU_het, FPGA_het}). In such heterogeneous systems, job assignment strategies like JSQ and the Power-of-d-choices (Pod), which were originally designed with homogeneous systems in mind, may exhibit notably suboptimal performance~\cite{gardiner_perf, bramson2012asymptotic, arpan_tcns}. Additionally, it is worth noting that in finite systems, the classical JSQ scheme is recognised as non-optimal when servers exhibit heterogeneity \cite{Krishnan1987joining, HYYTIA2017}. Therefore, for heterogeneous systems, speed-aware schemes that   assign incoming jobs based on both the queue lengths and the server speeds are essential to minimize the average response time of jobs.

In a  recent line of work by~\cite{Weng2020, bhambay2022asymptotic, debankur_constrained_2021}, a speed-aware variant of the JSQ scheme called the {\em Speed-Aware-Join-the-Shortest-Queue (SA-JSQ)}, also referred as {\em Join-the-Fastest-Shortest-Queue (JFSQ)}, has been introduced for systems with $n$ {\em heterogeneous} servers categorised into $M$ different {\em types} or {\em pools}. In this heterogeneous model, each server pool $j \in \{1,2,\ldots,M \}$ contains $O(n)$ servers with processing speed $\mu_j$ where $\mu_1 > \mu_2>\ldots>\mu_M$. Under the SA-JSQ scheme, incoming jobs are assigned to servers with the highest speed among those with the shortest queue length. Notably, ~\cite{bhambay2022asymptotic} prove that the SA-JSQ scheme asymptotically minimizes the average response time of jobs when $\lambda^{(n)}=\lambda\in[0,1)$ for each $n$ in the fluid limit. 
In the present work, we study the SA-JSQ load balancing scheme for heterogeneous systems in the {Halfin-Whitt regime}.  
Specifically, we study both the transient and the steady-state behaviour of the SA-JSQ scheme and show that SA-JSQ is also asymptotically optimal in the Halfin-Whitt regime among all load balancing policies which dispatch jobs based on queue lengths and server speeds. Our contributions are detailed below.
\subsection{Contributions}

{The main contributions of this paper are three-fold: (i) We present a transient analysis of a  system with heterogeneous servers and dedicated queues operating in the Halfin-Whitt regime under the SA-JSQ policy. Specifically, we prove that the appropriately scaled and centred queue lengths converge to a two-dimensional reflected Ornstein-Uhlenbeck process, exhibiting a state-space collapse property. (ii) We demonstrate that this limiting Ornstein-Uhlenbeck process has a stationary distribution which can be used to approximate the stationary behaviour of the system as the number of servers becomes large. (iii) We establish that under the Halfin-Whitt regime the SA-JSQ policy stochastically minimises both the diffusion-scaled total number of jobs and the diffusion-scaled number of waiting jobs in the system in the limit as $n\to \infty$. While our transient and stationary analyses broadly follow existing analyses for homogeneous systems with significant technical differences introduced due to the increased complexity of the underlying state-space, to the best of our knowledge, we provide the first
proof of asymptotic optimality of a JSQ-type load balancing policy in the Halfin-Whitt regime for a system where the servers are heterogeneous and each server has a dedicated queue. We note that the general problem of finding optimal load balancing policies for a system of heterogeneous parallel queues is a well known open problem in queuing theory while its homogeneous counterpart was solved almost four decades ago by~\cite{Winston1977optimality,Weber1978}. Our result, therefore, can be considered as an important step toward the general understanding of efficient policies for heterogeneous systems.}

{\textit{\textbf{Transient Analysis}}: 
To establish the diffusion limit for the conventional JSQ policy in homogeneous systems, \cite{Eschenfeldt2018} employs a three-step approach: (i) constructing a truncated system that differs from the original only in that the queues cannot exceed a certain bound determined by the limiting dynamics of the original system, (ii) proving the diffusion limit of the truncated system via an appropriate Skorohod mapping, and (iii) demonstrating that the difference between the truncated system and the original system vanishes as $n \to \infty$. We generalise this approach to heterogeneous systems. However, there are several technical challenges in doing so.}

{First, unlike homogeneous systems, the queue lengths in different server pools scale differently, resulting in a higher-dimensional state space and necessitating a more complex Skorohod mapping to establish weak convergence. Second, the truncated system used in our proof differs from the one used in \cite{Eschenfeldt2018}. Specifically, the buffer sizes are not uniform across all queues in our truncated system: While the queues in the fastest server pool can hold at most two jobs, the queues in the remaining pools are restricted to holding no more than one job. Finally, we show that under the SA-JSQ policy, the heterogeneous system exhibits a state space collapse, where the diffusion-scaled number of idle servers in the fastest 
$M-1$ pools converges to zero for all finite times (Proposition~\ref{prop:state_space_collapse}). This result has no direct analogue in homogeneous systems. }

{\textit{\textbf{Steady-State Analysis}}: 
The proof of the \emph{interchange of limits} is carried in two essential steps:  (i) by proving that the sequence of stationary distributions of the diffusion-scaled queue length process is tight (Theorem~\ref{thm:stationary_convergence}), and (ii) by proving that the limiting diffusion process is (exponentially) ergodic (Theorem~\ref{thm:stationarity_diffusion}). The proofs of both theorems use Stein's method which involves applying the generator of the underlying Markov process to a suitably chosen Lyapunov function. 
The key challenge in this approach is to identify a suitable Lyapunov function that will yield appropriate bounds on the expected value of the queue lengths in stationarity. The Lyapunov function is derived as the solution to a partial differential equation obtained through a generator expansion approach. In this method, the generator is approximated via a Taylor expansion, and it is shown that the error terms vanish using the gradient bounds of the solution. In a recent work, \cite{Braverman2020} used this approach to study the behaviour of the JSQ scheme for homogeneous systems in the Halfin-Whitt regime.}

{Even though we follow the framework in \cite{Braverman2020} closely to prove the interchange of limits, the complexity of the underlying state-space again poses several challenges in the heterogeneous setting. First, to obtain a tractable PDE, we need to introduce a suitable lifting operator that projects the state of the Markov process into $\RR^2$. In the homogeneous setting, this is done naturally by only considering the scaled total number of idle servers and the scaled total number of servers with two or more jobs. In the heterogeneous setting, however, we identify that the correct way to define the lifting operator is by setting the second dimension to  the scaled number of servers with two or more jobs \emph{only} in the fastest pool as opposed to that in the entire system. Once the PDE is obtained using this lifting operator, the next challenge is to solve it to find the Lyapunov function and the gradient bounds. The heterogeneous nature of the system results in a PDE different from the one used for homogeneous systems and this results in a notable difference in the analysis required to solve the PDE.  Moreover, the resulting error terms in the generator expansion are significantly more complicated in the heterogeneous setting and one needs to prove additional state space collapse results to derive bounds on these terms. }

{\textit{\textbf{Asymptotic Optimality}}: 
Finally, we prove the asymptotic optimality of the SA-JSQ scheme in terms of stochastically minimising the diffusion-scaled total number of jobs and the diffusion-scaled number of waiting jobs in the steady-state in the Halfin-Whitt regime (Theorem~\ref{thm:sajsq_optimality}). To do so, we first derive a lower bound on the performance achievable by any {\em admissible load balancing scheme} defined as a scheme that considers only the current queue lengths and server speeds to assign each incoming job in the original system (Proposition~\ref{prop:coupling_modf}). This lower bound is established by coupling the original system with a modified system where servers are not attached to queues, but are free to move between queues at any time, allowing faster servers to always serve longer queues. To complete the proof of asymptotic optimality, we finally prove that the modified system described above shares the same diffusion limit as the original system operating under the SA-JSQ scheme (Proposition~\ref{prop:modified_diffusion_limit}). We note that the modified system we use in our proofs differs fundamentally from the $M/M/n$ system having one central queue used in the proofs of optimality in prior works of~\cite{armony2005dynamic} and \cite{tezcan2008optimal}. This difference stems from the fact that the diffusion limit of an $M/M/n$ system with a single queue is one dimensional, while, in the parallel queue setting, both idle servers and queues with waiting jobs can co-exist yielding a two-dimensional diffusion limit.}

\subsection{Related Literature}
The JSQ policy was shown to minimize the average response time of jobs in finite systems comprising identical servers under the assumption of Poisson arrivals and exponential service times  in~\cite{Winston1977optimality}. This result was then generalised to service time distributions having non-decreasing hazard rates in~\cite{Weber1978}. Further extensions covered queues with state-dependent service rates, as demonstrated in~\cite{johri1989optimality}, and systems with finite buffers and general batch arrivals, as explored in~\cite{hordijk1990optimality}.

Recent works, such as~\cite{Mukherjee2016} and~\cite{Eschenfeldt2018}, analyzed the fluid and diffusion limits of the JSQ scheme. In the context of the fluid limit,~\cite{Mukherjee2016} showed that under a broader class of policies including the JSQ scheme, the proportion of servers with two or more jobs tends to zero as $n\to \infty$.
In the Halfin-Whitt regime, characterised by the normalised arrival rate $\lambda^{(n)}$ varying with system size $n$ as $\lambda^{(n)}=1-\beta/\sqrt{n}$ for some $\beta > 0$, ~\cite{Eschenfeldt2018} demonstrated that the diffusion-scaled process converges to a two-dimensional reflected OU process as $n \to \infty$. \cite{Braverman2020} showed that the many-server limit and the limit $t\to\infty$ can be interchanged and the stationary distribution of the limiting refected OU process can be used to approximate the stationary behaviour of the load balancing policy in the homogeneous setting. To do so, he showed that the steady-state proportion of idle servers as well as servers with precisely two jobs scale as $O(\sqrt{n})$, and the proportion of servers with more than three jobs scales as $O(1)$. ~\cite{gupta2019load} analyzed the JSQ scheme in the nondegenerate slowdown (NDS) regime introduced in~\cite{atar2012diffusion} where $\lambda^{(n)}=1-\beta/n$ for some fixed $\beta>0$ and established a diffusion limit for the total customer count process.

Note that the JSQ scheme faces challenges related to communication overhead. In addressing this issue, alternative scheduling schemes, like Power-of-d-Choices (Pod), were suggested \cite{Vvedenskaya1996queueing, Mitzenmacherthesis} where, as a new job arrives, $d\geq2$ servers are chosen uniformly at random from
the set of all servers, and the job is routed to the server with the shortest queue length among those $d$ servers. The problem of communication overhead is still there in the Pod scheme if $d$ is large. 
This can be eliminated with schemes relying only on the knowledge of idle servers in the system. In the Join-the-Idle-Queue (JIQ) scheme, first introduced by (\cite{Lu2011}), the dispatcher assigns the incoming arrival to an idle server (if available) chosen uniformly at random. \cite{Mukherjee2016} shows an equivalence between different load-balancing schemes in the diffusion regime. 

Relatively few works have explored load balancing in heterogeneous systems despite their importance in practical system design. The Pod scheme for heterogeneous systems has been studied in several papers, including \cite{Zhou2017, arpan_tcns, Makowski_SQd_2014}. While \cite{Makowski_SQd_2014} and \cite{shroff_heavy_traffic_2017} investigated the performance of the Pod policy under light and heavy traffic conditions for finite system sizes, respectively, \cite{arpan_tcns} explored its performance in the mean-field regime. It was demonstrated that the Pod scheme exhibits a reduced stability region in heterogeneous systems due to the infrequent sampling of faster servers. Subsequent works [see e.g.]~\cite{gardiner_perf, arpan_ssy} explored variations of the Pod scheme with the aim of enhancing its performance in heterogeneous systems while preserving the maximal stability region. ~\cite{Stolyar2015} analysed the JIQ scheme in a heterogeneous setting, revealing that the average waiting time of jobs under the JIQ scheme tends to zero in the fluid limit. 

The many-server systems with heterogeneous servers have received considerable attention in the $M/M/n$ system with single global queue setting. \cite{armony2005dynamic} showed that in the Halfin-Whitt regime the Fastest-Server-First (FSF) policy in which each job is assigned to fastest available server is asymptotically optimal in terms of stochastically minimising the diffusion-scaled total number of jobs in the system. {\cite{tezcan2008optimal} considered a similar model with heterogeneous servers where instead of a single central queue, each server pool consisting of servers of the same processing speed has a dedicated queue and an arriving job must be routed irrevocably to one of the server pools at its arrival instant. He showed that the Minimum Expected Delay-FSF (MED-FSF) routing policy in which a job is either  assigned to a queue with the minimum expected delay or to an available server having the highest speed asymptotically minimises the stationary distribution of the diffusion-scaled total number of jobs in the system in the Halfin-Whitt regime. In both works mentioned above, the diffusion-scaled total number of jobs in the system converges to the same one-dimensional diffusion process, unlike in our system, where the limit is a two-dimensional diffusion process.} In a more recent line of work, \cite{buke2019many} focused on multi-server queuing systems characterised by heterogeneous exponential servers and renewal arrivals where service rates are random and varies for each individual server and developed a framework has been to analyse the heavy traffic limit of these queues in a random environment, employing probability measure-valued stochastic processes.

\subsection{Notation}

Throughout the paper we use the following general notations. We denote  the set of non-negative integers as $\ZZ_+$ and  the set $\{1,\ldots, M\}$ by $[M]$. We use $\bar{\RR}$ to denote the set of extended reals, i.e., $\RR\cup \{\infty\}$, and define the sets $$\MM_{M,\rho} = \cbrac{(b_{j,i}\in \RR, i\in \ZZ_+, j\in [M]): \sum_{i\in \ZZ_+}\sum_{j\in [M]}\rho^{-i}|b_{j,i}|<\infty},$$ and $$\MM_{M,\rho}^1=\cbrac{(b_{j,i})\in \MM_{M,\rho}: \mbox{ there exists a $b_1\in \RR$ such that }b_{j,1}=b_1 \mbox{ for all }j\in [M]}.$$
We endow $\bar{\RR}$ with the order topology, where the neighborhoods of $\{\infty\}$ are the sets that contain a set $\{x>a\}$ for some $a\in\RR$. 
We denote all multi-dimensional terms, e.g., vectors, matrices and elements of $\MM_{M,\rho}$, with boldface letters.  
We also define $\DD_{\SSS}[0,\infty)$ as the set of all cadlag functions $f:[0,\infty)\to \SSS$. We denote the uniform norm $\norm{\mf x}_t=\sup_{0\leq s\leq t} |\mf x(s)|$ for $\mf x \in \DD_{\SSS}[0,\infty)$, where the $|\cdot|$ is considered to be the norm defined on space $\SSS$. For any $x\in \RR$, we define $(x)_+=\max\{x,0\}$. We assume that all random elements used in this work are defined on the probability space $\{\tilde{\Omega}, \cF, \PP\}$ and $\EE[\cdot]$ denotes the expectation with  respect to the probability measure $\PP$. We use upper-case letters to denote random elements in this probability space and corresponding lower-case  letters to denote their realisations. For any stochastic process $X=(X(t), t\geq 0)$, $X(\infty)$ denotes a random variable that follows the stationary distribution of the corresponding process. 
We use $\Rightarrow$ to denote weak convergence and $\overset{p}{\to}$ to denote convergence in probability. For any two non-negative real-valued functions $f$ and $g$ both depending on $n$, we write $f(n)=O(g(n))$ when $\limsup_{n \to \infty} \frac{f(n)}{g(n)}<\infty$ and $f(n)=o(g(n))$ when $\lim_{n \to \infty} \frac{f(n)}{g(n)} =0$.

\section{System Model} 
\label{Sec:system_model}
We consider a sequence of queuing systems, where the $n^{\textrm{th}}$ system consists of $n$ parallel servers, each with its own queue of infinite buffer size. Jobs are assumed to arrive at the $n$th system according to a Poisson process with rate $n\lambda^{(n)}$ where $\lambda^{(n)} < 1$ for all $n$. In this work, we consider the Halfin-Whitt regime, where 
 \begin{align}
 \tag{A1}\label{eqn:heavy_traffic_limit}
\lim_{n \to \infty}\sqrt{n}(1-\lambda^{(n)})= \beta,
\end{align}
for some constant $\beta >0$. Upon arrival, each job is assigned to a server where it either
receives service  immediately (if the server is idle at that instant) or waits in the
corresponding queue to be served later according to the First-Come-First-Served (FCFS) scheduling discipline. The servers are assumed to be {\em heterogeneous}, i.e., there are $M$ different server types and the service time of a job depends on the type of the server it is being served. In the $n$th system, there are $N_j^{(n)}$ servers of type $j\in [M]$ and the service time of a job being served by a server of type $j$ is exponentially distributed with rate $\mu_j$. Hence, 
\[\sum_{j=1}^M N_j^{(n)} = n, \mbox{ for all }n\in \ZZ_+,\]
and we assume that there exists $\gamma_j>0$ for each $j \in [M]$ such that 
\begin{align}
\tag{A2}\label{eq:proportion_assumption}
    \lim_{n \to \infty}\sqrt{n}\left|\frac{N_j^{(n)}}{n}-\gamma_j\right|= 0.
\end{align}
Furthermore, without loss of generality we assume that the normalised system capacity is unity, i.e., $\sum_{j=1}^M\mu_j\gamma_j=1$ and the server speeds satisfy $\mu_1 > \mu_2 > \ldots > \mu_M$. 


The classical JSQ policy is known to yield an asymptotically optimal performance for systems with homogeneous servers. Our main interest is to analyse
a modified version of the JSQ policy which we refer as the {\em Speed-Aware JSQ} {\em (SA-JSQ)} policy, and it is defined as follows:

\begin{definition}
Under the SA-JSQ policy, jobs are sent to a server with the minimum queue length 
among all the servers in the system upon arrival. Ties between servers of different types are broken by choosing the server type with the maximum speed and ties between servers of the same type are broken uniformly at random.   
\end{definition}

A load balancing scheme or policy $\Pi$ is called admissible if it considers only the present state of the system, i.e., the state and the queue sizes of each server, and the server speeds to dispatch jobs arriving at time $t\geq 0$. Let $\mathbf{Q}^{(n,\Pi)}(t)=\Big(Q^{(n, \Pi)}_{j,i}(t), j\in[M], i\geq1\Big)$ denote the system occupancy state at time $t\geq 0$,
where $Q^{(n,\Pi)}_{j,i}(t)$ is the number of type $j$ servers with at least $i$ jobs in the $n$th system under policy $\Pi$. It is clear that $Q^{(n, \Pi)}_{j,0}(t)=N_j^{(n)}$ for all $t\geq0$ and $j\in[M]$. Moreover, it is easy to verify that the process $\mathbf{Q}^{(n, \Pi)}=(\mathbf{Q}^{(n,\Pi)}(t):t\in \RR_+)$ is  Markov and  takes values in the space ${S}^{(n)}$
defined as 
\begin{equation*}
    {S}^{(n)}\define \{\mathbf{s}\in \ZZ_+^{M\times\infty}: N_j^{(n)}\geq s_{j, 1} \geq s_{j,2}\geq ... \geq 0, \mbox{ for all } j\in[M]\}.
\end{equation*}
We also define the scaled process $\mathbf{Y}^{(n, \Pi)}=(\mathbf{Y}^{(n, \Pi)}(t):t\geq 0)$, where $$\mathbf{Y}^{(n, \Pi)}(t)=\left(Y_{j,i}^{(n, \Pi)}(t),j\in[M], i\geq1\right)$$ with
\begin{equation}\label{eq:scaled_system_def}
Y_{j,1}^{(n, \Pi)}(t)=\frac{Q_{j,1}^{(n,\Pi)}(t)-N_j^{(n)}}{\sqrt{n}}, \ \   Y_{j,i}^{(n,\Pi)}(t)=\frac{Q_{j,i}^{(n,\Pi)}(t)}{\sqrt{n}}, \ i\geq2,j\in[M].
\end{equation}
Moreover, we also define a process $Y^{(n, \Pi)}_{[j_1,j_2],i}=(Y^{(n,\Pi)}_{[j_1,j_2],i}(t):t\geq 0)$ where
\begin{equation*}
Y^{(n, \Pi)}_{[j_1,j_2],i}(t)=\sum_{j=j_1}^{j_2}Y_{j,i}^{(n, \Pi)}(t).
\end{equation*}

Our main focus in this work is the analysis of load balancing systems under the SA-JSQ policy. To simplify the notation, we omit the superscript $\Pi$ when we consider the state of the system under the SA-JSQ policy.

\section{Main Results}


In this section, we state our main results and their consequences. We also provide intuitions for these results. The subsequent sections are devoted to the detailed proofs of these results.



We start with presenting the state-space collapse result in Proposition~\ref{prop:state_space_collapse} needed to prove the convergence of the stochastic processes. The proposition establishes that the diffusion-scaled number of idle servers in the fastest $M-1$ pools remains close to zero for all large $n$ and for all finite times. In other words, to keep track of the diffusion-scaled number of idle servers in the entire system in the limit as $n \to \infty$, it is sufficient to only keep track of the diffusion-scaled number of idle servers in the slowest pool (pool $M$). 

\begin{proposition}
\label{prop:state_space_collapse}
Suppose that~\eqref{eqn:heavy_traffic_limit} and~\eqref{eq:proportion_assumption} hold and $Y^{(n)}_{[1,M-1],1}(0) \Rightarrow Y_{[1,M-1],1}(0)$, where $Y_{[1,M-1],1}(0)$ is a proper random variable Then for any $\epsilon>0$, $t_0>0$ and $T\geq t_0$ we have
\begin{equation}
\label{eqn:ssc_1}
\mathbb{P}\brac{\sup_{s\in[t_0,T]} |Y^{(n)}_{[1,M-1],1}(s)|>\epsilon} \to 0, \ as \  n\to\infty.
\end{equation}
Moreover, if $Y^{(n)}_{[1,M-1],1}(0)\Rightarrow 0$, then \eqref{eqn:ssc_1} also holds for $t_0=0$.
\end{proposition}



\begin{figure}[h!]
  \centering
  \includegraphics[width=6cm]{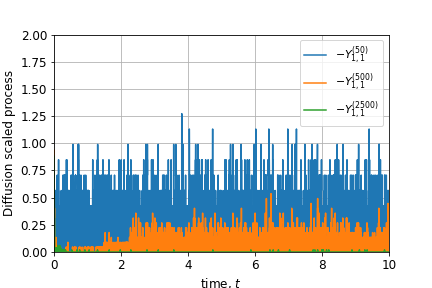}
\caption{Sample paths for $Y_{1,1}^{(n)}$ for varying $n$ with $M=2$, $\gamma_1=1/5$, $\gamma_2=4/5$, $\mu_1=20/8$, $\mu_2=5/8$ and $\beta=2$.}
\label{Fig:1}
\end{figure}
{Proposition~\ref{prop:state_space_collapse} states that, in the limit, almost all servers in the first $M-1$ pools should immediately become busy, even when there are idle servers in these pools at time $t=0$. This result can be intuitively explained by the fact that the total normalised capacity available in the fastest $M-1$ pools is strictly less than $1$ whereas the normalised arrival rate $\lambda^{(n)}$ into the system tends to $1$ as $n\to \infty$. Consequently, in the limit as $n \to \infty$, all servers in the fastest $M-1$ pools are used up by the SA-JSQ scheme as it prefers faster idle servers over slower ones.}
This can also be seen numerically in Figure~\ref{Fig:1}, where we plot the diffusion-scaled idle servers in the fastest pool as a function of time for different values of $n$ when $M=2$. It is apparent from Figure~\ref{Fig:1} that the fluctuations of the component $Y_{1,1}^{(n)}$ decreases to zero as $n$ increases. 

We now establish the limiting transient behaviour of the process $\mf Y^{(n)}$ as $n\to \infty$. To state the result, we assume that $Y_{j,1}^{(n)}(0) \Rightarrow 0$ for all $j\in[M-1]$, motivated by Proposition~\ref{prop:state_space_collapse}.

\begin{figure}
\centering
\subfigure[%
Idle Servers %
\label{Fig:diff_schemes}]{%
\includegraphics[width=6cm]{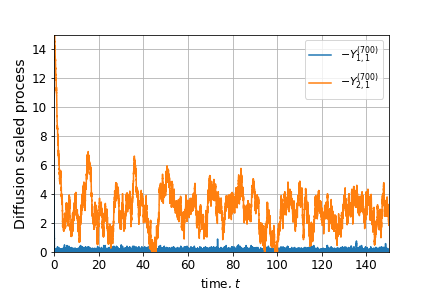}
}
\hfill
\centering
\subfigure[%
Queues with at least two jobs%
]{%
\includegraphics[width=6cm]{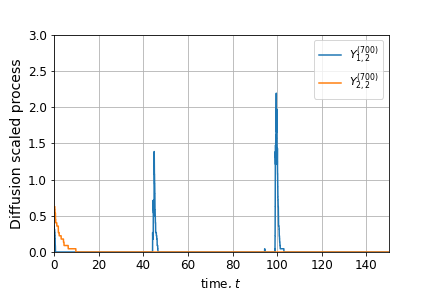}
}
\hfill
\centering
\subfigure[%
Queues with at least three jobs%
]{%
\includegraphics[width=6cm]{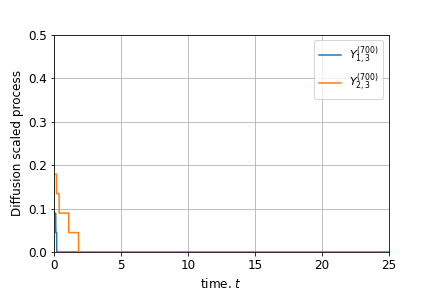}
}
\caption{Sample paths for $Y_{j,i}^{(n)}$ with $n=700$, $M=2$, $\gamma_1=1/5$, $\gamma_2=4/5$, $\mu_1=20/8$, $\mu_2=5/8$ and $\beta=2$.}
\label{Fig:comp}
\end{figure}

\begin{theorem}
\label{thm:SA-JSQ_diffusion_limit}
Suppose that~\eqref{eqn:heavy_traffic_limit} and~\eqref{eq:proportion_assumption} hold, and that there exists a random vector $\mf Y(0)=\Big( Y_{j,i}(0), i\geq 1, j\in[M]\Big) \in \MM_{M,\rho}$ for some $\rho>1$ with $Y_{j,1}(0)=0$ for all $j\in[M-1]$, such that $\mf Y^{(n)}(0) \Rightarrow \mf Y(0)$ as $n\to \infty$.
Then, we have
$\mf Y^{(n)} \Rightarrow \mf Y$ as $n\to \infty$, where $\mf Y =\brac{Y_{j,i},i\geq1,j\in[M]}$ is the unique solution in $\DD_{\MM_{M,\rho}}[0,\infty)$ of the following  stochastic integral equations:
\begin{align}
Y_{j,1}(t)&=0, \ j\in \cbrac{1,\dots,M-1},\\
\nonumber Y_{M,1}(t)&=Y_{M,1}(0) + \sqrt{2}W(t)-\beta t-\int_0^t (\mu_MY_{M,1}(s)-\mu_1Y_{1,2}(s))ds\\
&\qquad\qquad\qquad\qquad\qquad\qquad\qquad + \sum_{j=2}^M \int_0^t\mu_j {Y}_{j,2}(s)ds -U_1(t), \label{eq:dif_Y_1M} \\
\label{eq:dif_Y_12}
Y_{1,2}(t)&=Y_{1,2}(0)-\mu_1\int^t_0 (Y_{1,2}(s)-Y_{1,3}(s))ds+U_1(t),\\
Y_{j,2}(t)&=Y_{j,2}(0)-\mu_j\int^t_0 (Y_{j,2}(s)-Y_{j,3}(s))ds, \ j\in \cbrac{2,\dots,M},\label{eqn:Y_j2}\\
Y_{j,i}(t)&=Y_{j,i}(0)-\mu_j\int^t_0 (Y_{j,i}(s)-Y_{j,i+1}(s))ds, \ i\geq3 \ j\in [M]\label{eqn:Y_ji},
\end{align}
where $W$ is a standard Brownian motion and $U_1$ is the unique non-decreasing and non-negative process satisfying
\begin{equation}
 \label{eqn:diff_mjsq_fluct_term}
 \int_0^{\infty} \indic{Y_{M,1}(t) <0}dU_1(t) =0.
\end{equation}
 \end{theorem}
 
Equations \eqref{eqn:Y_j2} and \eqref{eqn:Y_ji} imply that the components $Y_{j,2}$ for $j\in \cbrac{2,\dots,M}$ and $Y_{j,i}$ for $i\geq3$, $j\in[M]$ decreases to zero exponentially fast as $t\to \infty$. More precisely, the explicit solution for these components can be written as
\begin{equation}
Y_{j,i}(t) = e^{-t} \cbrac{Y_{j,i}(0) + \sum_{k=1}^{\infty}\frac{1}{k!} t^k Y_{j,i+k}(0)}.
\end{equation}
Hence, if $Y_{j,2}(0)=0$ for $j\in \cbrac{2,\dots,M}$ and $Y_{j,i}(0)=0$ for  $j\in[M]$ and $i\geq3$, then $Y_{M,1}$ and $Y_{1,2}$ are the only non-zero components of the limiting diffusion process $\mf Y$, and equations \eqref{eq:dif_Y_1M} and \eqref{eq:dif_Y_12} reduce to
\begin{align}\label{eq:dif_process_0_1}
 Y_{M,1}(t)&=Y_{M,1}(0) + \sqrt{2}W(t)-\beta t+\int_0^t (\mu_1Y_{1,2}(s)-\mu_MY_{M,1}(s))ds-U_1(t),\\
\label{eq:dif_process_0_2}Y_{1,2}(t)&=Y_{1,2}(0)-\mu_1\int^t_0 Y_{1,2}(s)ds+U_1(t).   
\end{align}

{Equations \eqref{eq:dif_process_0_1} and \eqref{eq:dif_process_0_2} define a two-dimensional reflected Ornstein-Uhlenbeck (OU) process where the regulator process $U_1$ ensures that $Y_{M,1}$ remains non-positive and $Y_{1,2}$ remains non-negative at all times. Furthermore, we note that $Y_{1,2}$ can only increase when $Y_{M,1}$ is zero - a feature of any JSQ-type policy that ensures that jobs are sent to queues having an existing job only when there are no idle servers in the system. Furthermore,  both $Y_{M,1}$ and $Y_{1,2}$ can be non-zero at the same time. This is because even in the limit it can happen that jobs wait in queues belonging to the fastest pool when there are idle servers in the slowest server pool.
This is a fundamental drawback of having a separate queue for each server whose speed does not scale with $n$. In contrast, for systems where each queue is served by $O(n)$ servers of constant speed, in the limit as $n\to \infty$, no job has to wait in a queue while there are idle servers available. This yields a one dimensional diffusion limit as in~\cite{armony2005dynamic,tezcan2008optimal}.}

{Furthermore, the presence of the term $\mu_1Y_{1,2}(s)$ in~\eqref{eq:dif_process_0_1} can be intuitively explained as follows. 
In the Halfin-Whitt limit each departure from the fastest $M-1$ pools is immediately compensated by an arrival since in the limit the arrival rate of jobs is strictly above the rate at which the servers in the fastest $M-1$ pool can process jobs. 
This means that the moment a job leaves a queue with two jobs in the fastest pool a new job arrives which must be routed to an idle queue of the slowest pool since no idle queue exists in the fastest $M-1$ pools in the diffusion limit.}

The stability of a load balancing system operating under the SA-JSQ scheme for all $\lambda^{(n)} < 1$ was established in~\cite{bhambay2022asymptotic}. In our subsequent results, we investigate the steady-state characteristics of  the process $\mf Y^{(n)}$. We first establish the tightness of the sequence $(\mf Y^{(n)}(\infty))_{n\in \ZZ_+}$ in Theorem~\ref{thm:stationary_convergence}.

\begin{theorem}
\label{thm:stationary_convergence}
Suppose that~\eqref{eqn:heavy_traffic_limit} and~\eqref{eq:proportion_assumption} hold. Then, there exist positive constants $C_1, C_2$ and $C_3$ such that for all sufficiently large $n$ the following inequalities hold
\begin{align}
\label{eq:Expectation_1}&\EE[|Y_{j,1}^{(n)}(\infty)|] <C_1n^{-1/2}\mbox{ for all }1\leq j\leq M-1,\\
\label{eq:Expectation_2}&\EE[|Y_{M,1}^{(n)}(\infty)|]\leq C_2,\\
\label{eq:Expectation_3}&\EE[|Y_{1,2}^{(n)}(\infty)|]\leq C_3.
\end{align}
{Furthermore, for any  $a\in (1/2,1)$ there exists a constant $C_4(a)$ such that 
for all sufficiently large $n$ we have
\begin{equation}
\label{eq:Expectation_4}\EE[|Y_{j,i}^{(n)}(\infty)|]\leq C_4(a)n^{-(a-1/2)} \mbox{ for all }2\leq j\leq M, i=2 \mbox{ and } \mbox{ for all }1\leq j\leq M, i\geq 3.
\end{equation}}
\end{theorem}

{Theorem~\ref{thm:stationary_convergence} implies that the expected number of idle servers in the fastest $M-1$ pools are bounded above by constants independent of $n$. The expected number of servers with at least two jobs in pool $1$ and the expected number of idle servers in pool $M$ are $O(\sqrt{n})$. Moreover, the expected number of servers with at least three jobs in all pools and the expected number of servers with at least two jobs in the slowest $M-1$ pools all approaches to 0 faster than $n^{-(a-1/2)}$ for any $1/2<a<1$. As mentioned before, these results can be intuitively explained by the facts that (i) the spare capacity $n(1-\lambda^{(n)})$ in the Halfin-Whitt regime scales as $O(\sqrt{n})$ and (ii) this spare capacity  is mostly available through the idle servers in the slowest pool since the SA-JSQ scheme prefers faster servers over slower ones whenever they have the same queue lengths. By a simple application of the Markov inequality, it follows from the above theorem that the sequence $(\mf Y^{(n)}(\infty))_{n\in \ZZ_+}$ is tight.

\begin{remark}
From Theorem~\ref{thm:stationary_convergence} we can compute a bound on the probability with which a job has to wait in a queue before being processed in the steady-state. This probability is given by $\PP(\sum_{j\in [M]}Q_{j,1}^{(n)}(\infty)=n)\leq \PP(Q_{M,1}^{(n)}(\infty)=N_M^{(n)})\leq\EE[Q_{M,1}^{(n)}]/N_M^{(n)}=O(1/\sqrt{n})$ where the last equality follows by~\eqref{eq:Expectation_2} and~\eqref{eq:proportion_assumption}.
Thus, the stationary waiting probability of a job goes to zero at rate $O(1/\sqrt{n})$.
\end{remark}}

Next, we establish the positive recurrence of the limiting diffusion process $\mf Y$. This is the final ingredient needed to establish the interchange of limits as stated in Corollary~\ref{cor:stationary_convergence}.

\begin{theorem}
\label{thm:stationarity_diffusion}
The diffusion process $\mf Y$ defined in Theorem~\ref{thm:SA-JSQ_diffusion_limit}
is positive recurrent.
 \end{theorem}

\begin{corollary}\label{cor:stationary_convergence}
Let $\pi^{(n)}$  be the stationary distribution of the process $\mf Y^{(n)}$ and let $\pi$ denote the stationary distribution of $\mf Y$. Then, $\pi^{(n)}\Rightarrow \pi$.  
\end{corollary}

Hence, by Corollary~\ref{cor:stationary_convergence}, we can conclude that the stationary distribution of the diffusion-scaled queue length process $\mf Y^{(n)}$ can be well approximated by that of the limiting diffusion process $\mf Y$ for large $n$.


{In our final result, 
we focus on the asymptotic optimality of the SA-JSQ scheme in the Halfin-Whitt regime.
Note that optimality can be defined in various ways. In the present paper, we are interested in asymptotic optimality which is is defined in terms of minimising the steady-state distribution of the diffusion-scaled total number of jobs and the diffusion-scaled number of waiting jobs in the system as  in~\cite{armony2005dynamic,tezcan2008optimal}. To define a similar notion of optimality we define
\begin{equation}
Y^{(n)}_{+1}(t)=\frac{\sum_{j\in [M]}\sum_{i \geq 1}Q_{j,i}^{(n)}(t)-n}{\sqrt{n}}, \ \   \text{ and } Y^{(n)}_{+2}(t)=\sum_{j\in [M]}\sum_{i \geq 1}\frac{Q_{j,i}^{(n)}(t)}{\sqrt{n}}
\end{equation}
to be respectively the diffusion-scaled total number of customers in the system and  the diffusion-scaled number of waiting jobs in the system at time $t\in [0,\infty]$ under the SA-JSQ policy. The same quantities under a stable admissible policy $\Pi$ are denoted by $Y^{(n, \Pi)}_{+1}(t)$ and $Y^{(n, \Pi)}_{+2}(t)$, respectively. We are now ready to state our asymptotic optimality result.}

{\begin{theorem}
\label{thm:sajsq_optimality}
Assume~\eqref{eqn:heavy_traffic_limit} and~\eqref{eq:proportion_assumption} hold. Then, for any stable admissible policy $\Pi$ we have 
\begin{align}
&\lim_{n\to \infty}\mathbb{P}\brac{Y^{(n)}_{+1}(\infty)>y}\leq \liminf_{n\to \infty}\mathbb{P}\brac{Y^{(n,\Pi)}_{+1}(\infty)>y} \label{eqn:asmp_1}\\ 
&\lim_{n\to \infty}\mathbb{P}\brac{Y^{(n)}_{+2}(\infty)>y}\leq \liminf_{n\to \infty}\mathbb{P}\brac{Y^{(n,\Pi)}_{+2}(\infty)>y}\label{eqn:asmp_2},
\end{align}
for all $y\in \RR$.
\end{theorem}}

{Theorem~\ref{thm:sajsq_optimality} ensures that the SA-JSQ scheme asymptotically minimises the steady-state distribution of the diffusion-scaled total number of jobs and the diffusion-scaled number of waiting jobs in the system in the Halfin-Whitt regime. }

The remainder of the paper is organised as follows. The proof of the state space collapse result, stated in Proposition~\ref{prop:state_space_collapse}, is provided in Section~\ref{sec:ssc}.  Section~\ref{sec:diffusion_sajsq} is dedicated to the proof of Theorem~\ref{thm:SA-JSQ_diffusion_limit}.
Theorem~\ref{thm:stationary_convergence} is proved in Section~\ref{sec:stationary_convergence}. One essential component of the proof in Section~\ref{sec:stationary_convergence} is finding a suitable Lyapunov function. 
We do this in Section~\ref{sec:lyapunov_func}.
In Section~\ref{sec:mod_sys}, we prove the asymptotic optimality of the SA-JSQ scheme, as stated in Theorem~\ref{thm:sajsq_optimality}. Finally, we conclude the paper in Section~\ref{sec:conclusion}.

\section{State Space Collapse: Proof of Proposition~\ref{prop:state_space_collapse}}
\label{sec:ssc}

{In this section, we establish the state space collapse result in load-balancing systems under the SA-JSQ scheme. To prove Proposition~\ref{prop:state_space_collapse}, we first express the evolution of $Y^{(n)}_{[1,M-1],1}$ in terms of arrival and departure processes. We then examine the system's behavior in two parts: First, from $t=0$ up to the point where $Y^{(n)}_{[1,M-1],1}$ reaches zero (defined as $\theta_0^{(n)}$ below), and second, the behavior after hitting zero. We show that once $Y^{(n)}_{[1,M-1],1}$ hits 0, it stays very close to 0, i.e., the probability that it will deviate by an amount $\epsilon$ approaches 0 as $n \to \infty$. Finally, we show that the time $\theta_0^{(n)}$ to hit zero, also converges to 0 in probability. The following lemma from~\cite{Pang2007} is key in our proofs of Proposition~\ref{prop:state_space_collapse} and Theorem~\ref{thm:SA-JSQ_diffusion_limit}.}

{\begin{lemma}[Lemma~5.8 in ~\cite{Pang2007}]
\label{lem:stochastic_bounded}
Suppose that, for
each $n\geq1$, $M_n =(M_n(t) : t \geq0)$ is a square-integrable martingale (with respect
to a specified filtration) with predictable quadratic variation $\langle M_n\rangle = (\langle M_n\rangle(t) :
t \geq0)$, i.e., $M_n^2-\langle M_n\rangle=\brac{M_n^2(t)-\langle M_n\rangle(t):t\geq0}$ is a martingale by
the Doob-Meyer decomposition. If the sequence of random variables $\cbrac{\langle M_n\rangle (T)}_{n\geq1}$ is tight for each $T \geq 0$, then the sequence of
stochastic processes $\cbrac{M_n:n\geq1}$ is stochastically bounded in $\DD_\RR[0,\infty)$ i.e., the sequence $\cbrac{\sup_{0\leq t \leq T} {|M_n(t)|}}_{n\geq 1}$ is tight for every $T\geq0$.
\end{lemma}}

{\em Proof of Proposition~\ref{prop:state_space_collapse}.}
We first define
\begin{equation*}
\theta_0^{(n)}=\inf\cbrac{t\geq0: Y^{(n)}_{[1,M-1],1}(t)=0}.
\end{equation*}
Now, we can write~\eqref{eqn:ssc_1} as
\begin{align}
\nonumber \mathbb{P}\brac{\sup_{s\in[t_0,T]} \vert Y^{(n)}_{[1,M-1],1}(s)\vert >\epsilon}&=\mathbb{P}\brac{\sup_{s\in[t_0,T]} -Y^{(n)}_{[1,M-1],1}(s)>\epsilon,\theta_0^{(n)} \leq t_0}\\&\qquad\qquad\qquad+\mathbb{P}\brac{\sup_{s\in[t_0,T]}-Y^{(n)}_{[1,M-1],1}(s)>\epsilon,\theta_0^{(n)} > t_0} \nonumber\\
&\leq \mathbb{P}\brac{\sup_{s\in[t_0,T]}-Y^{(n)}_{[1,M-1],1}(s)>\epsilon,\theta_0^{(n)} \leq t_0} + \mathbb{P}\brac{\theta_0^{(n)} > t_0} \nonumber\\
&\leq \mathbb{P}\brac{\sup_{s\in[\theta_0^{(n)},T]} -Y^{(n)}_{[1,M-1],1}(s)>\epsilon} + \mathbb{P}\brac{\theta_0^{(n)} > t_0}\label{eqn:ssc_2}.
\end{align}
To prove~\eqref{eqn:ssc_1}, it is sufficient to show that each term on the right hand side of~\eqref{eqn:ssc_2} converges to $0$ as $n\to \infty$. 
For a given state $\mf q \in S^{(n)}$, the process $Y^{(n)}_{[1,M-1],1}$ decreases by $1/\sqrt{n}$ with rate $\sum_{j\in [M-1]}\mu_j (q_{j,1}-q_{j,2})$ and it {increases} by $1/\sqrt{n}$ with rate $n\lambda^{(n)} \indic{y_{[1,M-1],1} < 0}$. Let $A$ and $D$ be unit rate Poisson processes and define square integrable martingales
\begin{align}
\label{eqn:MA}
\hat{M}_A^{(n)}(t)&=\frac{A(n \lambda ^{(n)}t)-n\lambda^{(n)} t}{\sqrt{n}},\\
\hat{M}_D^{(n)}(t)&=\frac{1}{\sqrt{n}}\cbrac{D\brac{ \sum_{j\in[M-1]}\mu_jN_j^{(n)} t}-\sum_{j\in[M-1]}\mu_jN_j^{(n)} t}.
\end{align}
For any $0\leq t\leq \theta_0^{(n)}$, we can write
\begin{align}
-Y^{(n)}_{[1,M-1],1}(t)&=-Y^{(n)}_{[1,M-1],1}(0) - \frac{A\brac{n\lambda^{(n)}\int_0^t \indic{Y^{(n)}_{[1,M-1],1}(s)<0}ds}  }{\sqrt{n}}\nonumber \\
\nonumber&\quad \quad  +\frac{1}{\sqrt{n}} D\sbrac{\int_0^t \sum_{j\in[M-1]}\mu_j(Q_{j,1}^{(n)}(s) -  Q_{j,2}^{(n)}(s))ds} \\
&=-Y^{(n)}_{[1,M-1],1}(0) - \frac{A(n\lambda^{(n)}t)}{\sqrt{n}}\nonumber\\
&\quad+\frac{1}{\sqrt{n}} D\sbrac{\int_0^t \sum_{j\in[M-1]}\mu_j(Q_{j,1}^{(n)}(s) -  Q_{j,2}^{(n)}(s))ds} \nonumber\\
&\leq -Y^{(n)}_{[1,M-1],1}(0) - \frac{A(n\lambda^{(n)}t)}{\sqrt{n}} +\frac{1}{\sqrt{n}} D\brac{ \sum_{j\in[M-1]}\mu_jN_j^{(n)} t} \nonumber\\
&= -Y^{(n)}_{[1,M-1],1}(0) -\hat{M}_A^{(n)}(t) - \frac{ n\lambda^{n} t}{\sqrt{n}} + \hat{M}_D^{(n)}(t) + \frac{ \sum_{j\in[M-1]}\mu_jN_j^{(n)} t}{\sqrt{n}} \nonumber\\
&= -Y^{(n)}_{[1,M-1],1}(0) -\hat{M}_A^{(n)}(t) + \hat{M}_D^{(n)}(t) -c(t,n),
\label{eqn:ssc_3}
\end{align}
where 
\begin{equation}
\label{eq:c_n}
 c(t,n)=\sqrt{n}t\left(\lambda^{(n)}  - n^{-1}\sum_{j\in[M-1]}N_j^{(n)}\mu_j\right).   
\end{equation}

Here, the second equality follows as for any $0\leq s\leq \theta_0^{(n)}$, we have $\indic{Y^{(n)}_{[1,M-1],1}(s)<0}=1$ and the inequality follows from the fact that $Q_{j,1}^{(n)} \leq N_j^{(n)}$. Assumptions ~\eqref{eqn:heavy_traffic_limit} and~\eqref{eq:proportion_assumption} imply that $\lambda^{(n)}\to 1$ and $n^{-1}\sum_{j\in[M-1]}N_j^{(n)}\mu_j\to 1-\gamma_M\mu_M$, and hence, we have $c(t,n)=O(\sqrt{n})$ for any $t\geq0$. For large enough $n$, we observe that for any $t_0>0$ we can write 
\begin{align}
 \mathbb{P}\brac{\theta_0^{(n)} > t_0}&\leq  \mathbb{P}\brac{\sup_{0\leq t \leq t_0}Y^{(n)}_{[1,M-1],1}(t)<0} \nonumber \\
 &\leq \mathbb{P} \brac{\sup_{0\leq t \leq t_0}Y^{(n)}_{[1,M-1],1}(0) + \hat{M}_A^{(n)}(t) - \hat{M}_D^{(n)}(t) + c(t,n)< 0} \nonumber\\
  &\leq \mathbb{P} \brac{Y^{(n)}_{[1,M-1],1}(0) + \hat{M}_A^{(n)}(t_0) - \hat{M}_D^{(n)}(t_0) <-c(t_0,n)} \nonumber\\
 &\leq \mathbb{P} \brac{|Y^{(n)}_{[1,M-1],1}(0)|\geq c(t_0,n)/3} + \mathbb{P} \brac{|\hat{M}_A^{(n)}(t_0)|\geq c(t_0,n)/3} \nonumber\\
 & \quad \quad \quad \quad + \mathbb{P} \brac{|\hat{M}_D^{(n)}(t_0)|\geq c(t_0,n)/3},\label{eqn:ssc_4}
 \end{align}
%
where the second inequality follows from~\eqref{eqn:ssc_3}, and the last inequality follows from the union bound. 

The optional quadratic variations $\langle \hat{M}_A^{(n)}\rangle(t)=\lambda^{(n)}t$ and $$\langle \hat{M}_D^{(n)}\rangle(t)=\sum_{j \in [M-1]}\mu_j N_j^{(n)}t/n$$ are stochastically bounded due to Assumptions~\eqref{eqn:heavy_traffic_limit} and \eqref{eq:proportion_assumption}. Lemma~\ref{lem:stochastic_bounded} implies that the martingales $\hat{M}_A^{(n)}$ and $\hat{M}_D^{(n)}$ are also stochastically bounded. As $n\to \infty$, the sequence $\cbrac{Y^{(n)}_{[1,M-1],1}(0)}_{n\geq1}$ is  weakly convergent and $c(t_0,n)\to \infty$, each term on the right-hand side of~\eqref{eqn:ssc_4} converges to $0$. 

Next, we define
\begin{align*}
\theta_{\epsilon}^{n,1}&=\inf\cbrac{s>\theta_0^{(n)}:-Y^{(n)}_{[1,M-1],1}(s)>\epsilon},\\
\theta_{\epsilon}^{n,2}&=\sup\cbrac{\theta_0^{(n)}\leq s <\theta_{\epsilon}^{n,1}:-Y^{(n)}_{[1,M-1],1}(s)<\epsilon/2 }.
\end{align*}
Now, using a similar reasoning as in \eqref{eqn:ssc_3}, we can write 
\begin{align}
\nonumber&\mathbb{P}\brac{\sup_{s\in[\theta_0^{(n)},T]} -Y^{(n)}_{[1,M-1],1}(s)>\epsilon}\\&
\qquad\qquad\leq \mathbb{P}\brac{\theta_0^{(n)}<\theta_{\epsilon}^{n,2} <\theta_{\epsilon}^{n,1}<T} \nonumber\\  
&\qquad\qquad\leq \mathbb{P}\brac{\sup_{\theta_\epsilon^{n,2}\leq t\leq T}\frac{1}{\sqrt{n}} D\brac{ \sum_{j\in[M-1]}\mu_jN_j^{(n)} (t-\theta_\epsilon^{n,2})}- \frac{A(n\lambda^{(n)}(t-\theta_\epsilon^{n,2}))}{\sqrt{n}}>\epsilon/2} \nonumber\\&
\qquad\qquad\leq \mathbb{P}\left(\sup_{\theta_\epsilon^{n,2}\leq t\leq T}\hat{M}_D^{(n)}(t-\theta_\epsilon^{n,2})-\hat{M}_A^{(n)}(t-\theta_\epsilon^{n,2})
-c(t-\theta_\epsilon^{n,2},n)
>\epsilon/2\right),\label{eqn:ssc_5}
\end{align}
where the second inequality follows since the term in the supremum is the change in the process $-Y^{(n)}_{[1,M-1],1}$ from time $\theta_\epsilon^{n,2}$ to time $t$. Again using the fact that the martingales $\hat{M}_A^{n}$, $\hat{M}_D^{n}$ are stochastically bounded and $c(t-\theta_\epsilon^{n,2},n)$  diverges to $\infty$, the right-hand side converges to 0. When $Y^{(n)}_{[1,M-1],1}(0)\toP 0$, the second claim of the propositon can be proved using the same argument by writing
\begin{equation*}
\mathbb{P}\brac{\sup_{s\in[0,T]} -Y^{(n)}_{[1,M-1],1}(s)>\epsilon}\leq \mathbb{P}\brac{-Y^{(n)}_{[1,M-1],1}(0)>\epsilon/2} + \mathbb{P}\brac{0 <\theta_{\epsilon}^{n,1}<T},
\end{equation*}
and replacing $\theta_\epsilon^{n,2}$ with 0 in \eqref{eqn:ssc_5}. This completes the proof. 
\qed
\endproof

\section{Transient Analysis: Proof of Theorem~\ref{thm:SA-JSQ_diffusion_limit}}
\label{sec:diffusion_sajsq}

{The proof of Theorem~\ref{thm:SA-JSQ_diffusion_limit} relies on the continuous mapping theorem by showing that the system of integral equations in Proposition~\ref{thm:skorohod_mapping} defines a continuous mapping on $\DD_{\MM_{M,\rho}}[0,\infty)$. We first demonstrate that these equations possess a unique solution and hence correspond to a well-defined mapping. The proof of Proposition~\ref{thm:skorohod_mapping} can be found in Appendix~\ref{append:1}.
As the next step, we introduce a truncated system and formulate the evolution equations in terms of martingales specific to this truncated system. We establish the convergence of the martingale components appearing in the evolution equations for this truncated system, and finally, we prove that both the truncated and untruncated systems exhibit similar behavior in the limit.}

\begin{proposition}
\label{thm:skorohod_mapping} Consider the following set of integral equations
\begin{align}
y_1(t)&=b_1+z_1(t) -\int_0^t (\mu_M y_1(s)-\mu_1y_{1,2}(s))ds +\sum_{j=2}^M \mu_j \int_0^t y_{j,2}(s)ds-u_1(t),\label{eqn:map_1}\\
y_{1,2}(t)&=b_{1,2}+ z_{1,2}(t) -\int_0^t \mu_1(y_{1,2}(s)-y_{1,3}(s))ds +u_1(t) -u_2(t),\label{eqn:map_2}\\
y_{j,2}(t)&=b_{j,2} +z_{j,2}(t) -\int_0^t \mu_j(y_{j,2}(s)-y_{j,3}(s))ds, \ \forall j\in \{2,\dots,M\},\label{eqn:map_3}\\
y_{j,i}(t)&=b_{j,i} + z_{j,i}(t) -\int_0^t \mu_j (y_{j,i}(s)-y_{j,i+1}(s))ds, \ i\geq 3, \ j\in[M],\label{eqn:map_4}\\
y_1(t)&\leq 0, \ 0\leq y_{1,2}(t)\leq B, \ y_{j,i}(t)\geq0, \mbox{ for all }i\geq 2, j\in [M]\label{eqn:map_5}
\end{align}
where $u_1$ and $u_2$ are non-decreasing and non-negative functions such that 
\begin{align}
\int_0^{\infty} \indic{y_1(t)<0}du_1(t)&=0,\\
\int_0^{\infty} \indic{y_{1,2}(t)<B}du_2(t)&=0\label{eqn:map_6}.
\end{align}
Given $\mathbf{b}\in \MM_{M,\rho}^1$ and $\mathbf{z}\in \DD_{\MM_{M,\rho}^1}[0,\infty)$ the integral equations ~\eqref{eqn:map_1}-\eqref{eqn:map_3} has a unique solution $(\mathbf{y} ,\mathbf{u}) \in \DD_{\MM_{M,\rho}^1}[0,\infty) \times \DD_{\RR_+^2}[0,\infty)$. Moreover, there exists a well defined  function $(f,  g): \bar{\RR}_+ \times \MM_{M,\rho}^1 \times \DD_{\MM_{M,\rho}^1}[0,\infty) \to\DD_{\MM_{M,\rho}^1}[0,\infty) \times \DD_{\RR_+^2}[0,\infty)$ which maps $(B,\mf b,\mf z)$ in to $\mf y =f(B,\mf b,\mf z)$ and $\mf u=g(B,\mf b,\mf z)$. Furthermore, the function $(f,g)$ is continuous on $\bar{\RR}_+ \times \MM_{M,\rho}^1 \times \DD_{\MM_{M,\rho}^1}[0,\infty)$. Finally, $\mf y$ being continuous implies that $\mf x$ and $\mf u$ are also continuous. 
\end{proposition}


Next, we consider a truncated system where the queues with more than two jobs are not allowed to form in pool $1$ and no queue with more than one job will form in pools $\cbrac{2,\dots,M}$, even though there might be such queues initially present at time 0. We first analyse this system and then show that the truncated system and the original system are indistinguishable in the limit. For any $t\geq0$, let $\hat{\mathbf{Q}}^{(n)}(t)=(\hat{Q}^{(n)}_{j,i}(t),i\geq1,j\in[M])$, where $\hat{Q}^{(n)}_{j,i}(t)$ is the number of type $j$ servers with at least $i$ jobs in the truncated system. For all $t\geq0$, we can represent the dynamics of the truncated system as 
\begin{align}
\hat{Q}_{1,1}^{(n)}(t)&= \hat{Q}_{1,1}^{(n)}(0) + A(n\lambda^{(n)} t) - D_{1,1}\Big(\mu_1\int ^t _0 (\hat{Q}_{1,1}^{(n)}(s)-\hat{Q}_{1,2}^{(n)}(s))ds\Big )-\hat{U}_{1,1}^{(n)}(t),\label{eqn:truc_sys_1}\\
\hat{Q}_{j,1}^{(n)}(t)&= \hat{Q}_{j,1}^{(n)}(0) + \hat{U}_{j,1}^{(n)}(t)-D_{j,1}\Big(\mu_j\int ^t _0 (\hat{Q}_{j,1}^{(n)}(s)-\hat{Q}_{j,2}^{(n)}(s))ds\Big), \ 2\leq j\leq M,\label{eqn:truc_sys_2}\\
\hat{Q}_{1,2}^{(n)}(t)&= \hat{Q}_{1,2}^{(n)}(0) + \hat{U}_{1,2}^{(n)}(t)-D_{1,2}\Big(\mu_1\int ^t _0 (\hat{Q}_{1,2}^{(n)}(s)-\hat{Q}_{1,3}^{(n)}(s))ds\Big),\label{eqn:truc_sys_3}\\
\hat{Q}_{j,2}^{(n)}(t)&= \hat{Q}_{j,2}^{(n)}(0) -D_{j,2}\Big(\mu_j\int ^t _0 (\hat{Q}_{j,2}^{(n)}(s)-\hat{Q}_{j,3}^{(n)}(s))ds\Big), \ j\in \cbrac{2,\dots,M},\label{eqn:truc_sys_4}\\
\hat{Q}_{j,i}^{(n)}(t)&= \hat{Q}_{j,i}^{(n)}(0)-D_{j,i}\Big(\mu_j\int ^t _0 (\hat{Q}_{j,i}^{(n)}(s)-\hat{Q}_{j,i+1}^{(n)}(s))ds\Big), \ i\geq3, \ j\in [M],\label{eqn:truc_sys_5}
\end{align}
where $A$ and $D_{j,i}$ are mutually independent unit-rate Poisson processes and 
\begin{align}
\hat{U}_{1,1}^{(n)}(t)&=\int^t_0  \indic{\hat{Q}_{1,1}^{(n)}(s)=N_1^{(n)} }dA(n\lambda^{(n)}s),\label{eqn:u_11}\\
\hat{U}_{j,1}^{(n)}(t)&=\int^t_0  \indic{\hat{Q}_{l,1}^{(n)}(s)=N_l^n, \ \forall l\in[j-1], \ \hat{Q}_{j,1}^{(n)}(s)<N_j^{(n)}}dA(n\lambda^{(n)}s), \ 2\leq j \leq M, \label{eqn:u_1j}\\
\hat{U}_{1,2}^{(n)}(t)&=\int^t_0  \indic{ \hat{Q}_{j,1}^{(n)}(s)=N_j^{(n)}, \ \forall j\in [M],\hat{Q}_{1,2}^{(n)}(s)<N_1^{(n)}}dA(n\lambda^{(n)}s). \label{eqn:u_12}
\end{align}
It is easy to see that $\hat{U}_{1,1}^{(n)}(t)$ counts  the number of arrivals  which see all pool $1$ servers busy upon arrival within the time interval $[0,t]$. Similarly, $\hat{U}_{j,1}^{(n)}(t)$ for $2 \leq j\leq M$ counts  the number of arrivals which see all servers in pools in $[j-1]$ busy but finds at least one idle server in pool $j$ within the time interval $[0,t]$. Finally, $\hat{U}_{1,2}^{(n)}(t)$ counts the number of arrivals  which, upon arrival, see all servers in the system busy but find some servers in pool 1 with only one job within time interval $[0,t]$.  

Similar to~\eqref{eq:scaled_system_def}, we define the diffusion-scaled version of the process $\hat{\mathbf{Q}}^{(n)}$ as 
\begin{align*}
\hat{Y}_{j,1}^{(n)}(t)=\frac{\hat{Q}_{j,1}^{(n)}(t)-N_j^{(n)}}{\sqrt{n}}, \hat{Y}_{j,i}^{(n)}(t) =\frac{\hat{Q}_{j,i}^{(n)}(t)}{\sqrt{n}}, i\geq2, j\in [M],\mbox{ and }
\hat{Y}_{[1,M],1}^{(n)}(t)=\sum_{j\in [M]}\hat{Y}_{j,1}^{(n)}(t).
\end{align*}
In addition to $\hat{M}_{A}^{(n)}$ (defined in~\eqref{eqn:MA}) and $\hat{M}_{D}^{(n)}$, for any $i\geq 1$ and $j\in [M]$, we also define following martingales 
\begin{align}
\label{eqn:martingales_d_j_i}
\hat{M}_{D,j,i}^{(n)}(t)&=\frac{1}{\sqrt{n}}D_{j,i}\left(\mu_j\int ^t _0 (\hat{Q}_{j,i}^{(n)}(s)-\hat{Q}_{j,i+1}^{(n)}(s))ds\right)-\frac{1}{\sqrt{n}}\int ^t _0 \mu_j(\hat{Q}_{j,i}^{(n)}(s)-\hat{Q}_{j,i+1}^{(n)}(s))ds.
\end{align}
The martingale $\hat{M}_{D,j,i}^{(n)}(t)$ are square integrable martingale (c.f.~Lemma 3.2 in \cite{Pang2007}) with optional and predictable quadratic variations
\begin{align*}
    [\hat{M}_{D,j,i}^{(n)}](t)&=\frac{1}{n}D_{j,i}\left(\mu_j\int ^t _0 (\hat{Q}_{j,i}^{(n)}(s)-\hat{Q}_{j,i+1}^{(n)}(s))ds\right)\\
    \langle\hat{M}_{D,j,i}^{(n)}\rangle(t)&=\frac{1}{n}\mu_j\int ^t _0 (\hat{Q}_{j,i}^{(n)}(s)-\hat{Q}_{j,i+1}^{(n)}(s))ds.
\end{align*}

Combining~\eqref{eqn:truc_sys_1}, and~\eqref{eqn:truc_sys_2}, we have
%
\begin{multline*}
\hat{Y}_{[1,M],1}^{(n)}(t)=\hat{Y}_{[1,M],1}^{(n)}(0)+\hat{M}_{A}^{(n)}(t)-\sum_{j\in[M]}\hat{M}_{D,j,1}^{(n)}(t) - \sqrt{n}\Big(1-\lambda^{(n)}\Big) t 
\\
-\sum_{j\in[M]}\int_0^t\mu_j(\hat{Y}^{(n)}_{j,1}(s)-\hat{Y}^{(n)}_{j,2}(s))ds - \frac{\hat{U}_{1,1}^{(n)}(t)}{\sqrt{n}}+\frac{\sum_{j=2}^{M}\hat{U}_{j,1}^{(n)}(t)}{\sqrt{n}}.  
\end{multline*}
Slightly re-arranging the terms to obtain the form in \eqref{eqn:map_1}, we get 
\begin{align}
\nonumber\hat{Y}_{[1,M],1}^{(n)}(t)&=\hat{Y}_{[1,M],1}^{(n)}(0)+\hat{M}_{A}^{(n)}(t)-\sum_{j\in[M]}\hat{M}_{D,j,1}^{(n)}(t) - \sqrt{n}(1-\lambda^{(n)}) t - \sum_{j=1}^{M-1}\int_0^t\mu_j\hat{Y}^{(n)}_{j,1}(s)ds\\
\nonumber&\quad-\mu_M \int_0^t (\hat{Y}^{(n)}_{M,1}(s)-\hat{Y}^{(n)}_{[1,M],1}(s))ds - \int_0^t (\mu_M \hat{Y}^{(n)}_{[1,M],1}(s) -\mu_1 \hat{Y}^{(n)}_{1,2}(s))ds\\&\quad + \sum_{j=2}^M \int_0^t\mu_j \hat{Y}^{(n)}_{j,2}(s)ds-\hat{V}_1^{(n)}(t) , \label{eqn:diff_SA-JSQ_evolu_scaled_1_2}   
\end{align}
where,
\begin{align*}
    \hat{V}_1^{(n)}(t)&=\frac{\hat{U}_{1,1}^{(n)}(t)}{\sqrt{n}}-\frac{\sum_{j=2}^{M}\hat{U}_{j,1}^{(n)}(t)}{\sqrt{n}}.
\end{align*}
Similarly, using~\eqref{eqn:truc_sys_3}-\eqref{eqn:truc_sys_5}, we obtain
\begin{align}
\hat{Y}_{1,2}^{(n)}(t)&=\hat{Y}_{1,2}^{(n)}(0)-\hat{M}_{D,1,2}^{(n)}(t) -\int_0^t \mu_1(\hat{Y}^{(n)}_{1,2}(s)-\hat{Y}^{(n)}_{1,3}(s))ds
+\hat{V}_1^{(n)}(t)-\hat{V}_2^{(n)}(t)\label{eqn:diff_SA-JSQ_evolu_scaled_2_1},\\
\hat{Y}^{(n)}_{j,2}(t)&=\hat{Y}^{(n)}_{j,2}(0) -\hat{M}_{D,j,2}^{(n)}(t)-\int_0^t \mu_j(\hat{Y}^{(n)}_{j,2}(s)-\hat{Y}^{(n)}_{j,3}(s))ds,\ j \in \cbrac{2,\dots,M},\label{eqn:diff_SA-JSQ_evolu_scaled_3}\\
\hat{Y}^{(n)}_{j,i}(t)&=\hat{Y}^{(n)}_{j,i}(0) -\hat{M}_{D,j,i}^{(n)}(t)-\int_0^t \mu_j(\hat{Y}^{(n)}_{j,i}(s)-\hat{Y}^{(n)}_{j,i+1}(s))ds,\ i\geq 3, \ j\in[M],\label{eqn:diff_SA-JSQ_evolu_scaled_4}
\end{align}
where,
\begin{align*}
    \hat{V}_2^{(n)}(t) &= \hat{V}_1^{(n)}(t) - \frac{\hat{U}_{1,2}^{(n)}(t)}{\sqrt{n}}.
\end{align*}
From \eqref{eqn:u_11} and \eqref{eqn:u_1j}, we see that $\hat{V}_{1}^n(t)$ is a non-decreasing function increasing only when an arrival occurs and $\hat{Q}_{j,1}^{(n)}=N_j^{(n)}$ for all $j\in [M]$, i.e., for all $t\geq 0$, 
\begin{align}\label{eqn:diff_SA-JSQ_evolu_scaled_bound1}
\int_0^t\indic{\hat{Y}_{[1,M],1}^{(n)}(s)<0}d\hat{V}_1^{(n)}(s) = 0.
\end{align}
Also, by definition $\hat{Y}_{1,2}^{(n)}(t)\leq N_{1}^n/\sqrt{n}$ and $\hat{V}_{2}^n(t)$ is a non-decreasing function and only increases when there is an arrival, $\hat{Q}_{j,1}^{(n)}=N_j^{(n)}$ for all $j\in [M]$ and $\hat{Q}_{1,2}^{(n)}=N_1^{(n)}$. Hence, for all $t\geq 0$, 
\begin{align}\label{eqn:diff_SA-JSQ_evolu_scaled_bound2}
\int_0^t\indic{\hat{Y}_{1,2}^{(n)}(s)<\frac{N_1^{(n)}}{\sqrt{n}}}d\hat{V}_2^{(n)}(s) = 0.
\end{align}

Equations \eqref{eqn:diff_SA-JSQ_evolu_scaled_1_2}-\eqref{eqn:diff_SA-JSQ_evolu_scaled_bound2} have the same form as \eqref{eqn:map_1}-\eqref{eqn:map_6} with $B=B^{(n)}$, $z_1=\hat{Z}_1^{(n)}$ and $z_{j,i}=\hat{Z}_{j,i}^{(n)}$, where
\begin{align*}
 B^{(n)}&=N_1^{(n)}/\sqrt{n},\\
\nonumber \hat{Z}_1^{(n)}(t)&=\hat{M}_{A}^{(n)}(t)-\sum_{j\in[M]}\hat{M}_{D,j,1}^{(n)}(t) -\sqrt{n}(1-\lambda^{(n)}) t +\sum_{j=1}^{M-1}\int_0^t\mu_j\hat{Y}^{(n)}_{j,1}(s)ds\\&\quad  -\mu_M \int_0^t (\hat{Y}^{(n)}_{M,1}(s)-\hat{Y}^{(n)}_{[1:M],1}(s))ds,\\
 \hat{Z}_{j,i}^{(n)}(t)&=-\hat{M}_{D,j,i}^{(n)}(t), \ i\geq 2, \ j\in[M].
 \end{align*}
Once we show that these terms converge to the appropriate limits, we can use the continuous mapping theorem to prove Theorem~\ref{thm:SA-JSQ_diffusion_limit}. The convergence of the third term of $\hat{Z}_1^{(n)}(t)$ is implied by \eqref{eqn:asmp_1}. Moreover, Proposition~\ref{prop:state_space_collapse} implies that the fourth and the fifth terms converge to 0 as $n\to\infty$. Hence, we focus on the limits of the martingales $\hat{M}_{A}^{(n)}(t)$ and $\hat{M}_{D,j,1}^{(n)}(t)$ in Section~\ref{sec:martingale_convergence}. 
\subsection{Martingale Convergence}\label{sec:martingale_convergence}

The martingales $\hat{M}_A^{(n)}$ and $\hat{M}_{D,j,i}^{(n)}$ for all $i\in \ZZ_+$ and $j\in [M]$ are obtained by centering the time-changed unit Poisson processes with the corresponding random time changes. Hence, we can write the predictable quadratic variation process for the martingale $\hat{M}_A^{(n)}$ as 
\begin{align}
    \langle \hat{M}_A^{(n)}\rangle(t) &= \lambda^{(n)}t \label{eqn:M_A}
\end{align}
Assumption~\eqref{eqn:heavy_traffic_limit} clearly implies that $\langle \hat{M}_A^{(n)}\rangle(t)\to t$ as $n\to \infty$. The following lemma characterizes the limit of $\langle \hat{M}_{D,j,i}^{(n)}\rangle (t)$.
\begin{lemma}
\label{lem:fluid_limit}
For any $j\in [M]$ as $n\to \infty$ we have  
\begin{align}
\label{eqn:wc_1}
\frac{\hat{Q}_{j,1}^{(n)}}{n} \Rightarrow  \gamma_j e,
\end{align}
where $e(t)=1$ for all $t\geq0$.
Moreover, for $i\geq 2$ and for $j\in[M]$ as $n\to \infty$, we have 
\begin{align}
\label{eqn:wc_2}
\frac{\hat{Q}_{j,i}^{(n)}}{n} \Rightarrow  0.
\end{align}
\end{lemma}

To prove Lemma~\ref{lem:fluid_limit}, we need the following version of Lemma 9 in \cite{Eschenfeldt2018} modified by incorporating the necessary coefficients. We provide its proof in Appendix~\ref{append:2} for completeness. 

\begin{lemma}
\label{lem:stochastic_boundedness}
For each $n\geq 1$, the following equations 
\begin{align*}
&\hat{Y}_{[1,M],1}^{(n)}(t) =\hat{Y}_{[1,M],1}^{(n)}(0) +\hat{Z}_1^{(n)}(t) - \int_0^t(\mu_M \hat{Y}_{[1,M],1}^{(n)}(s) - \mu_1 \hat{Y}_{1,2}^{(n)}(s))ds -\hat{V}_1^{(n)}(t), \\
&\hat{Y}_{1,2}^{(n)}(t) =\hat{Y}_{1,2}^{(n)}(0) +\hat{Z}_{1,2}^{(n)}(t) - \int_0^t\mu_1 \hat{Y}_{1,2}^{(n)}(s)ds +\hat{V}_1^{(n)}(t)-\hat{V}_2^{(n)}(t),\\
&\hat{Y}_{[1,M],1}^{(n)}(t)\leq0, \ 0\leq \hat{Y}_{1,2}^{(n)}(t) \leq B_n, \ t\geq0, \\
&\int_0^t\indic{\hat{Y}_{[1,M],1}^{(n)}(s)<0}d\hat{V}^{(n)}_1(s)=0,\\
&\int_0^t\indic{\hat{Y}_{1,2}^{(n)}(s)<B_n}d\hat{V}^{(n)}_2(s)=0,
\end{align*}
have a unique solution $(\hat{Y}_{[1,M],1}^{(n)},\hat{Y}_{1,2}^{(n)},\hat{V}_1^{(n)},\hat{V}_2^{(n)})$. Moreover, if the sequences $\cbrac{\hat{Y}_{[1,M],1}^{(n)}(0)}_{n\geq 1}$, $\cbrac{\hat{Y}^{(n)}_{1,2}(0)}_{n\geq 1}$, $\cbrac{\hat{Z}_1^{(n)}}_{n\geq 1}$ and $\cbrac{\hat{Z}_{1,2}^{(n)}}_{n\geq 1}$ are stochastically bounded, then the sequences $\cbrac{\hat{Y}_{[1,M],1}^{(n)}}_{n\geq 1}$  and $\cbrac{\hat{Y}^{(n)}_{1,2}}_{n\geq 1}$ are also stochastically bounded. 
\end{lemma}

{\em Proof of Lemma~\ref{lem:fluid_limit}.}
Suppose $\{\hat{Y}_{[1,M],1}^{(n)}\}_n$ and $\{\hat{Y}_{j,i}^{(n)}\}_n$ for all $i\geq 2$ and $j\in [M]$ are stochatically bounded. Then by Lemma 5.9 of~\cite{Pang2007} we have
\[
\frac{\hat{Y}_{[1,M],1}^{(n)}}{\sqrt{n}}=\frac{\sum_{j\in[M]}(\hat{Q}_{j,1}^{(n)}-N_j^{(n)})}{n}\Rightarrow 0\mbox{ and } \frac{\hat{Y}_{j,i}^{(n)}}{\sqrt{n}} =\frac{\hat{Q}_{j,i}^{(n)}}{n}\Rightarrow 0.
\]
Applying \eqref{eq:proportion_assumption} to the above proves the lemma. Hence, it remains to establish that $\{\hat{Y}_{[1,M],1}^{(n)}\}_n$ and $\{\hat{Y}_{j,i}^{(n)}\}_n$ for all $i\geq 2$ and $j\in [M]$ are stochatically bounded.

Due to truncation, we know that 
\begin{align*}
\hat{Y}_{j,i}^{(n)}(t)&=\hat{Q}_{j,i}^{(n)}(t)/\sqrt{n}   \leq \hat{Q}_{j,i}^{(n)}(0)/\sqrt{n}.
\end{align*}
for all $j\in [M]$ and $i\geq 3$ and $j\in\{2,\ldots, M\}$ and $i\geq 2$. Hence, using Lemmas \ref{lem:stochastic_bounded} and ~\ref{lem:stochastic_boundedness}, we only need to prove that $\hat{Z}_{1}^{(n)}$ and $\hat{Z}_{1,2}^{(n)}$ are stochatiscally bounded, which follows if the respective predictive quadratic variations are stochastiscally bounded. From~\eqref{eqn:M_A} it follows immediately that the sequence $\{\langle \hat{M}_A^{(n)}\rangle(t)\}_{n \geq 1}$ is bounded for each $t\geq 0$, and we have 
\begin{align*}
\langle \hat{M}_{D,j,i}^{(n)}\rangle (t) &\leq \frac{1}{n} \int_0^t \mu_j\hat{Q}_{j,i}^{(n)}(s)ds + \frac{1}{n} \int_0^t \mu_j\hat{Q}_{j,i+1}^{(n)}(s)ds\\
&\leq \frac{2\mu_j}{n}N_j^{(n)}t.
\end{align*}
Hence, the stochastic boundedness of the sequence $\cbrac{\langle \hat{M}_{D,j,i}^{(n)}\rangle (t)}_{n\geq1}$ for each $t\geq0$ follows from~\eqref{eq:proportion_assumption}  and Proposition~\ref{prop:state_space_collapse}. Hence, the lemma follows. 
\qed
\endproof

The following lemma characterizes the limits of the martingales $\hat{M}_A^{(n)}$ and $\hat{M}_{D,j,i}^{(n)}$ as $n\to\infty$. 


\begin{lemma}
\label{lem:mart_convergence}
The sequence of scaled martingales converges weakly to following limits as $n\to \infty$
\begin{multline*}
\brac{\hat{M}_A^{(n)}, \hat{M}_{D,1,1}^{(n)},\dots,\hat{M}_{D,M,1}^{(n)},\dots,\hat{M}_{D,1,i}^{(n)},\dots,\hat{M}_{D,M,i}^{(n)},\dots} \Rightarrow \\\brac{W_A,\mu_1 \gamma_1W_{D_{1,1}},\dots,\mu_M \gamma_MW_{D_{M,1}},\dots,0,\dots,0,\dots},
\end{multline*}
where $W_A$ and $W_{D,j,1}$ for $j\in[M]$ are independent standard Brownian motions.
\end{lemma}
To prove this lemma, we utilise the continuous mapping theorem and functional central limit theorem (FCLT) for Poisson processes, which we restate here for completeness.

\begin{lemma}[Theorem 4.2~\cite{Pang2007}].
\label{lem:fclt_pois}
Let $A$, $D_{j,i}$ for $ i\geq1$ and $j\in[M]$, are independent unit rate Poisson processes and 
\begin{equation*}
\hat{M}_{C,n}=\frac{C(nt)-nt}{\sqrt{n}},
\end{equation*}
for $C=A,D_{j,i}$ for $ i\geq1$ and $j\in[M]$, then as $n\to \infty$ we have 
\begin{multline*}
\brac{\hat{M}_{A,n}, \hat{M}_{D_{1,1},n}, \dots, \hat{M}_{D_{M,1},n},\dots,\hat{M}_{D_{1,i},n}, \dots, \hat{M}_{D_{M,i},n},\dots} \Rightarrow \\ \brac{W_A, W_{D_{1,1}}, \dots, W_{D_{M,1}},\dots,W_{D_{1,i}}, \dots, W_{D_{M,i}},\dots},    
\end{multline*}
where $ W_A $ and $ W_{D_{j,i}} $ are independent standard Brownian motions.

\end{lemma}

{\em Proof of Lemma~\ref{lem:mart_convergence}.}
Using the predictable quadratic variation processes, we can express the martingales $ \hat{M}_A^{(n)} $ and $ \hat{M}_{D,j,i}^{(n)} $ for $ i \geq 1 $ and $ j \in [M] $ as the following compositions: 
\begin{equation*}
\hat{M}_A^{(n)}=\hat{M}_{A,n} \circ  \langle \hat{M}_A^{(n)}\rangle, \  
\hat{M}_{D_{j,i}}^{(n)}=\hat{M}_{D_{j,i},n} \circ \langle \hat{M}_{D,j,i}^{(n)}\rangle, \ i\geq1, \ j\in[M]. 
\end{equation*}
Now, applying CMT together with Lemma~\ref{lem:fclt_pois} and Lemma~\ref{lem:fluid_limit}, we get
\begin{multline*}
\brac{\hat{M}_A^{(n)}, \hat{M}_{D,1,1}^{(n)},\dots,\hat{M}_{D,M,1}^{(n)},\hat{M}_{D,1,2}^{(n)},\dots,\hat{M}_{D,M,2}^{(n)},\dots}\\
\Rightarrow \brac{W_A \circ e, W_{D_{1,1}} \circ \mu_1 \gamma_1 e, \dots, W_{D_{M,1}}\circ \mu_M \gamma_M e,W_{D_{1,2}} \circ 0, \dots, W_{D_{M,2}} \circ 0,\dots},
\end{multline*}
which proves the lemma. 
\qed
\endproof

Now, we  have all the necessary results to prove Theorem~\ref{thm:SA-JSQ_diffusion_limit}.
{\em Proof of Theorem~\ref{thm:SA-JSQ_diffusion_limit}.}
Note that $B_n\to\infty$ as $n\to \infty$. Therefore, we have $\hat{V}_2^{(n)}$ converges to $0$ in the limit. Furthermore, Lemma~\ref{lem:mart_convergence}, Proposition~\ref{prop:state_space_collapse}, and~\eqref{eqn:heavy_traffic_limit} imply 
\begin{align*}
&\hat{M}_{A}^{(n)}(t)-\sum_{j\in[M]}\hat{M}_{D,j,1}^{(n)}(t) -\sqrt{n}(1-\lambda^{(n)}) t+\sum_{j=1}^{M-1}\int_0^t\mu_j\hat{Y}^{(n)}_{j,1}(s)ds\\&\qquad -\mu_M \int_0^t (\hat{Y}^{(n)}_{M,1}(s)-\hat{Y}_{[1,M],1}^{(n)}(s))ds\\
&\Rightarrow W_A(t) -\sum_{j\in[M]}\mu_j \gamma_j W_{D_{j,1}}(t) -\beta t \\&\numeq{d}\sqrt{2}W(t) -\beta t,
\end{align*}
where $\numeq{d}$ indicates equivalence in distribution and $W$ is the standard Weiner process. Moreover, Lemma~\ref{lem:mart_convergence} also implies that martingales terms in~\eqref{eqn:diff_SA-JSQ_evolu_scaled_3}-\eqref{eqn:diff_SA-JSQ_evolu_scaled_4} converges weakly to $0$ as $n\to \infty$.
Now using CMT, as $n\to \infty$ we have $\hat{Y}_{[1,M],1}^{(n)}\Rightarrow Y_{M,1}$, $\hat{Y}_{1,2}^{(n)}\Rightarrow Y_{1,2}$,     
 $\hat{Y}_{j,2}^{(n)}\Rightarrow Y_{j,2}$, for $j\in \cbrac{2,\dots,M}$ and 
 $\hat{Y}_{j,i}^{(n)}\Rightarrow Y_{j,i}$  for $\ i\geq3, \ j\in [M]$.
Hence, for the truncated system, we have $\hat{\mf Y}^{(n)} \Rightarrow \mf Y$. 

As the final step, we prove that both the truncated and the original system have the same limiting dynamics.
To show this, we define 
\begin{equation*}
T_n= \inf \cbrac{t\geq0:Q_{1,2}^{(n)}(t)=N_1^{(n)}}.
\end{equation*}
It is easy to see that starting from the same initial state at $t=0$, the original system and that of the truncated system remain identical for all $t\in[0,T_n)$ since, in both systems, no arrival is sent to a busy server in the slowest $M-1$ pools until all servers in the fastest pool have two jobs. Therefore, it suffices to prove that $\mathbb{P}(T_n\leq t) \to 0, \ as \ n\to \infty$ for all $t\geq 0$.

Since the truncated and the original systems are identical up to $T_n$, we have
\begin{equation*}
\mathbb{P}\brac{T_n\leq t}=\mathbb{P}\brac{\sup_{0\leq s\leq t}\hat{Q}_{1,2}^{(n)}(s)\geq N_1^{(n)}}=\mathbb{P}\brac{\sup_{0\leq s\leq t}\hat{Y}_{1,2}^{(n)}(s)\geq N_1^{(n)}/\sqrt{n}}. 
\end{equation*}
Hence, we obtain
\begin{equation*}
\limsup_{n\to\infty}\mathbb{P}\brac{T_n\leq t}\leq \limsup_{n} \mathbb{P}\brac{\sup_{0\leq s\leq t}\hat{Y}_{1,2}^{(n)}(s)\geq N_1^{(n)}/\sqrt{n}}=0,
\end{equation*}
where the last equality follows as the sequence $\{{\hat Y}_{1,2}^{(n)}(s)\}_{n\in \ZZ_+}$ converges weakly, hence is stochastically bounded and  $\frac{N_1^{(n)}}{\sqrt{n}}$ is of the order of $O(\sqrt{n})$. This completes the proof. 
\qed
\endproof

\section{Stationary Analysis: Proof of Theorem~\ref{thm:stationary_convergence}}\label{sec:stationary_convergence}

The process level convergence in Theorem~\ref{thm:SA-JSQ_diffusion_limit} only states convergence on compact time intervals and does not imply the convergence of stationary measures. To prove the interchangeability of many-servers limit and the limit as time goes to infinity, and obtain the convergence of stationary distributions in Corollary~\ref{cor:stationary_convergence}, we need to prove the tightness of the stationary distributions for the scaled system sizes under consideration as stated in Theorem~\ref{thm:stationary_convergence}.

{ To prove~\eqref{eq:Expectation_1} of Theorem~\ref{thm:stationary_convergence}, we establish through an appropriate Lyapunov function that $\EE\sbrac{\sum_{j=1}^{M-1}|Y_{j,1}^{(n)}(\infty)|}=O(1/\sqrt{n})$. Then, leveraging rate conservation and~\eqref{eq:Expectation_1}, we establish $\EE [|Y_{M,1}^{(n)}(\infty)|]=O(1)$, thereby confirming~\eqref{eq:Expectation_2}. Additionally, employing Stein's method, we prove that $\EE[Y_{1,2}^{(n)}(\infty)]=O(1)$, which establishes~\eqref{eq:Expectation_3}. Finally using a bootstrapping argument as in~\cite{Braverman2020}, we prove~\eqref{eq:Expectation_4}.


Throughout the proof we shall be using the infinitesimal generator $\cG$ of the process $\mf Q^{(n)}$. To express the generator compactly, we write $(j,i)\prec (k,l)$ if $i<l$, or $i=l$ and $j<k$. We also use $(k,l)\succ (j,i)$ if $(j,i)\prec (k,l)$. The infinitesimal generator $\cG$ applied to a function $f:S^{(n)}\to \RR$ and evaluated at a state $\mf q \in S^{(n)}$ can be written as 
\begin{align}\label{eq:basic_generator}
    \nonumber\cG f(\mf q) &= n\lambda^{(n)} \indic{q_{1,1}<N_1^{(n)}}(f(\mf q+\mf e_{1,1})-f( \mf q))\\
    \nonumber&\quad + \sum_{(k,i)\succ (1,1)}n\lambda^{(n)}\indic{q_{l,j}=N_l^{(n)}\mbox{ for }(l,j)\prec (k,i),q_{k,i}<N_k^n}(f(\mf q+\mf e_{k,i})-f(\mf q))\\
    &\quad + \sum_{k=1}^\infty\sum_{i=1}^\infty \mu_k (q_{k,i}-q_{k,i+1})(f(\mf q-\mf e_{k,i})-f(\mf q)),
\end{align}
where $\mf e_{i,j}$ denotes the infinite dimensional 
unit vector with one
in the $(i,j)$th position. We will be using the following rate conservation result (see e.g. Lemma~1 of~\cite{Braverman2020}) which states that for any $f:S^{(n)}\to \RR$ with $\EE[|f(\mf Q^{(n)}(\infty))|]< \infty$ we have
\begin{equation}
    \EE[\cG f(\mf Q^{(n)}(\infty))]=0.
    \label{eq:rate_cons}
\end{equation}
The above equation plays a key role in the proof using Stein's approach which chooses an appropriate $f$ to plug into~\eqref{eq:rate_cons}.

To use~\eqref{eq:rate_cons} within Stein's approach, for each $t\in [0,\infty]$ we define $$X_1^{(n)}(t)=\left(\sum_{k=1}^MQ_{k,1}^{(n)}(t)-N_k^n\right)/n, \text{ and } X_2^{(n)}(t)=Q_{1,2}^{(n)}(t)/n.$$ 
Furthermore, for any function $f:\RR^2\to \RR$, we define its lifted version $Af: S^{(n)}\to \RR$ as $Af(\mf q)=f(\mf x)$ where $\mf q\in S^{(n)}$ and $\mf x$ is the value of $\mf X^{(n)}(t)$ when $\mf Q^{(n)}(t)=\mf q$. Then, the generator acting on the lifted function can be written as
\begin{align}\label{eq:generator_lyapunov}
    \nonumber\cG Af( \mf q)&=n\lambda^{(n)}\indic{x_1 < 0}\left(f(\mf x+\mf e_1/n)-f(\mf x)\right)
    +n\lambda^{(n)} \indic{x_1=0}\left(f(\mf x+\mf e_2/n)-f(\mf x)\right)\\
    \nonumber&\quad-n\lambda^{(n)} \indic{x_1=0, x_2=N_1^{(n)}/n}\left(f(\mf x+\mf e_2/n)-f(\mf x)\right)\\
    &\quad+\sum_{k=1}^M\mu_k(q_{k,1}-q_{k,2})\left(f(\mf x-\mf e_1/n)-f(\mf x)\right)+\mu_1(q_{1,2}-q_{1,3})\left(f(\mf x-\mf e_2/n)-f(\mf x)\right),
\end{align}
where $\mf e_i$ is the two-dimensional unit vector with one in the $i$th position. Furthermore, since $\mf X^{(n)}(\infty)$ can only take a finite many values for each $n$ we have $\EE[|f(\mf X^{(n)}(\infty))|]< \infty $ for any $f:\RR^2\to \RR$ which by~\eqref{eq:rate_cons} implies that 
\begin{equation}
\EE[\cG Af(\mf X^{(n)}(\infty))]=0
\label{eq:rate_lifted}
\end{equation}

To find the appropriate function $f$ to plug into the above equation we observe that
for any twice-differentiable $f:\RR^2\to \RR$,
the Taylor expansion yields
 \begin{align*}
    f(\mf x+\mf e_{1}/n)-f(\mf x) &=\frac{1}{n}f_1(\mf x)+\int_{x_1}^{x_1+1/n}(x_1+1/n-u)f_{11}(u,x_2)du,\\
    f(\mf x-\mf e_{1}/n)-f(\mf x) &=-\frac{1}{n}f_1(\mf x)+\int_{x_1-1/n}^{x_1}(u-x_1+1/n)f_{11}(u,x_2)du,
\end{align*}
where $f_i$ denotes the first-order partial derivative $\partial f/\partial x_i$ and $f_{ij}$ denotes the second-order partial derivative $\partial f/\partial x_i \partial x_j$.
Using the above expansion in~\eqref{eq:generator_lyapunov} we have that for any twice-differentiable $f:\RR^2\to \RR$

\begin{align}\label{eq:expanded_generator}
\cG Af(\mf q) & = \cL f(\mf x) - \lambda^{(n)} \indic{x_1=0} (f_1(\mf x)-f_2(\mf x)) + \varepsilon_f^n(\mf q),
\end{align}
where 
\begin{align}\label{eq:PDE_operator}
    \cL f(\mf x) &=\left(-\frac{\beta}{\sqrt{n}}-\mu_M x_1 + \mu_1 x_2\right)f_1(\mf x) - \mu_1 x_2 f_2(\mf x)
\end{align}
and 
\begin{align}
\varepsilon_f^n(\mf q) &= n\lambda^{(n)}\indic{x_1<0}\int_{x_1}^{x_1+1/n}(x_1+1/n-u)f_{11}(u,x_2)du \label{eqn:err1}\\
&\quad+  n\lambda^{(n)}\indic{x_1=0}\int_{x_2}^{x_2+1/n}(x_2+1/n-u)f_{22}(x_1,u)du\label{eqn:err2} \\
&\quad-  n\lambda^{(n)}\indic{x_1=0, x_2=N_1^{(n)}/n}\int_{N_1^{(n)}/n}^{(N_1^{(n)}+1)/n}f_2(0,u)du\label{eqn:err3}\\
&\quad + \sum_{k=1}^{M-1}(\mu_k-\mu_M)(N_k^n-q_{k,1})\int_{x_1-1/n}^{x_1}f_1(u,x_2)du\label{eqn:err4}\\
&\quad+ n\mu_Mx_1\int_{x_1+1/n}^{x_1}(u-x_1+1/n)f_{11}(u,x_2)du\label{eqn:err5}\\
&\quad + \sum_{k=1}^M N_k^n\mu_k\int_{x_1-1/n}^{x_1}(u-x_1+1/n)f_{11}(u,x_2)du\label{eqn:err6}\\
& \quad+ \sum_{k=2}^M\mu_kq_{k,2}\int_{x_1-1/n}^{x_1}f_1(u,x_2)du\label{eqn:err7}\\
&\quad - n\mu_1x_2 \int_{x_1-1/n}^{x_1}(u-x_1+1/n)f_{11}(u,x_2)du\label{eqn:err8}\\ &\quad +n\mu_1x_2\int_{x_2-1/n}^{x_2}(u-x_2+1/n)f_{22}(x_1,u)du\label{eqn:err9}\\
&\quad + \mu_1q_{1,3}\int_{x_2-1/n}^{x_2}f_{2}(x_1,u)du\label{eqn:err10}\\
&\quad+\left(\lambda^{(n)}-\frac{\beta}{\sqrt{n}}-\frac{\sum_{k=1}^{M}\mu_kN_k^{(n)}}{n}\right)f_1(\mf x)\label{eq:annoying}
\end{align}

Thus, the expansion~\eqref{eq:expanded_generator} of the generator $\cG$ shows that for any twice-differentiable $f:\RR^2\to \RR$ the operator $\cG$ acting on the lifted function $Af$ can be approximated by the linear operator $\cL$  (defined in~\eqref{eq:PDE_operator}) acting on $f$ as long as the error terms contained in $\varepsilon_f^n(\mf q)$ remain small and first order partial derivatives $f_1$ and $f_2$ coincide along the boundary $x_1=0$ (which makes the second term appearing the generator expansion~\eqref{eq:expanded_generator} zero). 

The key idea of our proof using Stein's approach is to restrict $f$ to a class of functions for which $f_1(\mf x)$ indeed matches with $f_2(\mf x)$ along the reflection boundary $x_1=0$ and  $\varepsilon_f^n(\mf q)$ remains small for all sufficiently large $n$. For $f$ belonging to such a class of functions we have from~\eqref{eq:rate_lifted}
$$\EE[-\cL f(\mf X^{(n)}(\infty)]=\EE[\varepsilon_f^n(\mf Q^{(n)}(\infty))]\approx 0$$ for all sufficiently large $n$. Thus, to show that a chosen metric $h$ of the steady-state $\mf X^{(n)}(\infty)$ remains close to zero in expectation for all sufficiently large $n$ , i.e., to establish
$$\EE[h(\mf X^{(n)}(\infty))]=\EE[\varepsilon_f^n(\mf Q^{(n)}(\infty))]\approx 0,$$
it is enough to find $f$ satisfying PDE $\cL f(\mf x)=-h(\mf x)$. This is the key idea we employ in our proof.



To prove~\eqref{eq:Expectation_1} and~\eqref{eq:Expectation_4} we will also make use the following lemma from~\cite{wang_etal2022} which provides steady-state bounds on expected value of a Lyapunov function when its drift is negative outside a compact set.

\begin{lemma}[Lemma 10 of~\cite{wang_etal2022}]
\label{lem:drift}
Let $Z=(Z(t), t\geq 0)$ be a continuous time Markov process with countable state-space $\cZ$ and infinitesimal generator $\cG_Z$.
Furthermore, assume that the chain $Z$ is irreducible, nonexplosive, and ergodic and  $Z(t) \Rightarrow Z(\infty)$ as $t\to \infty$.
Let $r_{i,j}$ denote the transition rate of $Z$ from state $i\in \cZ$ to state $j \in \cZ\setminus\{i\}$ and $r_{i,i}=-\sum_{j \in \cZ\setminus\{i\}} r_{i,j}$. Suppose that $ r_{\max}=\sup_{i \in \cZ} |r_{i,i}| < \infty$ and a Lyapunov function $V:\cZ\to \RR_+$ exists such that $\cG_Z V(i) \leq -\gamma$ for some $\gamma >0$ whenever $i\in \cZ$ is such that $V(i)> B$ for some $B >0$. Then,  for any integer $m\in \ZZ_+$ we have
\begin{equation*}
    \EE[V^m(Z(\infty))]\leq (2B)^m+(4\nu_{\max})^m\brac{\frac{r_{\max}\nu_{\max}+\gamma}{\gamma}}^m m!
\end{equation*}
where $\nu_{\max}=\sup_{i,j \in \cZ:r_{i,j}>0} |V(j)-V(i)|$.
\end{lemma}}

{\em Proof of Theorem~\ref{thm:stationary_convergence}.}
Define function 
\begin{equation}
\label{eqn:lyapunov_function}
 V(\mf q) = \sum_{j=1}^{M-1} (N_j^{(n)}-q_{j,1}).   
\end{equation}
 To prove \eqref{eq:Expectation_1}, it is enough to show that $\EE[V(\mf Q^{(n)}(\infty))]=O(1)$. Plugging $V(\mf q)$ into the generator in \eqref{eq:basic_generator}, we have 
\begin{align*}
    \cG V(\mf q) &= -n\lambda^{(n)} \indic{V(\mf q)>0} + \sum_{j=1}^{M-1}\mu_j(q_{j,1}-q_{j,2}).
\end{align*}
Fix $\epsilon>0$, \eqref{eqn:heavy_traffic_limit} and \eqref{eq:proportion_assumption} imply that for large enough $n$ and any $\mf q$ with $V(\mf q)\geq 1$ we have
\begin{align*}
    \cG V(\mf q) &\leq -n \lambda^{(n)} + \sum_{j=1}^{M-1}N_j^{(n)}\mu_j\\
    &\leq -n\ \gamma_M\mu_M + (\beta+\epsilon)\sqrt{n} <0.
\end{align*}
Hence, using Lemma~\ref{lem:drift} for $m\in \ZZ_+$ we have
\begin{equation*}
\EE \sbrac{V(\mf q)^m}\leq (2)^m + (4)^m\brac{\frac{\lambda^{(n)}+\sum_{k=1}^M \gamma_k\mu_k-(\beta+\epsilon)n^{-1/2}}{\gamma_M\mu_M-(\beta+\epsilon)n^{-1/2}}}^m m!.
\end{equation*}
Hence, for any $n>\left(\frac{2(\beta+\epsilon)} {\gamma_M\mu_M}\right)^2$ and $m\in \ZZ_+$
\begin{equation}
\label{eqn:v_q_m}
\EE \sbrac{V(\mf q)^m}\leq (2)^m + \brac{\frac{16}{\gamma_M\mu_M}}^mm!.  
\end{equation}
In particular, for $m=1$,~\eqref{eqn:v_q_m} becomes  
\[
\EE\sbrac{\sqrt{n}\sum_{j=1}^{M-1}|Y_{j,1}^{(n)}(\infty)|}\leq 2 + \frac{16}{\gamma_M\mu_M},
\]
which proves \eqref{eq:Expectation_1}.


Now we turn to the proofs of~\eqref{eq:Expectation_2} and~\eqref{eq:Expectation_3} for which we need the following lemma proved in Appendix~\ref{app6}. 

\begin{lemma}\label{lem:simple_stationary_bounds}
For any $n\geq 1$ we have
\begin{enumerate}
\item $\EE[\sum_{k=1}^M \mu_k Q_{k,1}^{(n)}(\infty)]=n\lambda^{(n)}$,
\item $\EE\left[\sum_{k=2}^M \mu_kQ_{k,2}^{(n)}(\infty)+\mu_1Q_{1,3}^{(n)}(\infty)\right]=n\lambda^{(n)}\PP(\sum_{k=1}^M Q_{k,1}^{(n)}(\infty)=n, Q_{1,2}^{(n)}(\infty)=N_1^{(n)})$.
\end{enumerate}
\end{lemma}

Now, we write
\begin{align*}
    \EE\sbrac{\sqrt{n}|Y_{M,1}^{(n)}(\infty)|}&\leq \mu_M^{-1}\brac{\sum_{k=1}^MN_k^{(n)}\mu_k - \EE\sbrac{\sum_{k=1}^M \mu_k Q_{k,1}^{(n)}(\infty)}}\\
    &= \mu_M^{-1}\brac{\sum_{k=1}^MN_k^{(n)}\mu_k - n\lambda^{(n)}}\\
    &=\sqrt{n}\mu_M^{-1}\brac{\sqrt{n}\sum_{k=1}^M(N_k^{(n)}/n-\gamma_k)\mu_k +\sqrt{n}(1 - \lambda^{(n)})},
\end{align*}
where the second line follows from Part 1 of Lemma~\ref{lem:simple_stationary_bounds}.  The proof of\eqref{eq:Expectation_2} now follows by applying \eqref{eqn:heavy_traffic_limit} and \eqref{eq:proportion_assumption} to the above inequality.

To establish~\eqref{eq:Expectation_3}, it is enough to show that
\begin{equation}
 \EE\sbrac{\brac{X_2^{(n)}(\infty)-\frac{\kappa}{\sqrt{n}}}_+}=O\brac{\frac{1}{\sqrt{n}}}   
 \label{eq:target}
\end{equation}
for any constant $\kappa > \beta/\mu_1$.
To establish \eqref{eq:target}, we 
use~\eqref{eq:rate_lifted} with an appropriate choice of $f$.
%
Lemma~\ref{lem:solution_PDE} below provides us with the appropriate Lyapunov function $f^*(\mf x)$ to plug in to \eqref{eq:rate_lifted} and prove \eqref{eq:target}. We construct the function $f^*(\mf x)$ in Section~\ref{sec:lyapunov_func} and prove that this function satisfies the conditions stated in Lemma~\ref{lem:solution_PDE} in Appendix~\ref{append:lem9}.

\begin{lemma}
\label{lem:solution_PDE}
Fix $\kappa>\beta/\mu_1$ and define $\Omega=(-\infty,0]\times [0,\infty)$. Then, there exists a function $f^*(\mf x)$ that solves the PDE 
\begin{align}\label{eq:PDE_Lyapunov1}
    \cL f^*(\mf x) &= -(x_2-\kappa/\sqrt{n})_+, \mbox{ for all } \mf x\in \Omega,\\
    \label{eq:PDE_Lyapunov2}f_1^*(0,x_2) &= f_2^*(0,x_2)=\frac{\sqrt{n}}{\beta}\left(x_2-\frac{\kappa}{\sqrt{n}}\right)_+, \mbox{ for all }x_2\geq 0,
\end{align}
with absolutely continuous first order derivatives $f_1^*(\cdot, x_2)$ and $f_2^*(x_1, \cdot)$. Moreover, we have
\begin{align}
f_1^*(\mf x)\geq 0, f_{12}^*(\mf x)\geq 0, f_{11}^*(\mf x)\geq 0, f_{22}^*(\mf x)\geq 0, \ \ \ \ \ \mf x\in \Omega, \label{eqn:second_derive_bound_nn}\\
f_{11}(\mf x)=f_{22}(\mf x)=0, \ \  x_2 \in [0,\kappa/\sqrt{n}],\\
f_{11}^*(\mf x)\leq C_5(\kappa)\sqrt{n}, f_{22}^*(\mf x)\leq C_6(\kappa)\sqrt{n}, \ \ \ \ \ \mf x\in \Omega,
\end{align}
where $C_5(\kappa)= \frac{1}{\beta}\frac{\mu_1+\mu_M}{\mu_1} \frac{\kappa}{\kappa-\beta/\mu_1}$ and $C_6(\kappa)=\frac{1}{\mu_1 \kappa}+ \frac{1}{ \beta}\brac{\frac{\mu_1}{\mu_1-\mu_M}}^2\frac{\kappa}{\kappa-\beta/\mu_1}+ \frac{1}{\beta}\frac{\mu_1(1+\mu_2)}{\mu_1-\mu_M}$.
\end{lemma}

Now, for $f^*$ as given in Lemma~\ref{lem:solution_PDE}, 
we have from~\eqref{eq:rate_lifted} 
\begin{align}
    \EE\sbrac{\brac{X_2^{(n)}(\infty)-\frac{\kappa}{\sqrt{n}}}_+}&=\EE[\epsilon_{f^*}^n(\mf Q^{(n)}(\infty))].
    \label{eq:target1}
\end{align}
Suppose that we can show that
the expected value of every individual term appearing in lines~\eqref{eqn:err1}-\eqref{eqn:err10} is $O(1/\sqrt{n})$. Then using~\eqref{eq:target1}, ~\eqref{eqn:err1}-\eqref{eq:annoying}, and the fact that $f_1^*(\mf x)\geq 0$
we have
\begin{align}
    \EE\sbrac{\brac{X_2^{(n)}(\infty)-\frac{\kappa}{\sqrt{n}}}_+}&\leq  O\brac{\frac{1}{\sqrt{n}}}
    \nonumber\\
    &\quad+\brac{\left\vert\lambda^{(n)}-1-\frac{\beta}{\sqrt{n}}\right\vert+\left\vert1-\frac{\sum_{k\in[M]}\mu_kN_k^{(n)}}{{n}}\right\vert}\EE[f^*_1(\mf X^{(n)}(\infty))],\nonumber\\
    &\leq O\brac{\frac{1}{\sqrt{n}}}+ o\brac{\frac{1}{\sqrt{n}}}\EE[f_1^*(0,X_2^{(n)}(\infty))],\nonumber\\
    &\leq O\brac{\frac{1}{\sqrt{n}}}+o(1)\EE\sbrac{\brac{X_2^{(n)}(\infty)-\frac{\kappa}{\sqrt{n}}}_+}
    \label{eq:last}
\end{align}
where the second inequality follows due to assumptions~\eqref{eqn:heavy_traffic_limit} and~\eqref{eq:proportion_assumption} and the fact that $f_{11}^*(\mf x)\geq 0$, and the last inequality follows since $f_1^*(0,x_2)=(\sqrt{n}/\beta)(x_2-\kappa/\sqrt{n})_+$.
From \eqref{eq:last} it follows that
$$\EE\sbrac{(1-o(1))\brac{X_2^{(n)}(\infty)-\frac{\kappa}{\sqrt{n}}}_+}= O\brac{\frac{1}{\sqrt{n}}}$$
which implies~\eqref{eq:target} and hence~\eqref{eq:Expectation_3}. Therefore, it suffices to show that the expected value of every individual term appearing in lines~\eqref{eqn:err1}-\eqref{eqn:err10} is $O(1/\sqrt{n})$.

Using the bounds on the second derivatives from Lemma~\ref{lem:solution_PDE}, we have
\begin{align*}
    n\lambda^{(n)}\indic{x_1<0}\int_{x_1}^{x_1+1/n}(x_1+1/n-u)f_{11}^*(u,x_2)du &\leq C_5(\kappa)\lambda^{(n)}/\sqrt{n}.
\end{align*}
All other terms involving second derivatives can be bounded in a similar fashion.

The second derivatives of $f^*(\mf x)$ being non-negative implies that the first derivatives are non-decreasing functions in $x_1$ and $x_2$ and hence,
\begin{align}\label{eq:monotonicity_f}
f_1^*(\mf x)\leq f_1^*\left(0,\frac{N_1^{(n)}}{n}\right)\mbox{ and }f_2^*(\mf x)\leq f_2^*\left(0,\frac{N_1^{(n)}}{n}\right)=f_1^*\left(0,\frac{N_1^{(n)}}{n}\right)   
\end{align}
for all $x_1\leq 0$ and $0\leq x_2\leq \frac{N_1^{(n)}}{n}$.
Hence for the term appearing in~\eqref{eqn:err4} we have
\begin{align*}
    &\sum_{k=1}^{M-1}(\mu_k-\mu_M)\EE\left[(N_k^n-Q_{k,1}^{(n)}(\infty))\int_{X_1^{(n)}(\infty)-1/n}^{X_1^{(n)}(\infty)}f_1^*(u,X_2^{(n)}(\infty))du\right]\\
    &\qquad\qquad\qquad\qquad\qquad\qquad\qquad\qquad\leq \frac{\mu_1}{\sqrt{n}} \EE\left[\left|\sum_{k=1}^{M-1}Y_{k,1}^{(n)}(\infty)\right|\right]f_1^*\left(0,\frac{N_1^{(n)}}{n}\right)\\
    &\qquad\qquad\qquad\qquad\qquad\qquad\qquad\qquad\leq \frac{\mu_1 C_1(M-1)}{\beta\sqrt{n}}\brac{\frac{N_1^{(n)}}{n}-\frac{\kappa}{\sqrt{n}}},
\end{align*}
where the first inequality follows from~\eqref{eq:monotonicity_f} and the last inequality follows from~\eqref{eq:Expectation_1} and~\eqref{eq:PDE_Lyapunov2}. From the above inequality it is clear that the right hand side converges to $0$ as $n\to \infty$. 

Using \eqref{eq:PDE_Lyapunov2}, we can bound the remaining terms appearing in \eqref{eqn:err3}, \eqref{eqn:err7}, and \eqref{eqn:err10} as
\begin{align*}
&-n\lambda^{(n)}\indic{x_1=0, x_2=\frac{N_1^{(n)}}{n}}\int_{\frac{N_k^{(n)}}{n}}^{\frac{N_k^{(n)}+1}{n}}f_2^*(0,u)du\\
&\qquad\qquad\qquad\qquad\qquad\qquad\qquad\qquad\qquad\leq-\lambda^{(n)}\indic{x_1=0, x_2=\frac{N_1^{(n)}}{n}}f_2^*\brac{0,\frac{N_1^{(n)}}{n}}\\
&\qquad\qquad\qquad\sum_{k=2}^M\mu_kq_{k,2}\int_{x_1-\frac{1}{n}}^{x_1}f_1^*(u,x_2)du \leq  \frac{1}{n}\sum_{k=2}^M\mu_kq_{k,2}f_2^*\brac{0,\frac{N_1^{(n)}}{n}}\\
&\qquad\qquad\qquad\qquad\mu_1q_{1,3}\int_{x_2-\frac{1}{n}}^{x_2}f_{2}^*\brac{x_1,u}du\leq 
\frac{1}{n}\mu_1q_{1,3}f_2^*\brac{0,\frac{N_1^{(n)}}{n}}.
\end{align*}
Now, replacing $\mf q$ with $\mf Q^{(n)}(\infty)$, summing the above inequalities and taking expectations both sides with respect to $\pi^{n}$ we see that the sum of the right-hand sides is equal to $0$ due to Part 2 of Lemma~\ref{lem:simple_stationary_bounds}. Hence, the contribution of these terms to $\EE[\epsilon_{f^*}^n(\mf Q^{(n)}(\infty))]$ is non-positive and can be safely ignored to obtain an upper bound. This completes the proof of~\eqref{eq:Expectation_3}.

{ To establish \eqref{eq:Expectation_4}, it is sufficient to prove that 
for any $1/2<a<1$ there exists a $C_8(a)>0$ such that for all sufficiently large $n$ 
\begin{equation*}
\EE\left[\sum_{k=2}^M \mu_kQ_{k,2}^{(n)}(\infty)+\mu_1Q_{1,3}^{(n)}(\infty)\right] \leq n^{1-a} C_8(a).
\end{equation*}
Using Part 2 of Lemma~\ref{lem:simple_stationary_bounds}, the above inequality follows if we can show that for any $1/2<a<1$ there exists a $C_8(a)$ such that for all sufficiently large $n$ the following holds
\begin{equation}
n\lambda^{(n)}\PP\brac{\sum_{k=1}^M Q_{k,1}^{(n)}(\infty)=n, Q_{1,2}^{(n)}(\infty)=N_1^{(n)}}\leq n^{1-a}C_8(a).
\end{equation}
By assumption~\eqref{eq:proportion_assumption} for any fixed $\epsilon >0$, we have $|N_1^{(n)}/n-\gamma_1| < \epsilon$ and $\max(\beta/(\mu_1\sqrt{n}),1/n)<\gamma_1-2\epsilon$ for all sufficiently large $n$.  Choose $\tilde{\kappa}$ in the interval $(\gamma_1-2\epsilon,\gamma_1-\epsilon)$. Hence, for all sufficiently large $n$ we have $\tilde\kappa \in (\max(\beta/(\mu_1\sqrt{n}),1/n),N_1^{(n)}/n)$. We have
\begin{align}
    \PP\brac{Q_{1,2}^{(n)}(\infty)=N_1^{(n)}}&=\PP\brac{X_2^{(n)}(\infty)=\frac{N_1^{(n)}}{n}}\\&=\frac{n}{N_1^{(n)}-\tilde{\kappa}n}\EE\sbrac{\brac{X_2^{(n)}(\infty)-\tilde{\kappa}}\indic{X_2^{(n)}(\infty)=\frac{N_1^{(n)}}{n}}}\\
    & \leq \frac{n}{N_1^{(n)}-\tilde{\kappa}n}\EE\sbrac{(X_2^{(n)}(\infty)-\tilde{\kappa})_+}.
\end{align}

Consider the $\tilde{f}^*(\mf x)$ and $\varepsilon_{\tilde{f}^*}^n(\mf{q})$ as in Lemma~\ref{lem:solution_PDE} defined by replacing $\kappa = \sqrt{n}\tilde{\kappa} > \beta/\mu_1$. Using  \eqref{eq:target1} with $\varepsilon_{\tilde{f}^*}^n(\mf{q})$, the result will follow if we can show that for any $1/2<a<1$ there exists a $C_8(a)$ such that $\EE[\varepsilon_{\tilde{f}^*}^n({\mf Q}^{(n)}(\infty))]\leq C_8(a)/n^a$ for all large enough $n$.

The terms in \eqref{eqn:err1}, \eqref{eqn:err2}, \eqref{eqn:err5}, \eqref{eqn:err6} and \eqref{eqn:err9} of $\varepsilon_{\tilde{f}^*}^n({\mf q}^{(n)}(\infty))$ are equal to 0 when $x_2 < \tilde{\kappa}-1/n$. On the other hand, when \( x_2 \geq \tilde{\kappa} -1/n \), we can find a $C_9>0$ such that for all large enough $n$
\begin{align*}
n\lambda^{(n)}\indic{x_1<0}\int_{x_1}^{x_1+1/n}(x_1+1/n-u)\tilde{f}^*_{11}(u,x_2)du&\leq \frac{C_5(\sqrt{n}\tilde{\kappa})}{\sqrt{n}}<\frac{C_9}{\sqrt{n}}, \\
n\lambda^{(n)}\indic{x_1=0}\int_{x_2}^{x_2+1/n}(x_2+1/n-u)\tilde{f}^*_{22}(x_1,u)du&\leq \frac{C_6(\sqrt{n}\tilde{\kappa})}{\sqrt{n}}<\frac{C_9}{\sqrt{n}},\\
n\mu_Mx_1\int_{x_1+1/n}^{x_1}(u-x_1+1/n)\tilde{f}^*_{11}(u,x_2)du &\leq \frac{C_5(\sqrt{n}\tilde{\kappa}) \mu_1}{\sqrt{n}} <\frac{C_9}{\sqrt{n}}, \\
\sum_{k=1}^M N_k^n\mu_k\int_{x_1-1/n}^{x_1}(u-x_1+1/n)\tilde{f}^*_{11}(u,x_2)du &\leq  \frac{C_5(\sqrt{n}\tilde{\kappa})}{\sqrt{n}} \left( 1 + o\left( \frac{1}{\sqrt{n}} \right) \right)<\frac{C_9}{\sqrt{n}}, \\
n\mu_1x_2\int_{x_2-1/n}^{x_2}(u-x_2+1/n)\tilde{f}^*_{22}(x_1,u)du &\leq  \frac{C_6(\sqrt{n}\tilde{\kappa}) \mu_1}{\sqrt{n}}<\frac{C_9}{\sqrt{n}}.
\end{align*}

Next, we prove a bound on the expected value of the third term of $\varepsilon_{\tilde{f}^*}^n(\mf q)$, \eqref{eqn:err4}. Taking the expectation of \eqref{eqn:err4} with respect to \( \pi^n \), we can write

\begin{align}
\label{eqn:56_bound}
    &\sum_{k=1}^{M-1}(\mu_k-\mu_M)\EE\left[(N_k^{(n)}-Q_{k,1}^{(n)}(\infty))\int_{X_1^{(n)}(\infty)-1/n}^{X_1^{(n)}(\infty)}\tilde{f}_1^*(u,X_2^{(n)}(\infty))du\right]\nonumber\\
 &\leq \mu_1 \EE\left[\sum_{k=1}^{M-1}(N_k^{(n)}-Q_{k,1}^{(n)}(\infty))\int_{X_1^{(n)}(\infty)-1/n}^{X_1^{(n)}(\infty)}\tilde{f}_1^*(u,X_2^{(n)}(\infty))du\right] \nonumber\\
 &\leq \frac{\mu_1}{n}\EE\left[\sum_{k=1}^{M-1}(N_k^{(n)}-Q_{k,1}^{(n)}(\infty))\tilde{f}_1^*\left(0,X_2^{(n)}(\infty)\right)\right]\nonumber\\
&= \frac{\mu_1}{\beta \sqrt{n}}\EE\left[\sum_{k=1}^{M-1}(N_k^{(n)}-Q_{k,1}^{(n)}(\infty))\left(X_2^{(n)}(\infty)-\tilde{\kappa}\right)_+\right]\nonumber\\
& =\frac{\mu_1}{\beta \sqrt{n}}\EE\left[\sum_{k=1}^{M-1}(N_k^{(n)}-Q_{k,1}^{(n)}(\infty))\left(X_2^{(n)}(\infty)-\tilde{\kappa}\right) \indic{X_2^{(n)}(\infty)>\tilde{\kappa}}\right]\nonumber\\
&\leq\frac{\mu_1}{\beta \sqrt{n}}\EE\left[\sum_{k=1}^{M-1}(N_k^{(n)}-Q_{k,1}^{(n)}(\infty))\indic{X_2^{(n)}(\infty)>\tilde{\kappa}}\right]\nonumber\\
&\leq\frac{\mu_1}{\beta \sqrt{n}} \brac{\EE\sbrac{V(\mf Q^{(n)}(\infty))^p}}^{1/p}
\brac{\PP\brac{X_2^{(n)}(\infty)>{\tilde{\kappa}}}}^{1/q},
\end{align}
where the first equality follows from~\eqref{eq:PDE_Lyapunov2}, the second inequality follows from~\eqref{eq:monotonicity_f},  the second to last inequality follows as $X_2^{(n)}(\infty)<1$, and in the last inequality we have used~\eqref{eqn:lyapunov_function}, and the Holder's inequality with $p,q>1$ satisfying $1/p+1/q=1$.  
Let $a=1/2+1/2q$ and define $p= \frac{1}{2-2a}$, then using~\eqref{eqn:v_q_m} for $m=\lceil p \rceil $, we have 
\begin{align*}
\EE \sbrac{V(\mf Q^{(n)}(\infty) )^p}\leq \EE[V(\mf Q^{(n)}(\infty))^{\lceil p \rceil }]^{\frac{p}{\lceil p \rceil}} \leq \brac{(2)^{\lceil p \rceil} + \brac{\frac{16}{\gamma_M\mu_M}}^{\lceil p \rceil} (\lceil p \rceil)!}^{\frac{p}{\lceil p \rceil}}=C_{10}(a). 
\end{align*}
Furthermore, using 
the Markov inequality and~\eqref{eq:Expectation_3}, for all sufficiently large $n$
$$\PP(X_2^{(n)}(\infty)\geq \tilde{\kappa})\leq \frac{\EE[X_2^{(n)}(\infty)]}{\tilde{\kappa}}\leq \frac{C_3}{\tilde{\kappa}\sqrt{n}}.$$
Hence, the right-hand side of~\eqref{eqn:56_bound} can be bounded as
\begin{equation*}
    \sum_{k=1}^{M-1}(\mu_k-\mu_M)\EE\left[(N_k^n-Q_{k,1}^{(n)}(\infty))\int_{X_1^{(n)}(\infty)-1/n}^{X_1^{(n)}(\infty)}f_1^*(u,X_2^{(n)}(\infty))du\right]\leq \frac{\mu_1C_3 C_{10}(a)}{\beta \tilde{\kappa}^{1/q} n^{a}}. 
\end{equation*}
for any $a \in (1/2,1)$ and all sufficiently large $n$.

In proving \eqref{eq:Expectation_3}, we have already shown that the expected sum of the terms appearing in~\eqref{eqn:err3}, \eqref{eqn:err7}, and \eqref{eqn:err10} is at most $0$ for any $\kappa$ and, therefore, we can safely ignore the contributions of these terms to obtain an upper bound on $\EE[\varepsilon_{f^*}^n(\mf Q^{(n)}(\infty))]$. Hence, combining all the bounds obtained, we have
\begin{equation}
\EE[\varepsilon_{\tilde{f}^*}^n(\mf Q^{(n)}(\infty))]\leq \frac{4C_9}{\sqrt{n}}\PP\brac{X_2^{(n)}(\infty)>\tilde{\kappa}-\frac{1}{n}} + \frac{\mu_1C_3 C_{10}(a)}{\beta \tilde{\kappa}^{1/q} n^{a}}.
\end{equation}
But using Markov inequality and~\eqref{eq:Expectation_3} again, we have
\begin{equation}
    \PP\brac{X_2^{(n)}(\infty)>\tilde{\kappa}-\frac{1}{n}} = \PP\brac{Y_2^{(n)}(\infty)>\sqrt{n}\tilde{\kappa}-\frac{1}{\sqrt{n}}}\leq \frac{\EE\sbrac{Y_2^{(n)}(\infty)}}{\sqrt{n}\tilde{\kappa}-\frac{1}{\sqrt{n}}}\leq \frac{C_3}{\tilde{\kappa}\sqrt{n}},
\end{equation}
for all sufficiently large $n$. Hence,
\[
\EE[\varepsilon_{\tilde{f}^*}^n(\mf Q(\infty))]\leq \frac{4C_9C_3}{n\tilde{\kappa}} + \frac{\mu_1C_3 C_{10}(a)}{\beta \tilde{\kappa} n^{a}}
\]
and the result follows. 
\qed
\endproof}

{\subsection{Construction of the Lyapunov Function for Lemma~\ref{lem:solution_PDE}}
\label{sec:lyapunov_func}

In this section, we construct a candidate Lyapunov function $f^*(\mf x)$ which solves the PDE~\eqref{eq:PDE_Lyapunov1}-\eqref{eq:PDE_Lyapunov2} using the drift-based fluid limit (DFL) approach of~\cite{Stol2015}. This approach relates to the method of characteristics for solving the PDE by reducing the problem to solving a set of ODEs (cf. e.g. Chapter 2 in~\cite{renardy2004introduction}). The key idea behind this approach is that a solution to the PDE 
\begin{equation}
\label{eqn:dfl}
 \mathcal{L} f^{(h)}(\mf x) =-h(\mf x),   
\end{equation}
where the operator $\mathcal{L}$ satisfies $\mathcal{L} f(\mf{x})=\nabla f(\mf x) \cdot \mf F(\mf x)$ for some {\em continuous} vector field $\mf F$ is given by 
\begin{equation}
    f^{(h)}(\mf x)=\int_0^\infty h(\mf v^{\mf x} (s))ds \label{eq:candidate},
\end{equation}
provided that the right-hand side of~\eqref{eq:candidate} is finite and the function $\mf v^{\mf x}(t)$ is the solution of the  ODE  $\dot{\mf{v}}(t)=\mf F(\mf v(t))$ with initial condition $\mf v(0)=\mf x$.

To apply the DFL approach to our problem, observe that a function $f^*$ solves the  PDE~\eqref{eq:PDE_Lyapunov1}-\eqref{eq:PDE_Lyapunov2}
if and only if it is also the solution of the following PDE
\begin{align}\label{eq:PDE_Lyapunov12}
    \mathcal{\bar{L}} f^*(\mf x) &= -\brac{x_2-\frac{\kappa}{\sqrt{n}}}_+, \\
    \label{eq:PDE_Lyapunov22}f_1^*(0,x_2) &= f_2^*(0,x_2)=\frac{\sqrt{n}}{\beta}\left(x_2-\frac{\kappa}{\sqrt{n}}\right)_+,
\end{align}
for all $\mf x\in \RR^2$ with $x_2\geq 0$, where the operator $\mathcal{\bar{L}}$ is defined as
\begin{align}\label{eq:PDE_operator_extended}
    \mathcal{\bar{L}} f(\mf x) &=\left(-\frac{\beta}{\sqrt{n}}-\mu_M x_1 + \mu_1 x_2\right)\indic{\mf x\in \Omega\setminus\bar{\Omega}}f_1(\mf x) \nonumber\\
    &- \mu_1 x_2 \indic{\mf x\in \Omega\setminus\bar{\Omega}}f_2(\mf x)-\frac{\beta}{\sqrt{n}}\indic{\mf x \in \bar \Omega}f_2(\mf x)
\end{align}
and $\bar \Omega=\left\{\mf x\in \Omega: x_1=0, x_2 >\frac{\beta}{\mu_1\sqrt{n}}\right\}$.
Clearly, the operator $\mathcal{\bar{L}}$ satisfies $\mathcal{\bar L} f(\mf{x})=\nabla f(\mf x) \cdot \mf F(\mf x)$
with $\mf F$ is defined as  
\begin{equation}
\mf F(\mf x)=\begin{cases}
                \left(0,-\frac{\beta}{\sqrt{n}}\right), & \mf x\in \bar \Omega\\
                \brac{-\frac{\beta}{\sqrt{n}} - \mu_M x_1 +\mu_1 x_2,-\mu_1x_2}, & \mf x\in \Omega\setminus\bar{\Omega}
            \end{cases}
            \label{eq:vector_field}
\end{equation}
Let $\mf v^{\mf x}(t)$ denote the solution to the ODE
\begin{align}
\dot{\mf v}(t)&=\mf F(\mf v(t)),\label{eq:ode1}\\
\mf v(0)&=\mf x\in \Omega\label{eq:ode2}.
\end{align}
Thus, taking $h(\mf x)=(x_2-\kappa/\sqrt{n})_+$, the DFL approach implies that 
\begin{align}\label{eq:lyapunov_integral}
    f^*(\mf x)=\int_0^\infty (v_2^{\mf x}(s)-\kappa/\sqrt{n})_+ds,
\end{align}
is a candidate for the solution to the PDE~\eqref{eq:PDE_Lyapunov12}-\eqref{eq:PDE_Lyapunov22}. In the remainder of this section, we derive the expression of $f^*$ given in~\eqref{eq:lyapunov_integral} by solving~\eqref{eq:ode1}-\eqref{eq:ode2} and computing the integral in~\eqref{eq:lyapunov_integral} explicitly.
However, it is to be noted that $f^*$ constructed in this way is not guaranteed to solve the PDE~\eqref{eq:PDE_Lyapunov12}-\eqref{eq:PDE_Lyapunov22} since the vector field $\mf F$ in~\eqref{eq:vector_field} is {\em not continuous}. Hence, we need to additionally verify that the constructed $f^*$ indeed solves~\eqref{eq:PDE_Lyapunov12}-\eqref{eq:PDE_Lyapunov22}.

\begin{remark}
We note that 
if the two-dimensional process $\mf v^{\mf x}(t)=(v_1^{\mf x}(t),v_2^{\mf x}(t))$ solves~\eqref{eq:ode1}-\eqref{eq:ode2} for each $\mf x=(x_1, x_2) \in \Omega$ then it must also solve the integral equations
 \begin{align}
    \label{eq:fluid_v1} &v_1^{\mf x}(t) = x_1 -\frac{\beta}{\sqrt{n}}t - \int_0^t(\mu_M v_1^{\mf x}(s) -\mu_1 v_2^{\mf x}(s))ds -U_1(t),\\
    \label{eq:fluid_v2} &v_2^{\mf x}(t) = x_2 -\int_0^t \mu_1 v_2^{\mf x}(s)ds + U_1(t),
\end{align}
where the function $U_1(t)$ is the unique non-decreasing function satisfying
\begin{align}
    \label{eq:fluid_reflection} &\int_0^\infty v_1^{\mf x}(s)dU_1(s) = 0, U_1(t)\geq 0, U_1(0)=0.
\end{align}
Proposition~\ref{thm:skorohod_mapping} ensures that the process $\mf v^{\mf x}(t)=(v_1^{\mf x}(t),v_2^{\mf x}(t)) \in \Omega$ is well-defined. 
Note that the equations~\eqref{eq:fluid_v1}-\eqref{eq:fluid_reflection} closely resemble those used to define the fluid model in \cite{Braverman2020}, except for the coefficients of $v_1$ and $v_2$. While this may appear to be a subtle difference, it significantly alters the solution, requiring a full re-derivation of the solution.   
\end{remark}

The drift of the process $\mf v^{\mf x}$, given by the function $\mf F$, clearly has a switching behaviour depending on whether the process $\mf v^{\mf x}$ enters the set $\bar \Omega$ or not which in turn is determined by the starting point $\mf x$. 
More specifically, if we define $\tau(\mf x) =\inf\{t\geq 0: v_1^{\mf x}(t)=0\}$ and set $\tau(\mf x)=\infty$ if $v_1^{\mf x}(t)<0$ for all $t\geq 0$, then 
$\tau(\mf x)=\infty$ implies that the process never enters the set $\bar \Omega$.
Furthermore, for all $t\in [0,\tau(\mf x))$, we have
%
\begin{align}
\label{eqn:s_1a}
    \frac{dv_1^{\mf x}(t)}{dt} &= -\frac{\beta}{\sqrt{n}} - \mu_M v_1^{\mf x}(t) +\mu_1 v_2^{\mf x}(t)\\
    \frac{dv_2^{\mf x}(t)}{dt} &= -\mu_1 v_2^{\mf x}(t). \label{eqn:s_1b}
\end{align}
Solving the second equation, we obtain
\begin{align}
    v_2^{\mf x}(t) = x_2 e^{-\mu_1t}
    \label{eq:sol_v2}
\end{align}
Plugging this into the first equation, we obtain
\begin{align}
    \label{eqn:v_1_Omega3a}
    v_1^{\mf x}(t) &= -\frac{\beta}{\mu_M \sqrt{n}}- \frac{\mu_1}{\mu_1-\mu_M}x_2e^{-\mu_1t} + \left(x_1 + \frac{\beta}{\mu_M\sqrt{n}} +\frac{\mu_1x_2}{\mu_1-\mu_M}\right)e^{-\mu_M t},\\
    \label{eqn:v_1_Omega3b}
    &=-\frac{\beta}{\mu_M \sqrt{n}}+ \left(x_1 + \frac{\beta}{\mu_M\sqrt{n}}\right)e^{-\mu_M t} +\frac{\mu_1x_2}{\mu_1-\mu_M}\left(e^{-\mu_Mt}-e^{-\mu_1t}\right).
\end{align}

In the lemmas below, we determine the conditions needed on the initial point $\mf x$ for $\tau(\mf x)$ to be finite and show that if $\tau(\mf x)$ is finite then $v_2^{\mf x}(\tau(\mf x))\geq \beta/\mu_1\sqrt{n}$ implying that the trajectory $\mf v^{\mf x}$ indeed enters the set $\bar \Omega$.

\begin{lemma}\label{lem:tau_solution}
Consider the nonlinear system of equations
\begin{align}
    \label{eq:tau_equation}&-\frac{\beta}{\mu_M \sqrt{n}} + \left(x_1 + \frac{\beta}{\mu_M\sqrt{n}}\right)e^{-\mu_M \tau} -\frac{\mu_1x_2}{\mu_1-\mu_M}\left(e^{-\mu_1\tau}-e^{-\mu_M\tau}\right)=0,\\
    &x_2e^{-\mu_1\tau} ={\kappa}/\sqrt{n}.\label{eq:tau_equation1}
\end{align}
The above system of equations has a unique solution in $(x_2,\tau) \in \left[{\kappa}/\sqrt{n}, \infty\right)\times [0,\infty)$ denoted by  $(x_2^*(x_1,{\kappa}),\tau^*(x_1,{\kappa}))$ for any $x_1\leq 0$ and any ${\kappa} \geq \beta/\mu_1$.    
In particular, we have $(x_2^*(0,{\kappa}),\tau^*(0,{\kappa}))=({\kappa}/\sqrt{n},0)$. 
\end{lemma}

Lemma~\ref{lem:tau_solution}, proved in Appendix~\ref{proof:tau_solution}, allows us to define the curve 
\begin{equation}
    \Gamma^{{\kappa}}=\{\mf x \in \Omega: x_2=x_2^*(x_1,{\kappa})\}
    \label{eq:gamma_kappa}
\end{equation}
for each ${\kappa} \geq \beta/\mu_1$. Clearly, we have $(0,\kappa/\sqrt{n}) \in \Gamma^{\kappa}$. 
We write $\mf x \geq \Gamma^{\kappa}$ when $x_2\geq x_2^*(x_1,\kappa)$ and define 
$\mf x > \Gamma^\kappa$, $\mf x \leq \Gamma^{\kappa}$ similarly. The following lemma relates the curves $\Gamma^{\kappa}$, $\kappa \geq \beta/\sqrt{n}$, to $\tau(\mf x)$.

 \begin{lemma}
\label{lem:tau_def_lem}
The following statements hold:
\begin{enumerate}
    \item If $\mf x \geq \Gamma^{{\kappa}}$ for some ${\kappa} \geq \beta/\mu_1$, then $\tau(\mf x) < \infty$ and $x_2e^{-\mu_1\tau(\mf x)} \geq {\kappa}/\sqrt{n}$. 
    
    \item If $\mf x\in \Gamma^{{\kappa}}$ for some ${\kappa} \geq \beta/\mu_1$, then $\tau(\mf x)=\tau^*(x_1,{\kappa})$, where $\tau^*(\cdot)$ is as defined in Lemma~\ref{lem:tau_solution}. Hence, $x_2e^{-\mu_1\tau(\mf x)}={\kappa}/\sqrt{n}$.

    \item We have that $\tau(\mf x)$ is differentiable if $\mf x \geq \Gamma^{{\kappa}}$ for some ${\kappa} \geq \beta/\mu_1$. Furthermore, the derivatives $\tau_1(\mf x)=\partial \tau(\mf x)/\partial x_1$ and $\tau_2(\mf x)=\partial \tau(\mf x)/\partial x_2$ satisfy
    \begin{align}
        \tau_1(\mf x) &=-e^{-\mu_M\tau(\mf x)}\left(\mu_1x_2e^{-\mu_1\tau(\mf x)}-\frac{\beta}{\sqrt{n}}\right)^{-1}\leq 0,\label{eq:tau1}\\
        \tau_2(\mf x) &=\frac{\mu_1}{\mu_1-\mu_M}\left(e^{-\mu_1\tau(\mf x)}-e^{-\mu_M\tau(\mf x)}\right)\left(\mu_1x_2e^{-\mu_1\tau(\mf x)}-\frac{\beta}{\sqrt{n}}\right)^{-1} \leq 0,\label{eq:tau2}
    \end{align}
    where $\tau_1(\mf x)$ is understood to be the left-derivative when $x_1=0$ and $\mf x \geq \Gamma^{{\kappa}}$ for some ${\kappa} \geq \beta/\mu_1$.
    
    \item If ${\kappa}_1 > {\kappa}_2 \geq \beta/\mu_1$, then $\Gamma^{{\kappa}_1} > \Gamma^{{\kappa}_2}$, i.e., $\mf x \geq \Gamma^{{\kappa}_1}$ implies $\mf x > \Gamma^{{\kappa}_2}$.

    \item If $\mf x < \Gamma^{\beta/\mu_1}$, then $\tau(\mf x)=\infty$.
\end{enumerate}


\end{lemma}

Thus, from Lemma~\ref{lem:tau_def_lem}, it follows that $\tau(\mf x) < \infty$ if $\mf x \in \Gamma^{\kappa}$, for some $\kappa\geq \beta/\mu_1$, or equivalently, if $\mf x \geq \Gamma^{\beta/\mu_1}$ and in this case
$v_2^{\mf x}(\tau(\mf x))=x_2 e^{-\mu_1\tau(\mf x)}\geq \beta/\mu_1\sqrt{n}$, i.e., at time $\tau(\mf x)$ the trajectory $\mf v^{\mf x}$ enters the set $\bar \Omega$. On the other hand if $\mf x < \Gamma^{\beta/\mu_1}$, then the trajectory $\mf v^{\mf x}$ never enters the set $\bar \Omega$. This behaviour is shown in Figure~\ref{fig:trag}. Once the trajectory $\mf v^{\mf x}$ enters $\bar \Omega$, we have from~\eqref{eq:vector_field}-\eqref{eq:ode2} that
\begin{align}
\label{eqn:s_2a}
    dv_1^{\mf x}(t) &= 0\\
    dv_2^{\mf x}(t) &= -\frac{\beta}{\sqrt{n}}dt,\label{eqn:s_2b}
\end{align}
The above equations hold as long as $\mf v^{\mf x}$ remains within $\bar \Omega$. Hence, for any $\kappa\geq\beta/\mu_1$, if we define $\tau^{\kappa}(\mf x)=\inf\{t\geq 0: v_2^{\mf x}(t)=\kappa/\sqrt{n}\}$, then 
$\tau^{\beta/\mu_1}(\mf x)\geq \tau(\mf x)$ if 
$\mf x\geq \Gamma^{\beta/\mu_1}$. Furthermore, for $t\in [\tau(\mf x),\tau^{\beta/\mu_1}(\mf x))$ we can characterise  $v_2^{\mf x}(t)$ using \eqref{eqn:s_2b} as
\begin{align}
    v_2^{\mf x}(t)=-\frac{\beta}{\sqrt{n}}(t-\tau(\mf x)) +x_2e^{-\mu_1\tau(\mf x)}.
\end{align}
After leaving the set $\bar \Omega$ at time $\tau^{\beta/\mu_1}(\mf x)$ the trajectory of $\mf v^{\mf x}$ never re-enters the set $\bar \Omega$  since $v_2^{\mf x}()$ is non-increasing in time.
From Lemma~\ref{lem:tau_def_lem}, it also follows that if $\mf x <\Gamma^{\beta/\mu_1}$, then $\tau(\mf x)=\infty$ and therefore $v_1^{\mf x}(t)$ and $v_2^{\mf x}(t)$ are given by~\eqref{eqn:v_1_Omega3b} and~\eqref{eq:sol_v2}, respectively, for all $t\geq 0$. 

\begin{figure}
    \centering
\begin{tikzpicture}
    \begin{scope}[shift={(-5,0)}]
    \draw[thick,<->] (-3,0) -- (3,0) node[anchor=north] {$x_1$};
    \draw[thick,<->] (0,-3) -- (0,3) node[anchor=west] {$x_2$};
        
        \filldraw[red] (-1,0) circle (2pt);
        \node at (-1.3,-0.3) {$(-\frac{\beta}{\mu_M\sqrt{n}},0)$};
        
        \draw[dashed] (-3,3) .. controls (0,2) and (0,1) .. (0,1) node[below right]{$(0,\frac{\beta}{ \mu_1 \sqrt{n}})$};
        \filldraw[red] (0,1) circle (2pt);

    \draw[blue,->] (-1.5,0.3) -- (-1.3,0.1);
    \draw[blue,->] (-1,0.3) -- (-1,0.1);
    \draw[blue,->] (-0.5,0.3) -- (-0.7,0.1);

    \draw[blue,->] (-1.5,1) -- (-1.1,0.5);
    \draw[blue,->] (-1,1) -- (-0.8,0.6);
    \draw[blue,->] (-0.5,0.8) -- (-0.7,0.4);

    \draw[blue,->] (-2,1) -- (-1.8,0.97);
    \draw[blue,->] (-2.4,1) -- (-2.1,0.95);
    \draw[blue,->] (-2.4,0) -- (-2.1,0);

    \draw[blue,->] (-1,2) -- (-0.5,1.5);
    \draw[blue,->] (-0.5,1.4) -- (-0.4,1.2);
    \end{scope}

    \begin{scope}[shift={(5,0)}]
    \draw[thick,<->] (-3,0) -- (3,0) node[anchor=north] {$x_1$};
    \draw[thick,<->] (0,-3) -- (0,3) node[anchor=west] {$x_2$};
        \filldraw[red] (-1,0) circle (2pt);
        \node at (-1.3,-0.3) {($-\frac{\beta}{\mu_M\sqrt{n}},0)$};

        \draw[dashed] (-3,3) .. controls (0,2) and (0,1) .. (0,1) node[below right]{$(0,\frac{\beta}{ \mu_1 \sqrt{n}})$};
         \filldraw[red] (0,1) circle (2pt);
        
        \draw[red,thick,->] (0,2) -- (0,1);

    \draw[blue,->] (-0.4,2.3) -- (-0.1,2);
     \draw[blue,->] (-0.8,2.4) -- (-0.5,2.2);
    \draw[blue,->] (0,1) -- (-0.6,0.5);

    \end{scope}

\end{tikzpicture}
    \caption{Any trajectory that begins below the dashed curve representing the curve $\Gamma^{\beta/\mu_1}$ will not intersect the vertical axis, while trajectories starting above the curve will eventually reach the axis and move downward until arriving at the point $(0, \beta/\mu_1\sqrt{n})$.}
    \label{fig:trag}
\end{figure}
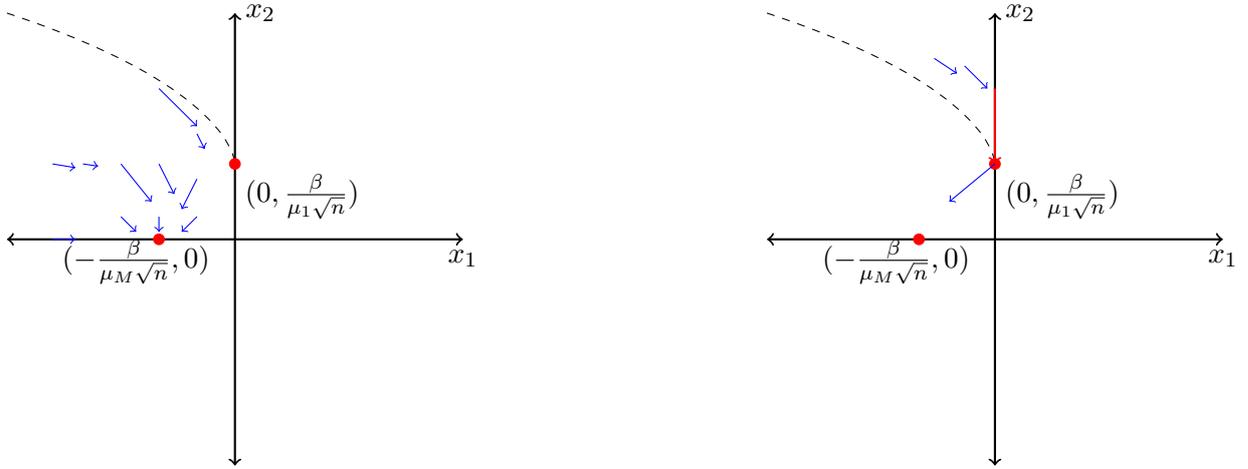
 
Taking $\kappa$ as in Lemma~\ref{lem:solution_PDE}, equation~\eqref{eq:tau_equation1} implies that the set $\{\mf x \in \Omega:x_2 \leq \kappa/\sqrt{n}, \mf x > \Gamma^\kappa\}$ 
is empty. Hence, the set $\Omega$ can be partitioned as
$\Omega=\Omega_1\cup\Omega_2\cup \Omega_3$, where
\[
\Omega_1=\{x_2\leq \kappa/\sqrt{n}\}, \Omega_2 = \{x_2>\kappa/\sqrt{n}\mbox{ and }\mf x\leq \Gamma^\kappa\}\mbox{ and }\Omega_3=\{\mf x>\Gamma^\kappa\}.
\]}

We are now ready to evaluate~\eqref{eq:lyapunov_integral} for each $\mf x \in \Omega_i$ for $i\in[3]$.
Since $v_2^{\mf x}$ is non-increasing in time we have $v_2^{\mf x}(s) \leq v_2^{\mf x}(0)=x_2\leq \kappa/\sqrt{n}$ for all $s\geq0$ and for all $\mf x\in \Omega_1$. Hence, from~\eqref{eq:lyapunov_integral} we have $f^*(\mf x)=0$ for all $\mf x\in \Omega_1$. 

For $\mf x\in \Omega_2$, we have 
$\tau^\kappa(\mf x)\leq \tau(\mf x)$ because otherwise we would have $x_2e^{-\mu_1\tau(\mf x)}> \kappa/\sqrt{n}$ which contradicts the fact that $\mf x \leq \Gamma^{\kappa}$ according to Lemma~\ref{lem:tau_def_lem}. Since $v_2^{\mf x}$ is given by~\eqref{eq:sol_v2} for all $t\in [0,\tau(\mf x))$ we have
$\tau^\kappa(\mf x)=\mu_1^{-1}\log(\sqrt{n}x_2/\kappa)$ and therefore
\begin{align*}
    f^*(\mf x)&=\int_0^{\frac{1}{\mu_1}\log\left(\frac{\sqrt{n}x_2}{\kappa}\right)} \left(x_2e^{-\mu_1t}-\frac{\kappa}{\sqrt{n}}\right)dt\\
    &=\left.-\frac{x_2}{\mu_1}e^{-\mu_1t}- \frac{\kappa}{\sqrt{n}}t \right|_0^{\frac{1}{\mu_1}\log\left(\frac{\sqrt{n}x_2}{\kappa}\right)}\\
    &=\frac{x_2}{\mu_1}-\frac{\kappa}{\mu_1\sqrt{n}}-\frac{\kappa}{\mu_1\sqrt{n}}\log\left(\frac{\sqrt{n}x_2}{\kappa}\right).
\end{align*}

Finally, for $\mf x\in \Omega_3$, we have $\tau^{\kappa}(\mf x) \geq \tau(\mf x)$ and therefore
\begin{align*}
    f^*(\mf x)&=\int_0^{\tau(\mf x)}\left(x_2e^{-\mu_1t}-\frac{\kappa}{\sqrt{n}}\right)dt + \int_0^{-\frac{\sqrt{n}}{\beta}\left(\frac{\kappa}{\sqrt{n}} + x_2e^{-\mu_1\tau(\mf x)}\right)}\left(-\frac{\beta}{\sqrt{n}}t +x_2e^{-\mu_1\tau(\mf x)}-\frac{\kappa}{\sqrt{n}}\right)dt\\
    &=\left.-\frac{x_2}{\mu_1}e^{-\mu_1t}- \frac{\kappa}{\sqrt{n}}t \right|_0^{\tau(\mf x)} - \left.\frac{\beta}{\sqrt{n}}t^2 +\left(x_2e^{-\mu_1\tau(\mf x)}-\frac{\kappa}{\sqrt{n}}\right)t\right|_0^{-\frac{\sqrt{n}}{\beta}\left(\frac{\kappa}{\sqrt{n}} + x_2e^{-\mu_1\tau(\mf x)}\right)}\\
    &=\frac{x_2}{\mu_1}-\frac{x_2}{\mu_1}e^{-\mu_1\tau(\mf x)}-\frac{\kappa}{\sqrt{n}}\tau(\mf x)+\frac{\beta}{2\sqrt{n}}\left(\frac{\sqrt{n}}{\beta}x_2e^{-\mu_1\tau(\mf x)}-\frac{\kappa}{\beta}\right)^2\\
    &=\frac{x_2}{\mu_1}-\frac{x_2}{\mu_1}e^{-\mu_1\tau(\mf x)}-\frac{\kappa}{\sqrt{n}}\tau(\mf x)+\frac{\sqrt{n}}{2\beta}\left(x_2e^{-\mu_1\tau(\mf x)}-\frac{\kappa}{\sqrt{n}}\right)^2.    
\end{align*}

Aggregating the above expressions, we obtain the candidate solution $f^*(\mf x)$ to the PDE~\eqref{eq:PDE_Lyapunov1}-\eqref{eq:PDE_Lyapunov2} as 
\begin{align}
    \label{eq:lyapunov_function}
    f^*(\mf x) &=\left\{\begin{array}{ll}
    0, &\mf x\in \Omega_1,\\
    \frac{x_2}{\mu_1}-\frac{\kappa}{\mu_1\sqrt{n}}-\frac{\kappa}{\mu_1\sqrt{n}}\log\left(\frac{\sqrt{n}x_2}{\kappa}\right), &\mf x\in \Omega_2,\\
    \frac{x_2}{\mu_1}-\frac{x_2}{\mu_1}e^{-\mu_1\tau(\mf x)}-\frac{\kappa}{\sqrt{n}}\tau(\mf x)+\frac{\sqrt{n}}{2\beta}\left(x_2e^{-\mu_1\tau(\mf x)}-\frac{\kappa}{\sqrt{n}}\right)^2, &\mf x\in \Omega_3.
    \end{array}\right.
\end{align}
In Appendix~\ref{append:lem9} we complete the proof of Lemma~\ref{lem:solution_PDE} by verifying that the above function indeed satisfies~\eqref{eq:PDE_Lyapunov1}-\eqref{eq:PDE_Lyapunov2}.

{ \section{Asymptotic optimality of the SA-JSQ scheme: Proof of Theorem~\ref{thm:sajsq_optimality} }
\label{sec:mod_sys}

In this section, we establish the asymptotic optimality of the SA-JSQ scheme in the Halfin-Whitt regime. To demonstrate this optimality, we first derive a lower bound on the performance achievable by any admissible load balancing scheme for the original system. This lower bound is obtained by coupling the original system with a modified system in which the servers are not attached to the queues but can freely move between queues ensuring that faster servers always handle the longer queues. We then show that this modified system shares the same diffusion limit as the original system when operating under the SA-JSQ scheme.

\subsection{A Modified System with Free Servers and Blocking}

To identify a lower bound on the achievable performance of admissible load balancing policies, we consider a {\em modified system} where, similar to the original system, in each pool $j \in [M]$, there are $N_j^{(n)}$ servers with rate $\mu_j$ and $\sum_{j=1}^M N_j^{(n)}=n$  queues. However, contrary to the original system, the queues are not attached to the servers. Instead, at any time $t\geq 0$, the servers re-arrange themselves in a preemptive fashion so that longer queues are always served by faster servers. We further assume that each queue has a buffer capacity of two including the job in service at that queue and arrivals occurring when all queues have two jobs are {\em rejected/blocked}. Each arriving job is irrevocably assigned to a queue according to the conventional JSQ policy which picks the shortest queue (breaking ties uniformly at random) for each incoming arrival. 

To better understand the dynamics of the modified system, consider a system where there are two type 1 servers with rate $\mu_1$, two type 2 servers with rate $\mu_2$ and one type 3 server with rate $\mu_3$. Suppose that right before time $t$ (i.e., at time $t-$), there are three queues with 2 jobs and 2 queues with 1 jobs as illustrated in Figure~\ref{fg:modified_before}, where at each queue the first job at the head-of-the-line (shaded job) is being served. As the servers always need to organize themselves so that the faster servers serve the longer queues, we need to have  two type 1 servers and one type two server serving queues with 2 jobs and one type 2 server and one type 3 server serving queues with one job at $t-$. Now, consider the following two scenarios. For the first scenario, suppose that a departure occurs at $t$ from the second queue with two jobs being served by a type 1 server. Then, to ensure that the longer queues are served by faster servers, the type 1 server serving this queue will change places with the type 2 server serving the queue with two jobs as illustrated in Figure~\ref{fg:modified_after_departure}. As another scenario, suppose that an arrival occurs at time $t$ and following the JSQ policy it is assigned to the fifth queue which was being served by a type 3 server. Then, at time $t$ the type 2 server serving the queue with one job should change places with the type 3 server as illustrated in Figure~\ref{fg:modified_after_arrival}. 

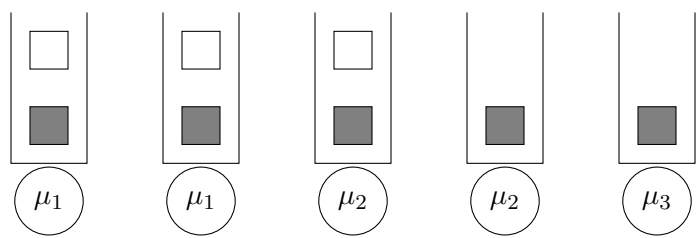
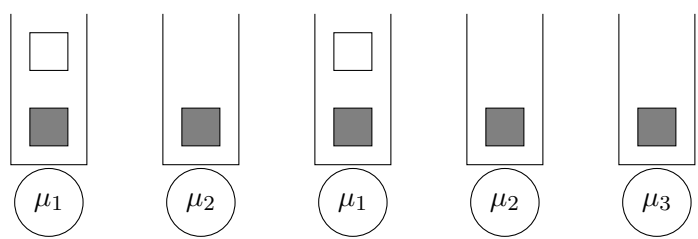
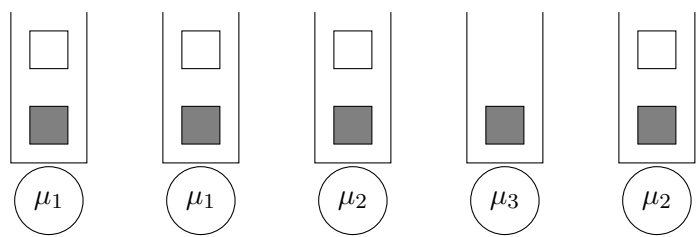
\begin{figure}[htbp]
    \centering
    
\subfigure[%
Snapshot of the modified system at $t-$. %
        \label{fg:modified_before}]{
        \centering
        \begin{tikzpicture}
            \draw (0,3) -- (0,1) -- (1,1) -- (1,3);
            \draw [fill=gray] (0.25,1.25) rectangle (0.75,1.75);
            \draw (0.25,2.25) rectangle (0.75,2.75);
            \draw (0.5,0.5) circle [radius = 0.45];
            \draw (0.5,0.5) node [scale = 1]{$\mu_1$};

            \draw (2,3) -- (2,1) -- (3,1) -- (3,3);
            \draw [fill=gray] (2.25,1.25) rectangle (2.75,1.75);
            \draw (2.25,2.25) rectangle (2.75,2.75);
            \draw (2.5,0.5) circle [radius = 0.45];
            \draw (2.5,0.5) node [scale = 1]{$\mu_1$};

            \draw (4,3) -- (4,1) -- (5,1) -- (5,3);
            \draw [fill=gray] (4.25,1.25) rectangle (4.75,1.75);
            \draw (4.25,2.25) rectangle (4.75,2.75);
            \draw (4.5,0.5) circle [radius = 0.45];
            \draw (4.5,0.5) node [scale = 1]{$\mu_2$};

            \draw (6,3) -- (6,1) -- (7,1) -- (7,3);
            \draw [fill=gray] (6.25,1.25) rectangle (6.75,1.75);
            \draw (6.5,0.5) circle [radius = 0.45];
            \draw (6.5,0.5) node [scale = 1]{$\mu_2$};

            \draw (8,3) -- (8,1) -- (9,1) -- (9,3);
            \draw [fill=gray] (8.25,1.25) rectangle (8.75,1.75);
            \draw (8.5,0.5) circle [radius = 0.45];
            \draw (8.5,0.5) node [scale = 1]{$\mu_3$};
\end{tikzpicture}
}

    \hfill
    
\subfigure[%
Scenario 1: A departure from the second queue occurs at time $t$. Servers serving second and third queues are swapped.%
\label{fg:modified_after_departure}]{
        \centering
                \begin{tikzpicture}
            \draw (0,3) -- (0,1) -- (1,1) -- (1,3);
            \draw [fill=gray] (0.25,1.25) rectangle (0.75,1.75);
            \draw (0.25,2.25) rectangle (0.75,2.75);
            \draw (0.5,0.5) circle [radius = 0.45];
            \draw (0.5,0.5) node [scale = 1]{$\mu_1$};

            \draw (2,3) -- (2,1) -- (3,1) -- (3,3);
            \draw [fill=gray] (2.25,1.25) rectangle (2.75,1.75);
            \draw (2.5,0.5) circle [radius = 0.45];
            \draw (2.5,0.5) node [scale = 1]{$\mu_2$};

            \draw (4,3) -- (4,1) -- (5,1) -- (5,3);
            \draw [fill=gray] (4.25,1.25) rectangle (4.75,1.75);
            \draw (4.25,2.25) rectangle (4.75,2.75);
            \draw (4.5,0.5) circle [radius = 0.45];
            \draw (4.5,0.5) node [scale = 1]{$\mu_1$};

            \draw (6,3) -- (6,1) -- (7,1) -- (7,3);
            \draw [fill=gray] (6.25,1.25) rectangle (6.75,1.75);
            \draw (6.5,0.5) circle [radius = 0.45];
            \draw (6.5,0.5) node [scale = 1]{$\mu_2$};

            \draw (8,3) -- (8,1) -- (9,1) -- (9,3);
            \draw [fill=gray] (8.25,1.25) rectangle (8.75,1.75);
            \draw (8.5,0.5) circle [radius = 0.45];
            \draw (8.5,0.5) node [scale = 1]{$\mu_3$};
\end{tikzpicture}}
        
    \hfill
    
\subfigure[%
Scenario 2: An arrival occurs at time $t$. Servers serving the fourth and the fifth queues are swapped.%
\label{fg:modified_after_arrival}]{
        \centering
        \begin{tikzpicture}
            \draw (0,3) -- (0,1) -- (1,1) -- (1,3);
            \draw [fill=gray] (0.25,1.25) rectangle (0.75,1.75);
            \draw (0.25,2.25) rectangle (0.75,2.75);
            \draw (0.5,0.5) circle [radius = 0.45];
            \draw (0.5,0.5) node [scale = 1]{$\mu_1$};

            \draw (2,3) -- (2,1) -- (3,1) -- (3,3);
            \draw [fill=gray] (2.25,1.25) rectangle (2.75,1.75);
            \draw (2.25,2.25) rectangle (2.75,2.75);
            \draw (2.5,0.5) circle [radius = 0.45];
            \draw (2.5,0.5) node [scale = 1]{$\mu_1$};

            \draw (4,3) -- (4,1) -- (5,1) -- (5,3);
            \draw [fill=gray] (4.25,1.25) rectangle (4.75,1.75);
            \draw (4.25,2.25) rectangle (4.75,2.75);
            \draw (4.5,0.5) circle [radius = 0.45];
            \draw (4.5,0.5) node [scale = 1]{$\mu_2$};

            \draw (6,3) -- (6,1) -- (7,1) -- (7,3);
            \draw [fill=gray] (6.25,1.25) rectangle (6.75,1.75);
            \draw (6.5,0.5) circle [radius = 0.45];
            \draw (6.5,0.5) node [scale = 1]{$\mu_3$};

            \draw (8,3) -- (8,1) -- (9,1) -- (9,3);
            \draw [fill=gray] (8.25,1.25) rectangle (8.75,1.75);
            \draw (8.25,2.25) rectangle (8.75,2.75);
            \draw (8.5,0.5) circle [radius = 0.45];
            \draw (8.5,0.5) node [scale = 1]{$\mu_2$};
        \end{tikzpicture}}
        
    
    \caption{A modified system with 2 type 1 servers, two type 2 servers and one type 3 server.}
\end{figure}

Similar to the original system, we define the following quantities for the modified system:
\begin{align*}
    \tilde{Q}_i^{(n)}(t)&:=\mbox{The number of queues with $i$ or more jobs at time $t$}\\
    \tilde{A}^{(n)}(t) &:=\mbox{The cumulative number of external arrivals by time $t$}\\
    \tilde{A}_1^{(n)}(t)&:= \mbox{The number of arrivals routed to queues with exactly 1 job by time $t$}\\
    \tilde{A}_R^{(n)}(t)&:=\mbox{The number of rejected jobs by time $t$}\\
    \tilde{D}_i^{(n)}(t) &:=\mbox{The total number of departures from servers with $i$ or more jobs by time $t$.}
\end{align*}
As the arrivals are routed in a JSQ fashion where the ties are broken randomly, we have
\begin{align*}
\tilde{A}_1^{(n)}(t) &:= \int_0^\infty \1\left(\tilde{Q}_1^{(n)}(s-)=\sum_{j=1}^MN_j^{(n)} \mbox{ and }\tilde{Q}_2^{(n)}(s-)<\sum_{j=1}^MN_j^{(n)}\right)d\tilde{A}^{(n)}(s)\\
\tilde{A}_R^{(n)}(t) &:= \int_0^\infty \1\left(\tilde{Q}_2^{(n)}(s-)=\sum_{j=1}^MN_j^{(n)}\right)d\tilde{A}^{(n)}(s)\\
\end{align*}
Then, the number of arrivals routed to idle servers is given by $\tilde{A}^{(n)}(t)-\tilde{A}_R^{(n)}(t)-\tilde{A}_1^{(n)}(t)$ for any time $t\geq 0$. Finally, as we always have the faster servers serving the longer queues, the departure rate for queues with $i$ or more jobs at time $t$ is given by
\begin{equation}
\label{eqn:rate_modified}
\tilde{\mu}^{(n)}(\tilde{Q}_i^{(n)}(t-)):=\sum_{j=1}^{M}\mu_j\left[\left(\tilde{Q}_i^{(n)}(t-) - \sum_{k=1}^{j-1}N_k^{(n)}\right)_+\wedge N_j^{(n)}\right].    
\end{equation}
From Lemma~16 of~\cite{bhambay2022asymptotic}
it follows that for any $0\leq Q_{j,i}^{(n)}(t-)\leq N_j^{(n)}$, we have 
\begin{equation}
    \tilde{\mu}^{(n)}\left(\sum_{j=1}^{M} Q_{j,i}^{(n)}(t-)\right)\geq \sum_{j=1}^M \mu_j Q_{j,i}^{(n)}(t-)
    \label{eq:rate_bound}
\end{equation}
Now, we have the following balance equations: 
\begin{align}
\tilde{Q}_1^{(n)}(t) &= \tilde{Q}_1^{(n)}(0) + \tilde{A}^{(n)}(t)-\tilde{A}_R^{(n)}(t) - \tilde{A}_1^{(n)}(t) - \tilde{D}_1^{(n)}(t) + \tilde{D}_2^{(n)}(t)\label{eq:balance_eq_m1}\\
\tilde{Q}_{2}^{(n)}(t) &=\tilde{Q}_{2}^{(n)}(0) + \tilde{A}_1^{(n)}(t) -\tilde{D}_2^{(n)}(t) \label{eq:balance_eq_m2}.
\end{align}
To present key results for the modified system, we define the following quantity for any $t \geq 0 $:
\[
{Q}_{+2}^{(n,\Pi)}(t) = \sum_{j=1}^M \sum_{i=2}^\infty Q_{j,i}^{(n)}(t),
\]
which represents the number of servers with at least two jobs in the original system under any admissible load balancing policy $\Pi$. With this definition in place, we now introduce the following coupling lemma.

\begin{proposition}
\label{prop:coupling_modf}
    Consider the original system under any admissible load balancing policy $\Pi$ and the modified system with initial states 
    \[
    Q^{(n,\Pi)}_1(0)=\tilde{Q}_1^{(n)}(0) \mbox{ and } Q_{+2}^{(n,\Pi)}(0) = \tilde{Q}_2^{(n)}(0) \leq \sum_{j=1}^M N_j^{(n)}.
    \]
    Then, there exists a coupling between the two systems such that for any $t\geq 0$, we have
    $\tilde{Q}_1^{(n)}(t) + \tilde{Q}_2^{(n)}(t)\leq Q_1^{(n,\Pi)}(t) + Q_{+2}^{(n,\Pi)}(t) \mbox{ and } \tilde{Q}_2^{(n)}(t)\leq Q_{+2}^{(n,\Pi)}(t)$.
\end{proposition}

The proposition above implies that if $\Pi$ is any stable admissible policy then $\mb{P}(\tilde{Q}_1^{(n)}(\infty) + \tilde{Q}_2^{(n)}(\infty) > y) \leq \mb{P}(Q_1^{(n,\Pi)}(\infty) + Q_{+2}^{(n,\Pi)}(\infty)>y)$ for all $y \in \RR$.

We now proceed in the next lemma to prove that the diffusion limit of the modified system is identical to that of the original system operating under the SA-JSQ scheme. To state the lemma we need the following additional notations.
Let $\mathbf{\tilde{Q}}^{(n)}(t)=\Big(\tilde{Q}^{(n)}_{j,i}(t), j\in[M], i\geq1\Big)$ denote the system occupancy state at time $t\geq 0$,
where $\tilde{Q}^{(n)}_{j,i}(t)$ is the number of queues with at least $i$ jobs at time $t$ being served by of type $j$ servers. We define the scaled process $\mathbf{\tilde{Y}}^{(n)}=(\mathbf{\tilde{Y}}(t):t\geq 0)$, where $\mathbf{\tilde{Y}}(t)=\left(\tilde{Y}_{j,i}^{(n)}(t),j\in[M], i\geq1\right)$ with
\begin{equation}
\tilde{Y}_{j,1}^{(n)}(t)=\frac{\tilde{Q}_{j,1}^{(n)}(t)-N_j^{(n)}}{\sqrt{n}}, \ \   \tilde{Y}_{j,i}^{(n)}(t)=\frac{\tilde{Q}_{j,i}^{(n)}(t)}{\sqrt{n}}, \ i\geq2,j\in[M].
\end{equation}

\begin{proposition}
\label{prop:modified_diffusion_limit}
Suppose that~\eqref{eqn:heavy_traffic_limit} and~\eqref{eq:proportion_assumption} hold, and there exists a random vector $\tilde{\mf Y}(0)=\Big( \tilde{Y}_{j,i}(0), i\in[2], j\in[M]\Big) \in \MM_{M,\rho}$ for some $\rho>1$ with $\tilde{Y}_{j,1}(0)=0$ for all $j\in[M-1]$ and $\tilde{Y}_{j,2}(0)=0$ for all $j\in\cbrac{2,\dots,M}$ such that
\begin{equation}
\label{eqn:HW_SA-JSQ_initial_condition} \tilde{\mf Y}^{(n)}(0) = \tilde{\mf Y}(0), \ \forall n\geq1.
\end{equation}
Then, we have
\begin{equation*}
\tilde{\mf Y}^{(n)} \Rightarrow \tilde{\mf Y}, \  \text{as} \ n\to \infty,
\end{equation*}
where $\tilde{\mf Y} =\brac{\tilde{Y}_{j,i},i\in[2],j\in[M]}$ is the unique solution in $\DD_{\MM_{M,\rho}}[0,\infty)$ of the following  stochastic integral equations

\begin{align}
\tilde{Y}_{j,1}(t)&=0, \ j\in \cbrac{1,\dots,M-1}, \label{eq:tl_Y_j1} \\
\label{eq:tl_Y_1M} 
\tilde{Y}_{M,1}(t)&=\tilde{Y}_{M,1}(0) + \sqrt{2}W(t)-\beta t-\int_0^t (\mu_M\tilde{Y}_{M,1}(s)-\mu_1\tilde{Y}_{1,2}(s))ds -\tilde{U}_1(t),\\
\label{eq:tl_Y_12}
\tilde{Y}_{1,2}(t)&=\tilde{Y}_{1,2}(0)-\mu_1\int^t_0 \tilde{Y}_{1,2}(s)ds+\tilde{U}_1(t),\\
\tilde{Y}_{j,2}(t)&=0, \ j\in \cbrac{2,\dots,M}.\label{eq:tl_Y_j2}
\end{align}
In the above, $W$ is a standard Brownian motion and $\tilde{U}_1$ is the unique non-decreasing and non-negative process satisfying
\begin{equation}
 \label{eqn:diff_modified_fluct_term}
 \int_0^{\infty} \indic{\tilde{Y}_{M,1}(t) <0}d\tilde{U}_1(t) =0.
\end{equation}
\end{proposition}

Note that the limiting diffusion process $\tilde{\mf Y}$ given in the above theorem is the same as the process $\mf Y$ defined in Theorem~\ref{thm:SA-JSQ_diffusion_limit} started with initial condition $Y_{j,1}=0$ for all $j \in [M-1]$ and $Y_{j,2}=0$ for all $j\in \{2,3,\ldots,M\}$.
Hence, we can now prove the asymptotic optimality of the SA-JSQ scheme as follows.

\proof{Proof of Theorem~\ref{thm:sajsq_optimality}.}
The diffusion-scaled total number of jobs $\tilde Y_{+1}^{(n)}$ and the diffusion-scaled number of waiting jobs $\tilde Y_{+2}^{(n)}$ in the modified system are given by $\tilde Y_{+1}^{(n)}=(\tilde{Q}_1^{(n)}(t) + \tilde{Q}_2^{(n)}(t)-n)/\sqrt{n}$ and $\tilde Y_{+2}^{(n)}=\tilde{Q}_2^{(n)}(t)/\sqrt{n}$, respectively. The same quantities for the original system operating under any general policy $\Pi$ are given by $\tilde Y_{+1}^{(n,\Pi)}=(\tilde{Q}_1^{(n,\Pi)}(t) + \tilde{Q}_{+2}^{(n,\Pi)}(t)-n)/\sqrt{n}$ and $\tilde Y_{+2}^{(n,\Pi)}=\tilde{Q}_{+2}^{(n,\Pi)}(t)/\sqrt{n}$, respectively.
Now, by Proposition~\ref{prop:coupling_modf}, we have
$\tilde{Q}_1^{(n)}(t) + \tilde{Q}_2^{(n)}(t)\leq_{\text{st}} Q_1^{(n,\Pi)}(t) + Q_{+2}^{(n,\Pi)}(t)$ for any $t\geq0$.
Therefore, assuming that the steady-state distribution exists under policy $\Pi$, we obtain after centring and scaling

\begin{equation*}
\mathbb{P}\brac{\tilde{Y}^{(n)}_{+1}(\infty)>y}\leq \mathbb{P}\brac{Y^{(n,\Pi)}_{+1}(\infty)>y}, \forall y\in \RR.
\end{equation*}
Hence, we have 
\begin{equation*}
\lim_{n\to \infty}\mathbb{P}\brac{\tilde{Y}^{(n)}_{+1}(\infty)>y}=\lim_{n\to \infty}\mathbb{P}\brac{Y^{(n)}_{+1}(\infty)>y}\leq \liminf_{n\to \infty}\mathbb{P}\brac{Y^{(n,\Pi)}_{+1}(\infty)>y}, \forall y\in \RR,
\end{equation*}
where the equality follows from Proposition~\ref{prop:modified_diffusion_limit}. This proves~\eqref{eqn:asmp_1}. The proof of~\eqref{eqn:asmp_2} follows similarly using stochastic ordering of $\tilde{Q}_2^{(n)}(t)\leq Q_{+2}^{(n,\Pi)}(t)$ for ant $t\geq0$.
\endproof

\subsection{A Coupling Construction: Proof of Proposition~\ref{prop:coupling_modf}}

In this section, we first construct a coupling that enables us to relate the modified system presented above to the original system operating under any load balancing policy.  To construct the original system under any admissible load balancing policy $\Pi$, we define
\begin{align*}
    A_{j,i}^{(n,\Pi)}(t) &:=\mbox{The number of arrivals assigned to queues with $i$ or more jobs in pool}\\
    &\qquad \mbox{$j\in [M]$ by time $t$}\\
    A_i^{(n,\Pi)}(t) &:=\mbox{The total number of arrivals assigned to queues with $i$ or more jobs by time $t$} \\
    D_{j,i}^{(n,\Pi)}(t) &:=\mbox{The number of departures from queues with $i$ or more jobs in pool $j\in [M]$}\\
      &\qquad \mbox{by time $t$}\\
    D_{i}^{(n,\Pi)}(t) &:=\mbox{The total number of departures from queues with $i$ or more jobs by time $t$}\\
    Q_1^{(n,\Pi)}(t)&:=\sum_{j=1}^M Q_{j,1}^{(n,\Pi)}(t)\\
    Q_{+2}^{(n,\Pi)}(t)&:= \sum_{j=1}^M \sum_{i=2}^\infty  Q_{j,i}^{(n,\Pi)}(t).
\end{align*}
We note that the processes $A_{j,i}^{(n,\Pi)}(t)$ is determined by the load balancing policy $\Pi$ and is not used explicitly in the remainder of this section. We only need that for all $t\geq 0$ the following equations hold
\[
A_i^{(n,\Pi)}(t) = \sum_{j=1}^M A_{j,i}^{(n,\Pi)}(t)\mbox{ and } D_i^{(n,\Pi)}(t)= \sum_{j=1}^M D_{j,i}^{(n,\Pi)}(t).
\]
Furthermore, we have the following balance equations
\begin{align}
Q_1^{(n,\Pi)}(t) &= Q_1^{(n,\Pi)}(0) + A_0^{(n,\Pi)}(t) - A_1^{(n,\Pi)}(t) - D_1^{(n,\Pi)}(t) + D_2^{(n,\Pi)}(t),\label{eq:balance_eq_o1}\\
Q_{+2}^{(n,\Pi)}(t) &=Q_{+2}^{(n,\Pi)}(0) + A_1^{(n,\Pi)}(t) -D_2^{(n,\Pi)}(t). \label{eq:balance_eq_o2}
\end{align}

Now, to couple the original system and the modified system, we take $A(t)$ and $D(t)$ to be two unit rate Poisson processes and assume that for both original and modified systems, we have
\begin{equation}
A_0^{(n,\Pi)}(t) = \tilde{A}^{(n)}(t) =A(n\lambda^{(n)}t),
\label{eq:arrival_coupling}
\end{equation}
i.e., both original and the modified systems have the same arrival streams. For both systems, we define the ``potential'' departure process as 
\[
D_P^{(n)}(t):= D\left(\sum_{j=1}^M \mu_j N_j^{(n)}\right), \mbox{ for all $t\geq 0$}
\]
and we let $\theta_l^{(n)}$ to be the $l$th event epoch for $D_P^{(n)}(t)$ and take $\{U_l^{(n)}\}_{l\in \ZZ_+}$ to be a sequence of independent uniform(0,1) random variables. 

To convert the potential departure process to the actual departure process in the original system, we first order the pairs of positive integers $(j,i)$ such that $(j,i)\geq (j',i')$ if $i>i'$ or $i=i'$ and $j<j'$. Then, the potential departure at epoch $\theta_l^{(n)}$ is assumed to contribute an actual departure for $D_{j,i}^{(n,\Pi)}(t)$, i.e., $D_{j, i}^{(n,\Pi)}(\theta_l^{(n)}) - D_{j, i}^{(n,\Pi)}(\theta_l^{(n)}-)=1$, if and only if $(j,i)$ is the supremum pair for which the following holds
\begin{equation}
U_l^{(n)}\leq \frac{\sum_{k=1}^j\mu_kQ_{k,i}^{(n,\Pi)}(\theta_l^{(n)}-)}{\sum_{j=1}^M\mu_jN_j^{(n)}},
\label{eq:service_coupling1}
\end{equation}
where the supremum should be understood in terms of the ordering above. If no pair $(j,i)$ satisfying the above inequality exists, then the potential departure assumed to be lost, i.e., it does not correspond to an actual departure.

In the modified system, the potential departure occuring at time $\theta_l^{(n)}$ is an actual departure from a queue with two customers if
\begin{equation}
U_l^{(n)}\leq \frac{\tilde{\mu}^{(n)}(\tilde{Q}_2^{(n)}(t-))}{\sum_{j=1}^M\mu_jN_j^{(n)}}, 
\label{eq:service_coupling2}
\end{equation}
or an actual departure from a queue with one customer if
\begin{equation}
\frac{\tilde{\mu}^{(n)}(\tilde{Q}_2^{(n)}(t-))}{\sum_{j=1}^M\mu_jN_j^{(n)}}< U_l^{(n)}\leq \frac{\tilde{\mu}^{(n)}(\tilde{Q}_1^{(n)}(t-))}{\sum_{j=1}^M\mu_jN_j^{(n)}}.
\label{eq:service_coupling3}
\end{equation}
The potential departure is assumed to be lost (i.e., does not correspond to an actual departure) if none of the above two conditions is satisfied.


{\em Proof of Proposition~\ref{prop:coupling_modf}.}
We prove the result with a contradiction. First, define $\tau_1^{(n)}:=\inf\{t\geq 0: \tilde{Q}_2^{(n)}(t)> Q_{+2}^{(n,\Pi)}(t)\}$. If $\tau_1^{(n)}<\infty$, we need to have $\tilde{Q}_2^{(n)}(\tau_1^{(n)}-)= Q_{+2}^{(n,\Pi)}(\tau_1^{(n)}-)$, which together with the assumption on initial conditions, \eqref{eq:balance_eq_m2} and \eqref{eq:balance_eq_o2}, yields
\begin{align}
    A_1^{(n,\Pi)}(\tau_1^{(n)}-) -D_2^{(n,\Pi)}(\tau_1^{(n)}-) &=   \tilde{A}_1^{(n)}(\tau_1^{(n)}-) -\tilde{D}_2^{(n)}(\tau_1^{(n)}-), \mbox{ and } \label{eq:q1_tau1}\\
    A_1^{(n,\Pi)}(t) -D_2^{(n,\Pi)}(t) &\geq   \tilde{A}_1^{(n)}(t) -\tilde{D}_2^{(n)}(t)\mbox{ for all }t<\tau_1^{(n)}.
 \label{eq:q1_tau2}
\end{align}
Equations \eqref{eq:rate_bound},\eqref{eq:service_coupling1} and \eqref{eq:service_coupling2} imply that the event occuring at time $\tau_1^{(n)}$ cannot be a departure. Hence, we need to have an arrival at $\tau_1^{(n)}$ and this arrival should be routed to a queue with one job in the modified system and it should be routed to an idle queue in the original system, which can only happen if
\[
Q_1^{(n,\Pi)}(\tau_1^{(n)}-)<\tilde{Q}_1^{(n)}(\tau_1^{(n)}-)=\sum_{j=1}^M N_j^{(n)}.
\]
Combining \eqref{eq:arrival_coupling}, \eqref{eq:q1_tau1},\eqref{eq:balance_eq_m1} and \eqref{eq:balance_eq_o1}, we need to have
\[
    \tilde{A}_R^{(n)}(\tau_1^{(n)}-) + \tilde{D}_1^{(n)}(\tau_1^{(n)}-) <D_1^{(n,\Pi)}(\tau_1^{(n)}-).
\]
Hence, there should exist a minimal $\tau_2^{(n)}<\tau_1^{(n)}$ where the equality $\tilde{A}_R^{(n)}(s) + \tilde{D}_1^{(n)}(s)  = D_1^{(n,\Pi)}(s)$ is violated (note that this equality trivially holds at $s=0$). We now show  that this is not possible even when $\tau_1^{(n)}=\infty$. Using a similar logic as above, if such a $\tau_2^{(n)}$ exists, we need to have 
\begin{equation*}
    \tilde{A}_R^{(n)}(\tau_2^{(n)}-) + \tilde{D}_1^{(n)}(\tau_2^{(n)}-)  = D_1^{(n,\Pi)}(\tau_2^{(n)}-).
\end{equation*}
This can only happen if there is a departure only in the original system at time $\tau_2^{(n)}$, i.e.,
\[
\tilde{\mu}^{(n)}\left(\tilde{Q}_1^{(n)}(\tau_2^{(n)}-)\right)< \sum_{j=1}^M \mu_j Q_{1,j}^{(n,\Pi)}(\tau_2^{(n)}-).
\]
However, using the inequality in \eqref{eq:q1_tau1} and \eqref{eq:q1_tau2} in combination with \eqref{eq:balance_eq_m1} and \eqref{eq:balance_eq_o1}, we see that $\tilde{Q}_1^{(n)}(\tau_2^{(n)}-)\geq Q_1^{(n)}(\tau_2^{(n)}-)$. The monotonicity of $\tilde{\mu}(\cdot)$ and \eqref{eq:rate_bound} implies that a departure only in the original system at $\tau_2^{(n)}$ is not possible. Hence, we can conclude first that $\tau_1^{(n)} =\infty$ and $\tilde{Q}_2^{(n)}(t)\leq Q_{+2}^{(n,\Pi)}(t)$. The non-existence of $\tau_2^{(n)}<\tau_1^{(n)}=\infty$ implies that $\tilde{A}_R^{(n)}(t) + \tilde{D}_1^{(n)}(t) \geq D_1^{(n,\Pi)}(t)$ for all $t\geq 0$ and hence summing \eqref{eq:balance_eq_m1}, \eqref{eq:balance_eq_m2}, \eqref{eq:balance_eq_o1} and \eqref{eq:balance_eq_o2}, we see that $\tilde{Q}_1^{(n)}(t) + \tilde{Q}_2^{(n)}(t)\leq Q_1^{(n,\Pi)}(t) + Q_{+2}^{(n,\Pi)}(t)$.
\qed
\endproof

\subsection{Diffusion Limit of the Modified System: Proof of Proposition~\ref{prop:modified_diffusion_limit}}

To prove Proposition~\ref{prop:modified_diffusion_limit}, we follow similar steps as in the proof of Theorem~\ref{thm:SA-JSQ_diffusion_limit}. In particular, (i) we first define a stopping time $\tilde T_n$ to be the first time the system reaches a state with at least $N_M^{(n)}/2$ idle servers or  at least $N_1^{(n)}$ queues with two jobs, (ii) we then derive the diffusion limit of the system for all $t \in [0, \tilde{T}_n]$, (iv) and finally show that, for any finite $t$, the probability of $T_n\leq t$ decreases to $0$ as $n\to \infty$ to obtain process-level convergence. Finally, we present a Lyapunov-type argument to argue tightness of the stationary distributions to conclude interchange of limits.  

{\em Proof of Proposition~\ref{prop:modified_diffusion_limit}.} To derive the diffusion limit for the modified system, we first define the stopping time $\tilde{T}_n$ as follows
\begin{equation}
\tilde{T}_n = \inf\{t\geq 0: n- \tilde{Q}_1^{(n)}(t)\geq N_M^{(n)}/2\text{ or }\tilde{Q}_2^{n}(t)=N_1^{(n)}\}.
\end{equation}
This definition implies before time $\tilde{T}_n$ only the slowest servers (i.e., servers in pool $M$) can be idle and no job can be blocked. 
Hence, using the unit rate Poisson processes $\tilde{A}(t)$, $\tilde{D}_{u1}(t)$ and $\tilde{D}_{u2}(t)$ the evolution of the modified system up to time $\tilde{T}_n$ can be described by the following equations 
\begin{align}
\tilde{Q}_1^{(n)}(t) &= \tilde{Q}_1^{(n)}(0) + \tilde{A}(n\lambda^{(n)}t)\nonumber\\
&\qquad-  \tilde{D}_{u1}\brac{\sum_{j=1}^{M}N_{j}^{(n)}\mu_j t-\int_{0}^t\brac{\mu_M(n-\tilde{Q}_1^{n}(s))+\mu_1\tilde{Q}_2^{(n)}(s)}ds} - \tilde{U}_1^{(n)}(t),\label{eq:balance_eq_m1a}\\
\tilde{Q}_{2}^{(n)}(t) &=\tilde{Q}_{2}^{(n)}(0) -\tilde{D}_{u2}\brac{\int_{0}^t\mu_1\tilde{Q}_2^{(n)}(s)ds} + \tilde{U}_1^{(n)}(t) \label{eq:balance_eq_m2a},
\end{align}
where $\tilde{U}_1^{(n)}(t)=\tilde{A}_1^{(n)}(t)=\int_0^t \indic{Q_1^{(n)}(s)=n}d\tilde{A}(n\lambda^{(n)}s)$ for all $0\leq t\leq \tilde{T}_n$. Upon centring and scaling the processes as in Section~\ref{sec:diffusion_sajsq}, we get
\begin{align}
\nonumber\tilde{Y}_1^{(n)}(t) &= \tilde{Y}_1^{(n)}(0) + \tilde{M}_A^{(n)}(t) -  \tilde{M}_{D,1}^{(n)}(t) +\frac{n\lambda^{(n)}-\sum_{j=1}^{M}\mu_jN_{j}^{(n)}}{\sqrt{n}}t\\
&\qquad - \mu_M\int_{0}^t\tilde{Y}_1^{n}(s)ds+\mu_1\int_0^t\tilde{Y}_2^{(n)}(s)ds - \tilde{V}_1^{(n)}(t)\label{eq:balance_eq_m1b}\\
\tilde{Y}_{2}^{(n)}(t) &=\tilde{Y}_{2}^{(n)}(0) -\tilde{M}_{D,2}^{(n)}(t) - \int_{0}^t\mu_1\tilde{Y}_2^{(n)}(s)ds + \tilde{V}_1^{(n)}(t) \label{eq:balance_eq_m2b},
\end{align}
where
\begin{align*}
    \tilde{M}_A^{(n)}(t) &= \frac{1}{\sqrt{n}}\brac{\tilde A(n\lambda^{(n)}t) -n\lambda^{(n)}t},\\
    \tilde{M}_{D,1}^{(n)}(t) & = \frac{1}{\sqrt{n}}\tilde{D}_{u1}\brac{\sum_{j=1}^{M}N_{j}^{(n)}\mu_j t -\int_{0}^t\brac{\mu_M(n-\tilde{Q}_1^{n}(s))+\mu_1\tilde{Q}_2^{(n)}(s)}ds}\\&\qquad-\frac{1}{\sqrt{n}}\brac{\sum_{j=1}^{M}N_{j}^{(n)}\mu_j t -  \int_{0}^t\brac{\mu_M(n-\tilde{Q}_1^{n}(s))+\mu_1\tilde{Q}_2^{(n)}(s)}ds},\\
    \tilde{M}_{D,2}^{(n)}(t) & = \frac{1}{\sqrt{n}}\brac{\tilde{D}_{u2}\brac{\mu_1\int_{0}^t\tilde{Q}_2^{(n)}(s)ds}-\int_{0}^t\mu_1\tilde{Q}_2^{(n)}(s)ds},\\
    \tilde{V}_1^{(n)}(t) &= \frac{\tilde{U}_1^{(n)}(t)}{\sqrt{n}}.
\end{align*}
Assumptions \eqref{eqn:heavy_traffic_limit} and \eqref{eq:proportion_assumption} imply that the following convergence holds  uniformly on compact intervals
\[
\frac{n\lambda^{(n)}-\sum_{j=1}^{M}\mu_jN_{j}^{(n)}}{\sqrt{n}}t\to -\beta e(t),
\]
where $e(t)=t$ for all $t$. Hence, the continuous mapping theorem and Lemma~\ref{lem:two-dim-sys} in the appendix imply that the solution of \eqref{eq:balance_eq_m1b} and \eqref{eq:balance_eq_m2b} converges to the solution of \eqref{eq:tl_Y_1M}-\eqref{eqn:diff_modified_fluct_term}, once we prove that $\tilde{M}_A^{(n)}(t)$ and $ \tilde{M}_{D,1}^{(n)}(t)$ each converges weakly to standard Brownian motions and $\tilde{M}_{D,2}^{n}(t)$ converges weakly to the zero function.  As in the proof of Lemma~\ref{lem:mart_convergence} we have that $\tilde{M}_A^{(n)}\Rightarrow W_{\tilde A}$ since $\langle \tilde{M}_A^{(n)}\rangle (t) =\lambda^{(n)}t \to t$ where $W_{\tilde A}$ is a standard Borwnian motion. Furthermore, the predictable quadratic variations $\langle \tilde{M}_{D,1}^{(n)}\rangle$ and   $\langle \tilde{M}_{D,2}^{(n)}\rangle$ of the martingales $\tilde{M}_{D,1}^{(n)}$ and $\tilde{M}_{D,2}^{n}$ satisfy
\begin{align*}
   \langle \tilde{M}_{D,1}^{(n)} \rangle(t) &=\frac{\sum_{j=1}^{M}N_{j}^{(n)}\mu_j t}{n} -  \int_{0}^t\brac{\frac{\mu_M(n-\tilde{Q}_1^{n}(s))}{n}+\frac{\mu_1\tilde{Q}_2^{(n)}(s)}{n}ds},\\
   &\leq \frac{\sum_{j=1}^{M}N_{j}^{(n)}\mu_j t}{n} \\
   \langle \tilde{M}_{D,2}^{(n)}\rangle (t)& = \frac{\int_{0}^t\mu_1\tilde{Q}_2^{(n)}(s)ds}{n}\leq \frac{\mu_1N_1^{(n)}t}{n},
\end{align*}
which by~\eqref{eq:proportion_assumption} and Lemma~\ref{lem:stochastic_bounded}
imply that  $\{\tilde{M}_{D,1}^{(n)}\}_n$ and $\{\tilde{M}_{D,2}^{(n)}\}_n$ are stochastically bounded. Hence, Lemma~\ref{lem:stochastic_boundedness} implies that $\{\tilde{Y}_1^{(n)}\}_n$ and $\{\tilde{Y}_2^{(n)}\}_n$ are stochastically bounded. Therefore, by Lemma 5.9 of~\cite{Pang2007} we have $\tilde{Y}_1^{(n)}/\sqrt{n}=(\tilde{Q}_1^{(n)}-n)/n\Rightarrow 0$ and $\tilde{Y}_2^{(n)}/\sqrt{n}=\tilde{Q}_2^{(n)}/n\Rightarrow 0$. Therefore, by~\eqref{eq:proportion_assumption} and continuous mapping theorem we have $\langle \tilde{M}_{D,1}^{(n)}\rangle \Rightarrow e$ and $\langle \tilde{M}_{D,2}^{(n)}\rangle \Rightarrow 0$. This proves that the solution of \eqref{eq:balance_eq_m1b} and \eqref{eq:balance_eq_m2b} converges to the solution of \eqref{eq:tl_Y_1M} and \eqref{eq:tl_Y_12}. 

To complete the proof we need to show that $\PP(\tilde{T}_n\leq t)\to 0$ for any finite $t\geq 0$. We have
\begin{align*}
    \PP\brac{\tilde{T}_n\leq t}&\leq \PP\brac{\sup_{0\leq s\leq t}|\tilde{Y}_1^{(n)}(s)|\geq\frac{N_M^{(n)}}{2\sqrt{n}}} + \PP\brac{\sup_{0\leq s \leq t}\tilde{Y}_2^{(n)}(s)\geq \frac{N_1^{(n)}}{\sqrt{n}}}.
\end{align*}
Assumption~\eqref{eq:proportion_assumption}
and the stochastic boundedness of $\tilde{Y}^{(n)}_1$ and $\tilde{Y}^{(n)}_2$, imply that the RHS of the above inequality goes to zero as $n \to \infty$. 

To prove the convergence of stationary distributions, we need to prove that the sequences $\{\tilde{Y}_1^{(n)}(\infty)\}_n$ and $\{\tilde{Y}_2^{(n)}(\infty)\}_n$ are tight. From Proposition~\ref{prop:coupling_modf} we have  $0\leq \tilde{Y}_2^{(n)}(\infty)\leq_{st} Y_{+2}^{(n)}(\infty)$. Furthermore, Theorem~\ref{thm:stationarity_diffusion} implies that the sequence $\{Y_{+2}^{(n)}(\infty)\}_n$ weakly converges to $Y_{1,2}(\infty)$ and is therefore tight.
Hence, $\{\tilde{Y}_2^{(n)}(\infty)\}_n$ is tight. To prove the tightness of $\{\tilde{Y}_1^{(n)}(\infty)\}_n$, we consider the Lyapunov function $\tilde V$ defined for any state $\tilde q$ of the modified system as
\[\tilde V( \tilde{\mf q})= -\tilde{y}_1\geq 0.\]
By definition, we have 
\[\tilde\mu^{(n)}\brac{\tilde{q}_1}\leq \sum_{j=1}^M N_j^{(n)}\mu_j-\sqrt{n}\mu_M\tilde{V}(\tilde{\mf q}).\]
Hence, applying the generator $\tilde{\cG}$ of the process $\tilde{\mf Q}^{(n)}$ to $\tilde{V}$ we obtain
\begin{align*}
    \tilde{\mathcal{G}}\tilde{V}( \tilde{\mf q}) &= -\sqrt{n}\lambda^{(n)}\indic{\tilde{y}_1<0} + \frac{\tilde\mu^{(n)}\brac{\tilde{q}_1}}{\sqrt{n}}\\
    &\leq -\sqrt{n}\lambda^{(n)} + \frac{\sum_{j=1}^M N_j^{(n)}\mu_j}{\sqrt{n}} - \mu_M\tilde{V}(\tilde{\mf q})\leq \beta-\mu_M\tilde{V}(\tilde{\mf q})+o(1),
\end{align*}
where the last inequality follows from~\eqref{eqn:heavy_traffic_limit} and~\eqref{eq:proportion_assumption}.
Hence, for $\tilde{\mf q}$ satisfying $\tilde{V}(\tilde{\mf q})>2\beta/\mu_M$ we have $\tilde{\mathcal{G}}\tilde{V}( \tilde{\mf q})\leq -(\beta-o(1))$.
Then,  applying Lemma 10 of \cite{wang_etal2022}, we have
\begin{equation*}
\EE \sbrac{\tilde{V}(\tilde{\mf {Q}}^{(n)}(\infty))}= \EE \sbrac{-\tilde{Y}_1^{(n)}(\infty)}\leq \frac{4\beta}{\mu_M} + 4\frac{(\lambda^{(n)}+\sum_{j \in [M]}\mu_j N_j^{(n)}/{n}) + \beta-o(1))}{\beta-o(1)}=O(1),
\end{equation*}
which shows that $\{\tilde{Y}_1^{(n)}(\infty)\}_n$ is tight. This completes the proof of interchange of limits.
\qed
\endproof
}

\section{Conclusion}
\label{sec:conclusion}
In this paper, we established that the SA-JSQ policy is asymptotically optimal in the Halfin-Whitt traffic regime for a system with heterogeneous servers. In order to prove this result, we needed to do both transient and steady-state analysis of the SA-JSQ policy in addition to finding a suitable system which provides asymptotically tight performance lower bounds for any admissible load balancing policy. The key difficulty in our transient and steady-state analyses arises from the heterogeneity of the servers which results in a more complex state-space for the underlying Markov process. This makes generalising existing analyses hard as we need to establish additional state-space collapse results and construct Lyapunov functions suited to heterogeneous systems. To the best of our knowledge, our proof of asymptotic optimality is the first in the literature for a heterogeneous system where each server has a dedicate queue. It uses a system which differs fundamentally from those used in earlier works where the number of queues did not scale with the number of servers.

The paper opens up several interesting avenues for further exploration. 
We stated our optimality result in terms of minimising the steady-state distribution of the total diffusion-scaled number of jobs and the diffusion-scaled number of waiting jobs in the system. Another natural way to define optimality is through the diffusion-scaled mean response time of jobs and the diffusion-scaled mean waiting time of jobs in the steady state. Although one can relate the expected number of waiting jobs and expected total number of jobs in the system at steady-state to the mean waiting time and mean response times of jobs, respectively, one cannot conclude the optimality in terms of (diffusion-scaled) averages from our result as this requires showing uniform integrability of stationary queue lengths which we have not established. An interesting direction of future research would be to show asymptotic optimality in terms of averages by establishing uniform integrability results as in~\cite{braverman2023join}. 

In this paper, we focused on analysing the SA-JSQ policy under the condition $\lambda^{(n)}=1-\beta/n^{\alpha}$, with $\alpha=1/2$. In a related analysis in~\cite{bhambay2022asymptotic}, the optimality of the SA-JSQ was established for $\alpha=0$, corresponding to the fluid limit. However, asymptotically optimal policies remain unknown for other values of $\alpha$.  Analysing the SA-JSQ scheme for $\alpha=1$, corresponding to the NDS regime, is of particular interest as it could reveal interesting dynamics.

Finding optimal load balancing policies for non-exponential service time distributions remains another challenging open problem.


\bibliographystyle{unsrt}
\bibliography{sample}

\newpage

\appendix

\section{Proof of Proposition~\ref{thm:skorohod_mapping}}
\label{append:1}

To prove Proposition~\ref{thm:skorohod_mapping}, we first prove an important result for first two dimensions involving reflection terms.

 \begin{lemma}
 \label{lem:two-dim-sys}
For a given $B\in \bar{\RR}_+$, $\mathbf{b}=(b_1,b_2)\in \RR^2$, and $\mathbf{z}=(z_1,z_2)\in \DD_{\RR^2}[0,\infty)$, consider the following integral equations
\begin{align}
y_1(t) &=b_1 +z_1(t) - \int_0^t(\mu_M y_1(s) - \mu_1 y_2(s))ds -u_1(t), \label{eqn:ref_1}\\
y_{1,2}(t)=y_2(t) &=b_2 +z_2(t) - \int_0^t\mu_1 y_2(s)ds +u_1(t)-u_2(t),\\
y_1(t)&\leq0, \ 0\leq y_2(t) \leq B, \ t\geq0, \label{eqn:ref_3}
\end{align}
where $u_1$ and $u_2$ are non-decreasing and non-negative functions satisfying
\begin{align}
\int_0^t\indic{y_1(s)<0}du_1(s)&=0,\label{eqn:ref4}\\
\int_0^t\indic{y_2(s)<B}du_2(s)&=0.\label{eqn:ref5}
\end{align}
Then~\eqref{eqn:ref_1}-\eqref{eqn:ref_3} has a unique solution in $(\mf y,\mf u) \in \DD_{\RR^2}[0,\infty) \times \DD_{\RR^2}[0,\infty)$. Moreover, there exists a well defined  function $(f,  g): \bar{\RR}_+ \times \RR^2 \times \DD_{\RR^2}[0,\infty) \to \DD_{\RR^2}[0,\infty) \times \DD_{\RR^2}[0,\infty)$ which maps $(B,\mf b,\mf z)$ in to $\mf y =f(B,\mf b,\mf z)$ and $\mf u=g(B,\mf b,\mf z)$. Furthermore, the function $(f,g)$ is continuous on $\bar{\RR}_+ \times \RR^2 \times \DD_{\RR^2}[0,\infty)$. Finally, $\mf z$ continuous implies that $\mf y$ and $\mf u$ are also continuous.
 \end{lemma}

To prove Lemma~\ref{lem:two-dim-sys}, we first define a one sided reflection map function with upper barrier $\kappa$ as $(\phi_\kappa,\psi_\kappa):\DD_{\RR}[0,\infty) \to \DD_{\RR^2}[0,\infty)$ where for any $x\in \DD_{\RR}[0,\infty)$ and for all $t\geq 0$ 
\begin{align*}
&\phi_\kappa(x)=x-\psi_\kappa(x)\leq \kappa,\\
&\int_0^\infty \indic{\phi_\kappa(x)<\kappa}d\psi_k(x)=0 \mbox{ and }\\
&\psi_\kappa(x)\mbox{ is nondecreasing and }\psi_\kappa(x(0))=0.
\end{align*}
It is well-known in the literature that the reflection map $(\phi_\kappa,\psi_\kappa)$ is well-defined (see, e.g. \cite{whitt2002stochastic} Sections 5.2 and 13.5) and can be expressed as 
\begin{align}
\psi_\kappa(x)(t)=\sup_{0\leq s\leq t} (x(s)-\kappa)^+, \label{eqn:phi}\\
\phi_\kappa(x)(t)=x(t)-\psi_\kappa(x)(t)\label{eqn:psi}.
\end{align}
The reflection mapping can be trivially extended to $\kappa=\infty$ as $(\phi_\infty(x),\psi_\infty(x))=(x,0)$. The Lipschitz continuity of the reflection mapping, expressed as
\begin{align}\label{eq:reflection_lipschitz1}
    ||\psi_\kappa(x)-\psi_{\kappa'}(x')||_t&\leq ||x-x'||_t +|\kappa-\kappa'|,\\
    \label{eq:reflection_lipschitz2}||\phi_\kappa(x)-\phi_{\kappa'}(x')||_t&\leq 2||x-x'||_t +|\kappa-\kappa'|.
\end{align}
is used in our proofs. To be able to use the continuous mapping theorem, we need the continuity to be extended to cover the order topology for extended reals. 
\begin{lemma}[\cite{Eschenfeldt2018}, Lemma 2]\label{lem:EG_continuity} The reflection mapping $(\phi,\psi):\bar{\RR}_+\times \DD_{\RR}[0,\infty)\to \DD_{\RR^2}[0,\infty)$ is continuous with respect to the product topology where $\bar{\RR}_+$  and $\DD_{\RR}[0,\infty)$ are equipped with the order topology and topology of uniform convergence on compact sets, respectively. 
\end{lemma}
We also need the following Gronwall-type inequality:
\begin{lemma}[Gronwall Inequality, \cite{Greene1977}, ~\cite{Das1979}]\label{lem:gronwall}
    Let $K_1, K_2\in \RR_+$ and $h_i\in \RR$ for $i=1,..,4$ and let $f$ and $g$ to be continuous nonnegative functions such that 
    \begin{align*}
        f(t) &\leq K_1 + h_1\int_0^tf(s)ds + h_2 \int_0^tg(s)ds,\\
        g(t) &\leq K_2 + h_3\int_0^tf(s)ds + h_4 \int_0^tg(s)ds
    \end{align*}
    for all $t\geq 0$. Then, 
    \begin{align*}
        f(t)\leq (K_1+K_2)e^{ht}\mbox{ and }g(t)\leq (K_1+K_2)e^{ht},
    \end{align*}
    for all $t\geq 0$, where $h=\max\{h_1+h_3,h_2+h_4\}$.
\end{lemma}
The following lemma is a version of Lemma 3 which addresses \eqref{eqn:ref4} and \eqref{eqn:ref5} internally using the reflection mapping. 
\begin{lemma}
\label{lem:two-dim-sys_2}
Consider the following equations
\begin{align}
w_1(t) &=b_1 +z_1(t) - \int_0^t(\mu_M \phi_0(w_1)(s) - \mu_1 \phi_B(w_2)(s))ds, \label{eqn:ref_1_1}\\
w_2(t) &=b_2 +z_2(t)+ \psi_0(w_1)(t) - \int_0^t\mu_1 \phi_B(w_2)(s)ds,\label{eqn:ref_3_1}
\end{align}
For given $B\in \bar{\RR}_+, \mathbf{b}\in \RR^2$ and $\mathbf{z}\in \DD_{\RR^2}[0,\infty)$, \eqref{eqn:ref_1_1} and \eqref{eqn:ref_3_1} has a unique solution $w$ and defines a function $\bs \xi=(\xi_1,\xi_2):\bar{\RR}_+\times \DD_{\RR}[0,\infty)\to \DD_{\RR^2}[0,\infty)$ such that $\mathbf{w}=\bs \xi(B,\mathbf{b},\mathbf{z})$. The function $\bs \xi$ is continuous and the continuity of $\mathbf{z}$ implies the continuity of $\mathbf{w}$. 
\end{lemma}
{\em Proof.}
In the first part of the proof, we show the existence of a solution using a fixed point argument. We consider the version of  \eqref{eqn:ref_1_1} and \eqref{eqn:ref_3_1} by replacing $(\tilde{w}_1, \tilde{w}_2) = (w_1, w_2-\psi_0(w_1))$ as   
\begin{align}
\tilde{w}_1(t) &=b_1 +z_1(t) - \int_0^t(\mu_M \phi_0(\tilde{w}_1)(s) - \mu_1 \phi_B(\tilde{w}_2+\psi_0(\tilde{w}_1))(s))ds, \label{eqn:ref_1_2}\\
\tilde{w}_2(t) &=b_2 +z_2(t)- \int_0^t\mu_1 \phi_B(\tilde{w}_2+\psi_0(\tilde{w}_1))ds,\label{eqn:ref_3_2}
\end{align}
For fixed $(B, \mathbf{b}, \mathbf{y})$,a solution to \eqref{eqn:ref_1_2} and \eqref{eqn:ref_3_2} is the fixed point of the operator $\cT:\DD_{\RR^2}[0,\infty) \to \DD_{\RR^2}[0,\infty)$ defined by 
\begin{align*}
\cT(\tilde{\mf w})_1(t) &=b_1 +z_1(t) - \int_0^t(\mu_M \phi_0(\tilde{w}_1)(s) - \mu_1 \phi_B(\tilde{w}_2+\psi_0(\tilde{w}_1))(s))ds,\\
\cT(\tilde{\mf w})_2(t) &=b_2 +z_2(t)- \int_0^t\mu_1 \phi_B(\tilde{w}_2+\psi_0(\tilde{w}_1))ds.
\end{align*}
For any $\mathbf{\tilde{w},\tilde{v}}\in \DD_{\RR^2}[0,\infty)$,
\begin{align*}
    ||\cT(\mathbf{\tilde{w}})_1-\cT(\mathbf{\tilde{v}})_1||_t &\leq \mu_Mt\left(||\phi_0(\tilde{w}_1)-\phi_0(\tilde{v}_1)||_t +||\phi_B(\tilde{w}_2+\psi_0(\tilde{w}_1))-\phi_B(\tilde{v}_2+\psi_0(\tilde{v}_1))||_t\right)\\
    &\leq \mu_1t\left(2||\tilde{w}_1-\tilde{v}_1||_t +2||\tilde{w}_2+\psi_0(\tilde{w}_1)-\tilde{v}_2-\psi_0(\tilde{v}_1)||_t\right)\\
    &\leq \mu_1t\left(4||\tilde{w}_1-\tilde{v}_1||_t +2||\tilde{w}_2-\tilde{v}_2||_t\right)\leq 6\mu_1t ||\mft{w}-\mft{v}||_t,
\end{align*}
and 

\begin{align*}
    ||\cT(\mathbf{\tilde{w}})_2-\cT(\mathbf{\tilde{v}})_2||_t &\leq \mu_Mt\left( ||\phi_B(\tilde{w}_2+\psi_0(\tilde{w}_1))-\phi_B(\tilde{v}_2+\psi_0(\tilde{v}_1))||_t\right)\\
    &\leq \mu_1t\left(2||\tilde{w}_2+\psi_0(\tilde{w}_1)-\tilde{v}_2-\psi_0(\tilde{v}_1)||_t\right)\\
    &\leq \mu_1t\left(2||\tilde{w}_1-\tilde{v}_1||_t +2||\tilde{w}_2-\tilde{v}_2||_t\right)\leq 4\mu_1t ||\mft{w}-\mft{v}||_t.
\end{align*}
Therefore, we have 
\begin{equation*}
     ||\cT(\mathbf{\tilde{w}})-\cT(\mathbf{\tilde{v}})||_t \leq  6\mu_1t ||\mft{w}-\mft{v}||_t.
\end{equation*}
This implies that the mapping $\cT$ is a contraction on $\DD_{\RR^2}[0,t_0]$ for any $t_0\leq (6\mu_1)^{-1}$. Now, repeating the same argument iteratively for intervals $[t_0,2t_0],[2t_0, 3t_0], \ldots$, and  using the contraction mapping theorem (c.f. \cite{Rudin87}, Theorem 9.23), we conclude that a fixed point to the operator $\cT$. Hence, a solution to \eqref{eqn:ref_1_2} and \eqref{eqn:ref_3_2} exists and $(w_1, w_2)=(\tilde{w}_1, \tilde{w}_2+\psi_0(\tilde{w}_1))$ solves \eqref{eqn:ref_1_1} and \eqref{eqn:ref_3_1}.
Now, we prove uniqueness. If $\mft{w}$ and $\mft{w}'$ both solve \eqref{eqn:ref_1_1} and \eqref{eqn:ref_3_1}, then
\begin{align*}
    ||w_1-w_1'||_t&\leq \mu_M  \int_0^t||\phi_0(w_1)-\phi_0(w_1')||_sds+\mu_1\int_0^t||\phi_B(w_2)-\phi_B(w_2')||_sds\\
    &\leq 2\mu_M  \int_0^t||w_1-w_1'||_sds+2\mu_1\int_0^t||w_2-w_2'||_sds\\
    &\leq 4\mu_1 \int_0^t||w_1-w_1'||_sds+2\mu_1\int_0^t(||w_2-w_2'||_s-||w_1-w_1'||_s)^+ds\\
\end{align*}
and 
\begin{align*}
    ||w_2-w_2'||_t&\leq ||\psi_0(w_1)-\psi_0(w_1')||_t +\mu_1\int_0^t||\phi_B(w_2)-\phi_B(w_2')||_sds\\
    ||w_2-w_2'||_t-||w_1-w_1'||_t &\leq  2\mu_1\int_0^t||w_1-w_1'||_sds+2\mu_1\int_0^t(||w_2-w_2'||_s-||w_1-w_1'||_s)^+ds\\
    (||w_2-w_2'||_t-||w_1-w_1'||_t)^+ &\leq  2\mu_1\int_0^t||w_1-w_1'||_sds+2\mu_1\int_0^t(||w_2-w_2'||_s-||w_1-w_1'||_s)^+ds,
\end{align*}
where the last inequality follows due to the nonnegativity of the right-hand side. Now, applying Lemma~\ref{lem:gronwall}, we have 
\[
||w_1-w_1'||_t=||w_2-w_2'||_t = 0,
\]
which implies uniqueness of the solution. 

To prove continuity, consider a sequence $(B^n,\mf{b}^n, \mf{z}^n)\to (B,\mf{b}, \mf{z})$. There exists solutions $\mf{w}^n$ and $\mf{w}$ corresponding to the solution of \eqref{eqn:ref_1_1} and \eqref{eqn:ref_3_1} for $(B^n,\mf{b}^n, \mf{z}^n)$ and $(B,\mf{b}, \mf{z})$, respectively. Lemma~\ref{lem:EG_continuity} implies that for any $\delta>0$, there exists an $N_\delta$ such that for $n>N_\delta$ implies
\[
|\mf{b}^n-\mf b|+||\mf{z}^n-\mf{z}||_t+||\phi_{B^n}(w_2)-\phi_{B}(w_2)||<\delta.
\]
For $n>N_\delta$,
\begin{align*}
    ||w_1^n-w_1||_t&\leq \delta+\mu_M  \int_0^t||\phi_0(w_1^n)-\phi_0(w_1)||_sds+\mu_1\int_0^t||\phi_{B^n}(w_2^n)-\phi_B(w_2)||_sds\\
    &\leq\delta+\mu_M  \int_0^t||\phi_0(w_1^n)-\phi_0(w_1)||_sds+\mu_1\int_0^t||\phi_{B^n}(w_2^n)-\phi_{B^n}(w_2)||_sds\\
    &\qquad + \mu_1\int_0^t||\phi_{B^n}(w_2)-\phi_{B}(w_2)||_sds\\
    &\leq\delta(1+\mu_1t)+2\mu_M  \int_0^t||w_1^n-w_1||_sds+2\mu_1\int_0^t||w_2^n-w_2||_sds\\
    &\leq\delta(1+\mu_1t)+4\mu_M  \int_0^t||w_1^n-w_1||_sds+2\mu_1\int_0^t(||w_2^n-w_2||_s-||w_1^n-w_1||_s)^+ds\\
\end{align*}
 and similarly, 
 \begin{align*}
     (||w_2^n-w_2||_t - ||w_1^n-w_1||_t)^+ &\leq \delta(1+\mu_1t) + 2\mu_1\int_0^t||w_1^n-w_1||_sds \\&\qquad +2\mu_1\int_0^t(||w_2^n-w_2||_s-||w_1^n-w_1||_s)^+ds.
 \end{align*}
 Choosing $\delta=(4+4\mu_1t)^{-1} e^{-6\mu_1t}\epsilon$ and applying Lemma~\ref{lem:gronwall}, 
 \[
 ||w_1^n-w_1||_t<\epsilon \mbox{ and }||w_2^n-w_2||_t <2\epsilon, 
 \]
 which proves the continuity of the mapping $\xi$.  Finally, to prove that continuity of $\mf z$ implies the continuity of $\mf w$, we first show that the function $w$ is bounded using a similar approach as above. Choosing $x=w$ and $x'$ as the zero function, for any $\kappa\geq 0$, \eqref{eq:reflection_lipschitz1} and \eqref{eq:reflection_lipschitz2} implies
 \begin{align*}
 ||\psi_\kappa(w)||_t\leq ||w||_t, ||\phi_\kappa(w)||_t\leq 2||w||_t.
 \end{align*}
 Now, re-organising \eqref{eqn:ref_1_1} and \eqref{eqn:ref_3_1}, we have
\begin{align*}
||w_1||_t &\leq |b_1| +||z_1||_t + (\mu_1+\mu_M)\int_0^t ||w_1||_sds + \mu_1\int_0^t (||w_2||_s-||w_1||_s)ds,\\
||w_2||_t -||w_1||_s &\leq |b_2| +||z_2||_t+ \mu_1\int_0^t||w_1||_sds +\mu_1\int_0^t (||w_2||_s-||w_1||_s)ds.
\end{align*}
Now, applying Lemma~\ref{lem:gronwall} again, we get
\begin{align*}
&||w_1||_t\leq (|b_1| +||z_1||_t +|b_2| +||z_2||_t)e^{(2\mu_1+\mu_M)t}\mbox{ and }\\
&||w_2||_t-||w_1||_t\leq (|b_1| +||z_1||_t +|b_2| +||z_2||_t)e^{(2\mu_1+\mu_M)t},
\end{align*}
which implies that $\mf{w}$ is bounded. Then, again using \eqref{eqn:ref_1_1} and \eqref{eqn:ref_3_1}, we get 
 \begin{align*}
 |w_1(t+s)-w_1(t)|&\leq |z_1(t+s)-z_1(t)|+\mu_1\int_t^{t+s}|\phi_0(w_1)(z)|+|\phi_B(w_2)(z)|dz,\\
 |w_2(t+s)-w_2(t)|&\leq |z_1(t+s)-z_1(t)|+|w_1(t+s)-w_1(t)|+ \mu_1\int_t^{t+s}|\phi_B(w_2)(z)|dz,
 \end{align*}
 and along with the boundedness of $\mf{w}$ this implies the desired continuity. 
\endproof

{\em Proof of Lemma~\ref{lem:two-dim-sys}.} The lemma follows using the same arguments in the proof of Lemma 1 in \cite{Eschenfeldt2018} by taking 
\begin{align*}
    &y_1=(\phi_0\circ \xi_1)(\mf{b},\mf{z}), \quad u_1=(\psi_0\circ \xi_1)(\mf{b},\mf{z})\\
    &y_2=(\phi_B\circ \xi_2)(\mf{b},\mf{z}), \quad u_2=(\psi_B\circ \xi_2)(\mf{b},\mf{z}).
\end{align*}
\endproof

{\em Proof of Proposition~\ref{thm:skorohod_mapping}.}
The existence, uniqueness and continuity of $y_{j,i}$ solving \eqref{eqn:map_3}-\eqref{eqn:map_5} can be shown by fixing $j$ and applying Lemma 5 in \cite{Eschenfeldt2018}. Defining 
\begin{align*}
    \hat{z}_{1,2}(t)=z_{1,2}(t) + \mu_1\int_{0}^ty_{1,3}(s)ds,
\end{align*}
and applying Lemma~\ref{lem:two-dim-sys} with $\hat{z}_{1,2}$, we can prove existence, uniqueness and continuity of the solution. To show the continuity of the mapping $(f,g)$, suppose that $(B^n, \mf{b}^n, \mf{z}^n)\to (B, \mf{b}, \mf{z})$ as $n\to \infty$ and $\mf{y}^n$ and $\mf{y}$ solve the respective \eqref{eqn:map_1}-\eqref{eqn:map_3}. Lemma~5 in \cite{Eschenfeldt2018} implies that $y_{j,i}^n\to y_{j,i}$ for all $i\geq 3$ and $(j,2)$ for $j\geq 2$, which implies
\begin{align*}
    z_{1,2}^n(t) + \mu_1\int_{0,t}y_{1,3}^n(s)ds\to z_{1,2}(t) + \mu_1\int_{0,t}y_{1,3}(s)ds.
\end{align*}
Combining this with the continuity in Lemma~\ref{lem:two-dim-sys}, we obtain the desired continuity. 
\qed
\endproof

\section{Proof of Lemma~\ref{lem:stochastic_boundedness}}

\label{append:2}
 Let $\hat{\mf Y}^{(n)}=(\hat{Y}_{1,2}^{(n)},\hat{Y}_{[1,M],1}^{(n)})$,  $\hat{\mf Z}^{(n)}= (\hat{Z}_1^{(n)},\hat{Z}_{1,2}^{(n)})$ and fix $t\geq0$. To prove the result of proposition we establish the following bound
\begin{equation}
\label{eqn:bound}
\norm{\hat{\mf Y}^{(n)}}_t \leq 8e^{6 \mu_1t}\brac{|\hat{\mf Y}^{(n)}(0)| + \norm{\hat{\mf Z}^{(n)}}_t }.
\end{equation}
To prove the above bound, we show that the similar bound holds for the unreflected process $\mf w^{(n)}$ defined in Lemma~\ref{lem:two-dim-sys_2}. Then~\eqref{eqn:bound} follows from the Lipschitz continuity of the reflected maps $\phi_0$ and $\phi_{B_n}$. We write $\hat{Y}_{[1,M],1}^{(n)}(t)=\phi_0(w_1^{(n)}(t))$ and $\hat{Y}_{1,2}^{(n)}(t)=\phi_{B_n}(w_2^{(n)}(t))$ where $w_1^{(n)}(t)$ and $w_2^{(n)}(t)$ satisfies
\begin{align*}
w_1^{(n)}(t) &=\hat{Y}_{[1,M],1}^{(n)}(0) +\hat{Z}_1^{(n)}(t) - \int_0^t(\mu_M \phi_0(w_1^{(n)})(s) - \mu_1 \phi_B(w_2^{(n)})(s))ds, \\
w_2^{(n)}(t) &=\hat{Y}_{1,2}^{(n)}(0) +\hat{Z}_{1,2}^{(n)}(t)+ \psi_0(w_1^{(n)})(t) - \int_0^t\mu_1 \phi_B(w_2^{(n)})(s)ds.
\end{align*}

Now using Grownwall's inequality from Lemma~\ref{lem:gronwall} and using Lipschitz property of $\phi_0$, $\phi_{B_n}$, $\psi_0$ we have for any $t\geq0$
\begin{align*}
\norm{w_1^{(n)}}_t &\leq |\hat{Y}_{[1,M],1}^{(n)}(0)| +\norm{\hat{Z}_1^{(n)}}_t +2 \mu_1 \int_0^t\brac{\norm{w_2^{(n)}}_s + \norm{w_1^{(n)}}_s}ds, \\
\norm{w_2^{(n)}}_t &\leq|\hat{Y}_{1,2}^{(n)}(0)| +\norm{\hat{Z}_1^{(n)}}_t + \norm{\psi_0(w_1^{(n)})}_t +\int_0^t\mu_1 \norm{w_2^{(n)}}_sds.
\end{align*}
Note that we have $||\psi_0(w_1^{(n)})||_t\leq ||w_1^{(n)}||_t$. Moreover, we define $r_1(t)=\norm{w_1^{(n)}}_t$ and $r_2(t)=\brac{\norm{w_2^{(n)}}_t-\norm{w_1^{(n)}}_t}^+$ and note that $\norm{w_2^{(n)}}_t \leq r_1(t)+r_2(t)$. Therefore, we can write above inequalities as
\begin{align*}
r_1(t) &\leq |\hat{Y}_{[1,M],1}^{(n)}(0)| +\norm{\hat{Z}_1^{(n)}}_t + 4\mu_1 \int_0^t r_1(s)ds + 2\mu_1 \int_0^t r_2(s)ds, \\
r_2(t) &\leq|\hat{Y}_{1,2}^{(n)}(0)| +\norm{\hat{Z}_1^{(n)}}_t + \int_0^t \mu_1 r_1(s)ds +\int_0^t\mu_1 r_2(s)ds.
\end{align*}
Let $|\hat{\mf Y}^{(n)}(0)| + \norm{\hat{\mf Z}^{(n)}}_t =m$. Now from Lemma~\ref{lem:gronwall}, we have $r_1(t)\leq 2m e^{6\mu_1 t} $ and $r_2(t)\leq 2m e^{6\mu_1 t}$. Finally, using the definition of $r_1$, $r_2$ and from the fact that $\phi_0$, $\phi_{B_n}$ are Lipschitz with constant $2$ we get
\begin{equation*}
 \norm{\hat{Y}_{[1,M],1}^{(n)}}_t \leq 4 m e^{6\mu_1 t}, \   \norm{\hat{Y}_{1,2}^{(n)}}_t \leq 8 m e^{6\mu_1 t},  
\end{equation*}
which proves~\eqref{eqn:bound}. Hence, the proof is complete.
\qed
\endproof

\section{Proof of Lemma~\ref{lem:simple_stationary_bounds}}
\label{app6}
To prove part 1 of the lemma, we define a Lyapunov function $g(\mf q)=\min\{T, \sum_{k=1}^M\sum_{i=1}^\infty q_{k,i}\}$. Then, 
\begin{align*}
    \cG g(\mf q) &=n\lambda^{(n)}1\left(\sum_{k=1}^M\sum_{i=1}^\infty q_{k,i}<T\right) - \sum_{k=1}^M\mu_k q_{k,1}1\left(\sum_{k=1}^M\sum_{i=1}^\infty q_{k,i}\leq T\right).
\end{align*}
Now taking the expectation of both sides with respect to the appropriate stationary measure, applying~\eqref{eq:rate_cons} and the monotone convergence theorem, we have $\EE[\sum_{k=1}^M \mu_k Q_{k,1}^{(n)}(\infty)]=n\lambda^{(n)}$. Hence,
\begin{align*}
0\leq \mu_M \EE\left[n^{-1/2}\left(\sum_{k=1}^M(N_k^n-Q_{k,1}^{(n)}(\infty))\right)\right]&\leq \EE\left[n^{-1/2}\left(\sum_{k=1}^M\mu_k(N_k^n-Q_{k,1}^{(n)}(\infty))\right)\right]\\
&=n^{1/2}\left((1-\lambda^{(n)})+\sum_{k=1}^M\mu_k\left(\frac{N_k^n}{n}-\gamma_k\right)\right),
\end{align*}
which implies tightness using \eqref{eqn:heavy_traffic_limit} and \eqref{eq:proportion_assumption}. Part 2 of the lemma can be proved in a similar fashion by using the  Lyapunov function 
\begin{equation*}
g(\mf q)=\min\left\{T, \sum_{k=2}^M\sum_{i=2}^\infty q_{k,i} + \sum_{i=3}^\infty q_{1,i}\right\}.
\end{equation*}
\endproof

\section{Proof of Lemma~\ref{lem:tau_solution}}
\label{proof:tau_solution}


For fixed $x_1\leq 0$, plugging in $x_2e^{-\mu_1\tau}={\kappa}/\sqrt{n}$, \eqref{eq:tau_equation} becomes
\begin{equation}\label{eq:solution_uniqueness}
-\frac{\beta}{\mu_M \sqrt{n}} + \left(x_1 + \frac{\beta}{\mu_M\sqrt{n}}\right)e^{-\mu_M \tau} -\frac{\mu_1{\kappa}}{(\mu_1-\mu_M)\sqrt{n}}\left(1-e^{-(\mu_M-\mu_1)\tau}\right)=0.
\end{equation}
When $\tau=0$, the left-hand side is equal to $x_1$, which is non-positive, and as $\tau\to\infty$ the left-hand side converges to $\infty$.  Furthermore, the derivative of the left-hand side of \eqref{eq:solution_uniqueness} with respect to $\tau$ is given by
\begin{align*}
    \left(-\mu_Mx_1 - \frac{\beta}{\sqrt{n}}+\frac{\mu_1{\kappa} e^{\mu_1\tau}}{\sqrt{n}}\right)e^{-\mu_M \tau}&= \left(-\mu_Mx_1 + \frac{\mu_1}{\sqrt{n}}\left({\kappa} e^{\mu_1\tau}-\frac{\beta}{\mu_1}\right)\right)e^{-\mu_M \tau},
\end{align*}
which is clearly non-negative for any ${\kappa} \geq \beta/\mu_1$. 
Hence, we can conclude that there exists a unique  $\tau^*\in [0,\infty)$ such that equation~\eqref{eq:solution_uniqueness} holds, which in turn implies the uniqueness of $x_2^*$. The second statement of the lemma follows directly by solving for $(x_1,\tau)$ as a function of $x_2$.
\qed
\endproof

{\section{Proof of Lemma~\ref{lem:tau_def_lem}}
\label{proof:tau_def_lem}

If $\mf x \geq \Gamma^{{\kappa}}$ for some ${\kappa} \geq \beta/\mu_1$, then by definition of $\Gamma^{{\kappa}}$ and Lemma~\ref{lem:tau_solution} there exist $\delta \geq 0$ and $\tau'=\tau^*(x_1,{\kappa}) \in [0,\infty)$ such that the following equations hold
\begin{align}
    \label{eq:tau_equation2}&-\frac{\beta}{\mu_M \sqrt{n}} + \left(x_1 + \frac{\beta}{\mu_M\sqrt{n}}\right)e^{-\mu_M \tau'} -\frac{\mu_1(x_2-\delta)}{\mu_1-\mu_M}\left(e^{-\mu_1\tau'}-e^{-\mu_M\tau'}\right)=0,\\
    &(x_2-\delta)e^{-\mu_1\tau'} ={\kappa}/\sqrt{n}.\label{eq:tau_equation12}
\end{align}
To simplify the notation, we define
\begin{equation}
    g_{\mf x}(\tau)=-\frac{\beta}{\mu_M \sqrt{n}} + \left(x_1 + \frac{\beta}{\mu_M\sqrt{n}}\right)e^{-\mu_M \tau} -\frac{\mu_1x_2}{\mu_1-\mu_M}\left(e^{-\mu_1\tau}-e^{-\mu_M\tau}\right)
    \label{eq:defg}
\end{equation}
Clearly, $g_{\mf x}(0)=x_1 \leq 0$ and using \eqref{eq:tau_equation2} we also have $g_{\mf x}(\tau')=-\frac{\mu_1\delta}{\mu_1-\mu_M}\left(e^{-\mu_1\tau'}-e^{-\mu_M\tau'}\right) \geq 0$. Hence, by the continuity of $g_{\mf x}(\cdot)$ there exists $\tau \in [0,\tau']$ where $g_{\mf x}(\tau)=0$. Furthermore, for such a $\tau \leq \tau'$ we also have $x_2e^{-\mu_1 \tau}\geq x_2e^{-\mu_1 \tau'}\geq {\kappa}/\sqrt{n}$. This proves the first part of the lemma since $\tau(\mf x)$ is defined as the smallest non-negative root of $g_{\mf x}(\cdot)$. 

To prove the second part, observe that 
$\mf x \in \Gamma^{{\kappa}}$ implies $\delta=0$ in equations~\eqref{eq:tau_equation2}-\eqref{eq:tau_equation12}. Hence, we have $g_{\mf x}(\tau')=0$ where $\tau'=\tau^*(x_1,{\kappa})$. 
Now, we shall prove that there exists no root of $g_{\mf x}(\cdot)$ in $[0,\tau')$ which will establish that $\tau(\mf x)=\tau'$, thereby completing the proof of this part. Let us suppose that it is not true. Hence, there exists at least one point $0\leq \tau_1 < \tau'$ such that $g_{\mf x}(\tau_1)=0$ and 
$g'_{\mf x}(\tau_1)g'_{\mf x}(\tau') \leq 0$ where $g'_{\mf x}(\cdot)$ denotes the derivative of $g_{\mf x}(\cdot)$ w.r.t $\tau$.
But from~\eqref{eq:defg} it is easy to see that $g'_{\mf x}(\tau)=-\mu_M g(\tau)+\mu_1 x_2 e^{-\mu_1 \tau}-\beta/\sqrt{n}$. Hence, $g'_{\mf x}(\tau_1)=\mu_1 x_2 e^{-\mu_1 \tau_1}-\beta/\sqrt{n}$ and $g'_{\mf x}(\tau')=\mu_1 x_2 e^{-\mu_1 \tau'}-\beta/\sqrt{n}$. 
Furthermore, $0\leq \tau_1 < \tau'$ implies that $x_2e^{-\mu_1 \tau_1}> x_2e^{-\mu_1 \tau'}= {\kappa}/\sqrt{n}> \beta/\mu_1\sqrt{n}$. Hence, $g'_{\mf x}(\tau_1)g'_{\mf x}(\tau_2) > 0$ which contradicts $g'_{\mf x}(\tau_1)g'_{\mf x}(\tau') \leq 0$. 


To prove the third part of the lemma, we take the partial derivatives of the equation $g_{\mf x}(\tau(\mf x))=0$ with respect to $x_1$ and $x_2$. This yields the expressions given by~\eqref{eq:tau1} and~\eqref{eq:tau2}, respectively, after manipulation using $g_{\mf x}(\tau(\mf x))=0$. The inequalities $\tau_1(\mf x)\leq 0 \text{ and } \tau_2(\mf x) \leq 0$ follow by observing that $x_2e^{-\mu_1 \tau(\mf x)}\geq {\kappa}/\sqrt{n}\geq \beta/\mu_1\sqrt{n}$.

To prove the fourth part of the lemma, we fix $\mf x$ and assume that $\mf x \geq \Gamma^{{\kappa}_1}$. From first part we know that $x_2e^{-\mu_1\tau(\mf x)} \geq {\kappa}_1/\sqrt{n}>{\kappa}_2/\sqrt{n}$. Taking derivative of $x_2e^{-\mu_1\tau(\mf x)}$ with respect to $x_2$ we get
\begin{equation*}
    \frac{d}{dx_2} x_2e^{-\mu_1\tau(\mf x)}=e^{-\mu_1\tau(\mf x)} (1-\mu_1 x_2 \tau_2(\mf x)) >0,
\end{equation*}
where the inequality follows from~\eqref{eq:tau2}. Therefore, $x_2e^{-\mu_1\tau(\mf x)}$ is non-decreasing in $x_2$ and we can write for any $\delta\geq0$ 
\begin{equation}
(x_2+\delta)e^{-\mu_1\tau(x_1,x_2+\delta)}\geq {\kappa}_1/\sqrt{n}>{\kappa}_2/\sqrt{n}\label{eq:violation}.
\end{equation}
The above inequality  implies that $(x_1,x_2+\delta) \notin \Gamma^{{\kappa}_2}$ for all $\delta\geq0$ since by Part 2 of this lemma $(x_1,x_2+\delta) \in \Gamma^{{\kappa}_2}$ will imply $\tau(x_1,x_2+\delta)=\tau^*(x_1,{\kappa}_2)$ which violates~\eqref{eq:violation}. On the other hand, from Lemma~\ref{lem:tau_solution}, there must exist $x_2'\geq0$ such that $(x_1,x_2')\in \Gamma^{{\kappa}_2}$, implies that $x_2'=x_2-\delta'$ for some $\delta'>0$. Hence, we have $\mf x >\Gamma^{{\kappa}_2}$.

To prove the last part of the lemma, let us assume $\tau(\mf x) <\infty$ for some $x < \Gamma^{\beta/\mu_1}$. Then we must have $\mu_1 x_2 e^{-\mu_1 \tau(\mf x)}<\beta/\sqrt{n}$ since otherwise we would have $\mu_1 x_2 e^{-\mu_1 \tau(\mf x)}=\kappa/\sqrt{n}$ for some $\kappa\geq \beta/\mu_1$ and combined with $g_{\mf x}(\tau(\mf x))=0$ it would imply that $\mf x\in \Gamma^{\kappa}$ for some $\kappa \geq \beta/\mu_1$ which contradicts with $x< \Gamma^{\beta/\mu_1}$ due to the fourth part of the lemma. Furthermore, since $\tau(\mf x)$ is the smallest non-negative root of $g_{\mf x}(\cdot)$ and $g_{\mf x}(0)=x_1\leq 0$ we must have $g_{\mf x}'(\tau(\mf x))\geq 0$. But from the proof of the previous part we know that  $g'_{\mf x}(\tau(\mf x))=-\mu_M g(\tau(\mf x))+\mu_1 x_2 e^{-\mu_1 \tau(\mf x)}-\beta/\sqrt{n}=\mu_1 x_2 e^{-\mu_1 \tau(\mf x)}-\beta/\sqrt{n}<0$ which contradicts with $g_{\mf x}'(\tau(\mf x))\geq 0$.
\qed
\endproof}

\section
{Proof of Lemma~\ref{lem:solution_PDE}}
\label{append:lem9}

In this section, we show that the Lyapunov function given by~\eqref{eq:lyapunov_function} solves the PDE \eqref{eq:PDE_Lyapunov1}-\eqref{eq:PDE_Lyapunov2}
and satisfies the properties stated in Lemma~\ref{lem:solution_PDE}.

For convenience, we repeat the expression of $f^*(\mf x)$ below
\begin{align}
    \label{eq:lyapunov_function_rep}
    f^*(\mf x) &=\left\{\begin{array}{ll}
    0, &\mbox{for } \mf x\in \Omega_1,\\
    \frac{x_2}{\mu_1}-\frac{\kappa}{\mu_1\sqrt{n}}-\frac{\kappa}{\mu_1\sqrt{n}}\log\left(\frac{\sqrt{n}x_2}{\kappa}\right), &\mbox{for } \mf x\in \Omega_2,\\
    \frac{x_2}{\mu_1}-\frac{x_2}{\mu_1}e^{-\mu_1\tau(\mf x)}-\frac{\kappa}{\sqrt{n}}\tau(\mf x)+\frac{\sqrt{n}}{2\beta}\left(x_2e^{-\mu_1\tau(\mf x)}-\frac{\kappa}{\sqrt{n}}\right)^2, &\mbox{for } \mf x\in \Omega_3.
    \end{array}\right.
\end{align}
where the sets $\Omega_i$, $i\in [3]$ are defined as
\[
\Omega_1=\{x_2\leq \kappa/\sqrt{n}\}, \Omega_2 = \{x_2>\kappa/\sqrt{n}\mbox{ and }\mf x\leq \Gamma^\kappa\}\mbox{ and }\Omega_3=\{\mf x>\Gamma^\kappa\}
\]
and by Lemma~\ref{lem:tau_def_lem} we have $(-\infty,0]\times[0,\infty)=\Omega=\Omega_1\cup\Omega_2\cup\Omega_3$. Also recall that the operator $\cL$ is defined as
\begin{align}\label{eq:PDE_operator_rep}
    \cL f(\mf x) &=\left(-\frac{\beta}{\sqrt{n}}-\mu_M x_1 + \mu_1 x_2\right)f_1(\mf x) - \mu_1 x_2 f_2(\mf x),
\end{align}
for any differentiable function $f$.

We first calculate the derivatives of $f^*$ on the sets $\Omega_1$, $\Omega_2$ and $\Omega_3$. When $\mf x\in \Omega_1$, we clearly have $f_1^*(\mf x)=f_2^*(\mf x)=0$. By~\eqref{eq:PDE_operator}, this implies that $\cL f^*(x)=0=-\left(x_2-\kappa/\sqrt{n}\right)_+$ when $\mf x\in \Omega_1$.

When $\mf x \in \Omega_2$, we have $f_1^*(\mf x)=0$ and 
\begin{align*}
    f_2^*(\mf x)&=\frac{1}{\mu_1}-\frac{\kappa}{\mu_1\sqrt{n}}\frac{\kappa}{\sqrt{n}x_2}\frac{\sqrt{n}}{\kappa}=\frac{1}{\mu_1}-\frac{\kappa}{\mu_1\sqrt{n}x_2},
\end{align*}
and therefore we have
\begin{align*}
    \cL f^*(\mf x)& = -\mu_1x_2f_2(\mf x)=-x_2+\frac{\kappa}{\sqrt{n}}=-\left(x_2-\kappa/\sqrt{n}\right)_+.
\end{align*}
When $\mf x \in \Omega_3$ we have
\begin{align}
    f_1^*(\mf x)&=\left(x_2e^{-\mu_1\tau(\mf x)}-\frac{\kappa}{\sqrt{n}}-\frac{\sqrt{n}}{\beta}\left(x_2e^{-\mu_1\tau(\mf x)}-\frac{\kappa}{\sqrt{n}}\right)x_2\mu_1e^{-\mu_1\tau(\mf x)}\right)\tau_1(\mf x),\nonumber\\
    &=-\frac{\sqrt{n}}{\beta}\tau_1(\mf x)\left(x_2\mu_1e^{-\mu_1\tau(\mf x)}-\frac{\beta}{\sqrt{n}}\right)\left(x_2e^{-\mu_1\tau(\mf x)}-\frac{\kappa}{\sqrt{n}}\right),\nonumber\\
    &=\frac{\sqrt{n}}{\beta}e^{-\mu_M\tau(\mf x)}\left(x_2e^{-\mu_1\tau(\mf x)}-\frac{\kappa}{\sqrt{n}}\right),\label{eq:f1}
\end{align}
where the last line follows by using the expression of $\tau_1(\mf x)$ given in Lemma~\ref{lem:tau_def_lem}. Furthermore, in the same region we have
\begin{align}
    f_2^*(\mf x) &= \frac{1}{\mu_1}\brac{1-e^{-\mu_1\tau(\mf x)}} + \brac{x_2e^{-\mu_1\tau(\mf x)}-\frac{\kappa}{\sqrt{n}}}\tau_2(\mf x)\nonumber\\
    &\quad +\frac{\sqrt{n}}{\beta}\left(x_2e^{-\mu_1\tau(\mf x)}-\frac{\kappa}{\sqrt{n}}\right)\left(e^{-\mu_1\tau(\mf x)}-\mu_1x_2e^{-\mu_1\tau(\mf x)}\tau_2(\mf x)\right)\nonumber\\
    &=\frac{1}{\mu_1}\brac{1-e^{-\mu_1\tau(\mf x)}}\nonumber\\
    &\quad+\brac{x_2e^{-\mu_1\tau(\mf x)}-\frac{\kappa}{\sqrt{n}}}\brac{\brac{1-\frac{\mu_1\sqrt{n}}{\beta}x_2e^{-\mu_1\tau(\mf x)}}\tau_2(\mf x)+\frac{\sqrt{n}}{\beta}e^{-\mu_1\tau(\mf x)}}\nonumber\\
    &=\frac{1}{\mu_1}\brac{1-e^{-\mu_1\tau(\mf x)}}\nonumber\\
    &\quad+\brac{x_2e^{-\mu_1\tau(\mf x)}-\frac{\kappa}{\sqrt{n}}}\frac{\sqrt{n}}{\beta}\brac{\frac{\mu_1}{\mu_1-\mu_M}\brac{e^{-\mu_M\tau(\mf x)}-e^{-\mu_1\tau(\mf x)}}+e^{-\mu_1\tau(\mf x)}}\label{eq:f2}
\end{align}
Also, by definition of $\tau(\mf x)$ we have 
\begin{align}
    -\frac{\beta}{\mu_M \sqrt{n}} + \left(x_1 + \frac{\beta}{\mu_M\sqrt{n}}\right)e^{-\mu_M \tau(\mf x)} -\frac{\mu_1x_2}{\mu_1-\mu_M}\left(e^{-\mu_1\tau(\mf x)}-e^{-\mu_M\tau(\mf x)}\right)=0\label{eq:tau}.
\end{align}
After some straightforward algebraic manipulations using~\eqref{eq:f1},~\eqref{eq:f2}, and~\eqref{eq:tau}
it follows that $\cL f^*(\mf x)=-x_2+\kappa/\sqrt{n}$ when $\mf x\in \Omega_3$. Hence, we have verified that 
\eqref{eq:PDE_Lyapunov1} is satisfied for all $\mf x \in \Omega$.

Now, we verify the boundary conditions for the derivatives. It is easy to verify from~\eqref{eq:tau} that when $x_1=0$ and $\mf x\in \Omega_3$
we have $\tau(\mf x)=0$. Hence, from~\eqref{eq:f1} and~\eqref{eq:f2} we have 
\begin{align}
    f_1^*(0,x_2)=f_2^*(0,x_2) = \frac{\sqrt{n}}{\beta}\left(x_2-\frac{\kappa}{\sqrt{n}}\right),
\end{align}
which shows that \eqref{eq:PDE_Lyapunov2} holds for all $\mf x\in \Omega$. 

Now, we show that the second-order derivatives of $f^*(\mf x)$ satisfy the bounds given in Lemma~\ref{lem:solution_PDE}. 
To summarize our development above, we have derived the following expressions for the first derivatives of $f^*$.
\begin{align}
\label{eq:lyapunov_x1_der}
    f_1^*(\mf x)=\left\{\begin{array}{ll}
    0 &\mbox{if }\mf x\in \Omega_1\cup \Omega_2\\
    \frac{\sqrt{n}}{\beta}e^{-\mu_M\tau(\mf x)}\left(x_2e^{-\mu_1\tau(\mf x)}-\frac{\kappa}{\sqrt{n}}\right) &\mbox{if } \mf \mf \mf x\in\Omega_3
    \end{array}\right.,\mbox{ and }
\end{align}
\begin{align}
\label{eq:lyapunov_x2_der}
    f_2^*(\mf x) &=\left\{\begin{array}{ll}
    0 &\mbox{if }\mf x\in \Omega_1\\
    \frac{1}{\mu_1}-\frac{\kappa}{\mu_1\sqrt{n}x_2} &\mbox{if }\mf x\in\Omega_2\\
    \frac{1}{\mu_1}\brac{1-e^{-\mu_1\tau(\mf x)}}+\brac{x_2e^{-\mu_1\tau(\mf x)}-\frac{\kappa}{\sqrt{n}}}\frac{\sqrt{n}}{\beta}\theta(\mf x)&\mbox{if }\mf x\in \Omega_3\\
    \end{array}\right.,
\end{align}
where
\begin{equation}
    \theta(\mf x)=\frac{\mu_1}{\mu_1-\mu_M}\brac{e^{-\mu_M\tau(\mf x)}-e^{-\mu_1\tau(\mf x)}}+e^{-\mu_1\tau(\mf x)}\geq 0
\end{equation}

Taking the derivative of \eqref{eq:lyapunov_x1_der} with respect to $x_1$ for $\mf x\in \Omega_3$
\begin{align*}
    f_{11}^*(\mf x) &=\frac{\sqrt{n}(\mu_1+\mu_M)}{\beta}\tau_1(\mf x)e^{-\mu_M \tau(\mf x)}\left(-x_2e^{-\mu_1\tau(\mf x)}+\frac{\kappa}{\sqrt{n}}\frac{\mu_M}{\mu_1+\mu_M}\right)
\end{align*}
Now, $\mf x> \Gamma^\kappa$ implies that $ x_2e^{-\mu_1\tau(\mf x)} \geq \kappa/\sqrt{n}$ which in turn implies that the term inside the parenthesis in the above expression is negative. Furthermore, since by Lemma~\ref{lem:tau_def_lem} we have $\tau_1(\mf x)\leq 0$, we conclude that $f_{11}^*(\mf x)\geq 0$ for all $\mf x$. Plugging in the expression of $\tau_1(\mf x)$ from Lemma~\ref{lem:tau_def_lem} in the above equation we obtain
\begin{align*}
    f_{11}^*(\mf x) &=\frac{\sqrt{n}(\mu_1+\mu_M)}{\beta}\frac{e^{-2\mu_M\tau(\mf x)}}{\mu_1x_2e^{-\mu_1\tau(\mf x)}-\frac{\beta}{\sqrt{n}}}\left(x_2e^{-\mu_1\tau(\mf x)}-\frac{\kappa}{\sqrt{n}}\frac{\mu_M}{\mu_1+\mu_M}\right),\\
    &\leq \frac{\sqrt{n}(\mu_1+\mu_M)}{\beta}\frac{e^{-2\mu_M\tau(\mf x)}x_2e^{-\mu_1\tau(\mf x)}}{\mu_1x_2e^{-\mu_1\tau(\mf x)}-\frac{\beta}{\sqrt{n}}},\\
    &\leq \frac{\sqrt{n}(\mu_1+\mu_M)}{\mu_1\beta}\frac{x_2e^{-\mu_1\tau(\mf x)}}{x_2e^{-\mu_1\tau(\mf x)}-\frac{\beta}{\mu_1\sqrt{n}}},\\
    &\leq \frac{\sqrt{n}}{\beta}\frac{\mu_1+\mu_M}{\mu_1}\frac{\kappa}{\kappa-\frac{\beta}{\mu_1}},
\end{align*}
where the last inequality follows since $ x_2e^{-\mu_1\tau(\mf x)} \geq \kappa/\sqrt{n}$ for $\mf x > \Gamma^\kappa$ and $z\mapsto\frac{z}{z-c}$ is a non-increasing function for any constant $c\geq 0$.

Now, taking derivative of $f_1^*(\mf x)$ with respect to $x_2$ for $\mf x\in \Omega_3$ we obtain
\begin{align*}
    f_{12}^*(\mf x) &=-\frac{\sqrt{n}(\mu_1+\mu_M)}{\beta}e^{-\mu_M\tau(\mf x)}\tau_2(\mf x)\left(x_2e^{-\mu_1\tau(\mf x)}-\frac{\kappa}{\sqrt{n}}\frac{\mu_M}{\mu_1+\mu_M}\right).
\end{align*}
As $\tau_2(\mf x)\leq 0$ (by Lemma~\ref{lem:tau_def_lem}) and $\mf x> \Gamma_\kappa$, it follows that $f_{12}^*(\mf x)\geq 0$.

Now, we concentrate on $f_{22}^*(\mf x)$ which is trivially zero when $\mf x\in \Omega_1$. When $\mf x \in \Omega_2$ from~\eqref{eq:lyapunov_x2_der} we have
\[
f_{22}^*(\mf x) = \frac{\kappa}{\mu_1\sqrt{n}x_2^2}\leq \frac{\sqrt{n}}{\mu_1\kappa}
\]
and hence $f_{22}^*(\mf x)\in \sbrac{0,\frac{\sqrt{n}}{\mu_1\kappa}}$.
When $\mf x\in \Omega_3$ again from~\eqref{eq:lyapunov_x2_der} we have
\begin{align}
    f_{22}^*(\mf x)&=e^{-\mu_1\tau(\mf x)}\brac{\tau_2(\mf x)\brac{1-\frac{\sqrt{n}}{\beta}\mu_1x_2\theta(\mf x)}+\frac{\sqrt{n}}{\beta}\theta(\mf x)}\nonumber\\
    &\quad\quad+\frac{\sqrt{n}}{\beta}\brac{x_2e^{-\mu_1\tau(\mf x)}-\frac{\kappa}{\sqrt{n}}}\theta_2(\mf x),\label{eq:f22}
\end{align}
where $\theta_2(\mf x)$ is the derivative of $\theta(\mf x)$ with respect to $x_2$ and is given by
\begin{equation}
    \theta_2(\mf x)=-\frac{\mu_1\mu_M}{\mu_1-\mu_M}\brac{e^{-\mu_M\tau(\mf x)}-e^{-\mu_1\tau(\mf x)}}\tau_2(\mf x),\label{eq:theta2}
\end{equation}
for $\mf x \in \Omega_3$.
Since $\tau_2(\mf x) \leq 0$, we have $\theta_2(\mf x) \geq 0$. Furthermore, by the definition of $\theta(\mf x)$ we have $\theta(\mf x)\geq e^{-\mu_1\tau(\mf x)}$. Now, since $\mf x \in \Omega_3$ we have $x_2\theta(\mf x)\geq x_2e^{-\mu_1 \tau(\mf x)}\geq \kappa/\sqrt{n}\geq \beta/\mu_1\sqrt{n}$. In other words, we have $1-\frac{\sqrt{n}}{\beta}\mu_1x_2\theta(\mf x) \leq 0$ which implies that $f_{22}^*(\mf x) \geq 0$ for $\mf x\in \Omega_3$. Hence, we have established that $f_{22}^*(\mf x)\geq 0$ for all $\mf x \in \Omega$.

Now we proceed to find an upper bound on $f_{22}^*(\mf x)$ on $\Omega_3$. To do this, we first note that $$\theta(\mf x)= \frac{\mu_1}{\mu_1-\mu_M}e^{-\mu_M\tau(\mf x)}-\frac{\mu_M}{\mu_1-\mu_M}e^{-\mu_1\tau(\mf x)}\leq \frac{\mu_1}{\mu_1-\mu_M}.$$
Now, plugging in the expression of $\tau_2(\mf x)$ from Lemma~\ref{lem:tau_def_lem} into~\eqref{eq:theta2} we obtain $$\theta_2(\mf x)=\frac{\mu_1\mu_M}{\mu_1-\mu_M}\frac{\brac{e^{-\mu_M\tau(\mf x)}-e^{-\mu_1\tau(\mf x)}}^2}{x_2e^{-\mu_1\tau(\mf x)}-\frac{\beta}{\mu_1\sqrt{n}}}\leq \frac{\mu_1\mu_M}{\mu_1-\mu_M}\frac{1}{x_2e^{-\mu_1\tau(\mf x)}-\frac{\beta}{\mu_1\sqrt{n}}}.$$
From the above and the fact that $\kappa>\beta/\mu_1$, the second term in~\eqref{eq:f22} can be bounded above by $$\frac{\sqrt{n}}{\beta}\frac{\mu_1\mu_M}{\mu_1-\mu_M}.$$
Using the expression of $\tau_2(\mf x)$ from Lemma~\ref{lem:tau_def_lem} again we have
\begin{align}
  e^{-\mu_1\tau(\mf x)}\tau_2(\mf x)\brac{1-\frac{\sqrt{n}}{\beta}\mu_1x_2\theta(\mf x)}&=\frac{e^{-\mu_1\tau(\mf x)}}{\mu_1-\mu_M}\frac{\left(e^{-\mu_M\tau(\mf x)}-e^{-\mu_1\tau(\mf x)}\right)}{\left(x_2e^{-\mu_1\tau(\mf x)}-\frac{\beta}{\mu_1\sqrt{n}}\right)}\brac{\frac{\sqrt{n}}{\beta}\mu_1x_2\theta(\mf x)-1}\nonumber\\
  &\leq \frac{\sqrt{n}}{\beta}\brac{\frac{\mu_1}{\mu_1-\mu_M}}^2\frac{\left(e^{-\mu_M\tau(\mf x)}-e^{-\mu_1\tau(\mf x)}\right)}{\left(x_2e^{-\mu_1\tau(\mf x)}-\frac{\beta}{\mu_1\sqrt{n}}\right)}x_2e^{-\mu_1\tau(\mf x)}\nonumber\\
  &\leq \frac{\sqrt{n}}{\beta}\brac{\frac{\mu_1}{\mu_1-\mu_M}}^2\frac{x_2e^{-\mu_1\tau(\mf x)}}{\left(x_2e^{-\mu_1\tau(\mf x)}-\frac{\beta}{\mu_1\sqrt{n}}\right)}\nonumber\\
  &\leq \frac{\sqrt{n}}{\beta}\brac{\frac{\mu_1}{\mu_1-\mu_M}}^2\frac{\kappa}{\kappa-\frac{\beta}{\mu_1}}
\end{align}
Hence, combining the bounds we have
\begin{equation}
f_{22}^*(\mf x)\leq \frac{\sqrt{n}}{\beta}\brac{\brac{\frac{\mu_1}{\mu_1-\mu_M}}^2\frac{\kappa}{\kappa-\frac{\beta}{\mu_1}}+\frac{\mu_1(1+\mu_M)}{\mu_1-\mu_M}},
\end{equation}
for $\mf x\in \Omega_3$. The bound on $f_{22}^*(\mf x)$ stated in Lemma~\ref{lem:solution_PDE} now follows by combining the bounds for the sets $\Omega_2$ and $\Omega_3$.

The absolute continuity of $f_1^*(\cdot, x_2)$ and $f_2^*(x_1,\cdot)$ follows using the same argument as in \cite{Braverman2020}. 
\qed
\endproof

\section{Positive Recurrence of The Limiting Diffusion: Proof of Theorem~\ref{thm:stationarity_diffusion}}
\label{append:4}

The diffusion limit is expressed as 
\begin{align}
 Y_{M,1}(t)&=Y_{1,M}(0) + \sqrt{2}W(t)-\beta t+\int_0^t (\mu_1Y_{2,1}(s)-\mu_MY_{M,1}(s))ds-V_1(t),\\
Y_{1,2}(t)&=Y_{1,2}(0)-\mu_1\int^t_0 Y_{1,2}(s)ds+V_1(t),   
\end{align}

Using Ito's lemma, for any $f(\mf x) \in C^2(\Omega)$
\begin{align*}
    &\EE_x[f(Y_{M,1}(t), Y_{1,2}(t))]-E_x[f(Y_{M,1}(0), Y_{1,2}(0))]\\
    &=\EE_x\left[\int_0^t\left(\mu_1Y_{1,2}(s)-\mu_MY_{M,1}(s)-\beta\right)f_1(Y_{M,1}(s), Y_{1,2}(s))-\mu_1Y_{1,2}(s)f_1(Y_{M,1}(s), Y_{1,2}(s))ds\right]\\
    &\quad +\EE_x\left[\int_0^tf_{1,1}(Y_{M,1}(s), Y_{1,2}(s))ds\right] + \EE_x\left[\int_0^t\left(f_1(0, Y_{1,2}(s))-f_2(0, Y_{1,2}(s))\right)dU(s)\right].
\end{align*}
So for any function $f(\mf x)$ with $f_1(0,x_2)=f_2(0,x_2)$, the extended generator becomes
\begin{align*}
    \cG_Yf(\mf x) = (-\beta-\mu_Mx_1+\mu_1x_2)f_1(\mf x) -\mu_1x_2f_2(\mf x)+f_{11}(\mf x), \mf x\in \Omega.
\end{align*}

The geometric ergodicity of the diffusion will follow if we can find a suitable Lyapunov function $V(\mf x)$, constants $c,d>0$ and a compact set $K$ such that 
\[
\cG_YV(\mf x) \leq -cV(\mf x) +d\indic{x\in K}.
\]
In the rest of this section, our goal is to prove the theorem by identifying the function $V(\mf x)$, the compact set $K$, and the constants $c$ and $d$ following the same steps as \cite{Braverman2020}. The compact set $K$ has the form 
\begin{align*}
    K=[-\kappa_2,0]\times [0, \kappa_2],
\end{align*}
for some $\kappa_2>\beta$. Unfortunately, the indicator function $\indic{ \mf x\notin K}$ is not smooth and hence, we need to define a function $\phi(\mf x)\geq \indic{\mf x\notin K}$ for all $\mf x\in \Omega$ with regularity properties to approximate the indicator function. The function we need has the form
\begin{align*}
    \phi(\mf x) = \left\{\begin{array}{ll}
    0 &\mbox{for }\mf x\leq \kappa_1\\
    (\mf x-\kappa_1)^2\left(\frac{-(x-\kappa_1)}{((\kappa_2+\kappa_1)/2-\kappa_1)^2(\kappa_2-\kappa_1)}+\frac{2}{((\kappa_2+\kappa_1)/2-\kappa_1)(\kappa_2-\kappa_1)}\right) &\mbox{for }\mf x\in [\kappa_1,(\kappa_2+\kappa_1)/2]\\
    1-(\mf x-\kappa_2)^2\left(\frac{-(\mf x-\kappa_1)}{((\kappa_2+\kappa_1)/2-\kappa_1)^2(\kappa_2-\kappa_1)}-\frac{2}{((\kappa_2+\kappa_1)/2-\kappa_1)(\kappa_2-\kappa_1)}\right)&\mbox{for } \mf x\in [(\kappa_2+\kappa_1)/2,\kappa_2]\\
    1 &\mbox{for }\mf x\geq \kappa_2,
    \end{array}\right. 
\end{align*}
for some $\kappa_2>\kappa_1>\beta$. We need the first and second derivatives of this function and as shown in \cite{Braverman2020}, we have
\begin{align*}
    &\phi'(\kappa_1)=\phi'(\kappa_2)=0,\\
    &|\phi'(x)|\leq \frac{4}{\kappa_2-\kappa_1}\mbox{ and }|\phi''(x)|\leq \frac{12}{(\kappa_2-\kappa_1)^2}.
\end{align*}
We choose the Lyapunov function $V(\mf x)$ of the form 
\begin{align}\label{eq:Lyapunov_stationary}
    V(\mf x) &= \exp\left(\alpha(f^{(1)}(\mf x) + f^{(2)}(\mf x))\right),
\end{align}
where $\alpha>0$ and 
\begin{align}
    \label{eq:pde_stationary1}\cL_1f^{(1)}(\mf x)&:=(-\beta -\mu_Mx_1 +\mu_1x_2)f_1^{(1)}(x)-\mu_1x_2f_2^{(1)}(\mf x)=-\phi(-x_1), \mf x\in \Omega,\\
    \nonumber f_1^{(1)}(0,x_2)&=f_2^{(1)}(0,x_2), x_2\leq 0\\
    \label{eq:pde_stationary2}\cL_1f^{(2)}(\mf x)&=-\phi(x_2), \mf x\in \Omega,\\
    \nonumber f_1^{(2)}(0,x_2)&=f_2^{(2)}(0,x_2), x_2\leq 0.
\end{align}
It is easy to check $V_1(0,x_2)=V_2(0,x_2)$ for all $x_2\geq 0$. For convenience, we use the notation
\begin{align*}
 f^{(\Sigma)}(\mf x) &= f^{(1)}(\mf x)+ f^{(2)}(\mf x).   
\end{align*}
We need the following derivatives of $V(\mf x)$
\begin{align*}
    V_1(\mf x) &= \alpha f_1^{(\Sigma)}(\mf x)V(\mf x),\\
    V_2(\mf x) &=\alpha f_2^{(\Sigma)}(\mf x)V(\mf x),\\
    V_{11}(\mf x)&=\alpha\left(f_{11}^{(\Sigma)}(\mf x)\right)V(\mf x) + \alpha^2 \left(f_1^{(\Sigma)}(\mf x)\right)^2V(\mf x).\\
\end{align*}

Now, we have
\begin{align*}
    \cG V(\mf x)& = \left(\cL f^{(\Sigma)}(\mf x) \right)\alpha V(\mf x) + \frac{\alpha}{n} \left(f_{11}^{(\Sigma)}(\mf x)\right)V(\mf x) + \frac{\alpha^2}{n} \left(f_1^{(\Sigma)}(\mf x)(\mf x)\right)^2V(\mf x)\\
    &=(-\phi(-x_1)-\phi(x_2))\alpha V(\mf x) + \left(\alpha\left(f_{11}^{(\Sigma)}(\mf x)\right) + \alpha^2 \left(f_1^{(\Sigma)}(x)\right)^2\right)V(\mf x)\\
    &\leq -\alpha V(\mf x) 1(\mf x\notin K) + \left(\alpha\left(f_{11}^{(\Sigma)}(\mf x)\right) + \alpha^2 \left(f_1^{(\Sigma)}(\mf x)\right)^2\right)V(\mf x)\\
    &=\alpha\left(-1 + f_{11}^{(\Sigma)}(\mf x) + \alpha \left(f_1^{(\Sigma)}(x)\right)^2\right)V(\mf x)\\
    &\quad+ \left(\alpha\left(f_{11}^{(\Sigma)}(\mf x)\right) + \alpha^2 \left(f_1^{(\Sigma)}(\mf x)\right)^2\right)V(\mf x)1(\mf x\in K)
\end{align*}

Our goal now is to find $f^{(1)}(\mf x)$ and $f^{(2)}(\mf x)$ that solve the two PDEs in~\ref{eq:pde_stationary1} and $\ref{eq:pde_stationary2}$, and then show that $\kappa_1$ and $\kappa_2$ can be chosen so that \begin{align}
   \label{eq:drift_bound} &\alpha\left(-1 + f_{11}^{(\Sigma)}(\mf x) + \alpha \left(f_1^{(\Sigma)}(\mf x)\right)^2\right) \leq -c, \mbox{ for all }\mf x\in \Omega\\
    \label{eq:drift_bound_conditional}&\left(\alpha\left(f_{11}^{(\Sigma)}(\mf x)\right) + \alpha^2 \left(f_1^{(\Sigma)}(\mf x)\right)^2\right)V(x)\leq d \mbox{ for all } \mf x\in K.
\end{align} 
\begin{lemma}\label{lem:stationary_pde_solution}
For any $\kappa_2>\kappa_1>\beta$, there exist solutions $f^{(1)}(\mf x), f^{(2)}\in C^2(\Omega)$ that solves~\ref{eq:pde_stationary1} and \ref{eq:pde_stationary2}. Moreover, $\kappa_1$ and $\kappa_2$ can be chosen so that
\begin{align}
\label{eq:f1_bounds_stationary}&|f^{(1)}(\mf x)|\leq \frac{1}{\mu_M}\log 2 \qquad\qquad \mf x\in [-\kappa_2, 0]\times[0,\kappa_2]\\
\label{eq:f2_bounds_stationary}&|f^{(2)}(\mf x)|\leq \frac{\log2}{\mu_1} + \frac{\epsilon}{\beta} \qquad\quad\; \mf x\in [-\kappa_2, 0]\times[0,\kappa_2]\\
\label{eq:f1_der_bounds_stationary}&|f_1^{(1)}(\mf x)| \leq \frac{4\log2}{\epsilon},\qquad\qquad |f_{11}^{(1)}(\mf x)|\leq \frac{12\log2}{\mu_M\epsilon^2}\\
\label{eq:f2_der_bounds_stationary}&|f_1^{(2)}(\mf x)| \leq \frac{1}{\beta}, \qquad\qquad\qquad |f_{11}^{(2)}(\mf x)|\leq \frac{\mu_M+4}{\beta\epsilon}.
\end{align}
\end{lemma}
\proof{Proof.}
    For $\kappa_2>\kappa_1>\beta$, we first identify functions $f^{(1)}(\mf x)$ and $f^{(2)}(\mf x)$ that solves~\ref{eq:pde_stationary1} and \ref{eq:pde_stationary2} as 
\begin{align*}
    f^{(1)}(\mf x) = \int_{0}^\infty \phi(-v_1^x(t))dt
\end{align*}
\begin{align*}
    f^{(2)}(\mf x) = \int_{0}^\infty \phi(v_2^x(t))dt,
\end{align*} 
where $v_1^x(t)$ and $v_2^x(t)$ are as in Section~\ref{sec:lyapunov_func} with $n=1$. To evaluate the first integral, we need $t\geq 0$ that solves  
\begin{align*}
    \frac{\beta}{\mu_M}-\kappa &= \left(x_1 + \frac{\beta}{\mu_M} +\frac{\mu_1x_2}{\mu_1-\mu_M}\left(1-e^{-(\mu_1-\mu_M)t}\right)\right)e^{-\mu_M t}
\end{align*}
There exists a solution to this equation when $x_1<-\kappa$, as the right hand side is $x_1$ for $t=0$ and 0 as $t\to \infty$. However, the solution may not be unique. Hence, we define 
\begin{align*}
    \tau^{\kappa}(\mf x) &= \min\left\{t\geq 0: \frac{\beta}{\mu_M }-\kappa \leq \left(x_1 + \frac{\beta}{\mu_M} +\frac{\mu_1x_2}{\mu_1-\mu_M}\left(1-e^{-(\mu_1-\mu_M)t}\right)\right)e^{-\mu_M t}\right\}.
\end{align*}
The continuity of the right-hand side implies that for any $x_1<-\kappa_2<-\kappa_1<0$, we have $\tau^{\kappa_2}(\mf x)<\tau^{\kappa_1}(\mf x)$. We also have
\begin{align}
    \nonumber\frac{\mu_M\kappa_2-\beta}{\mu_M\kappa_1-\beta} &= \frac{\left(x_1 + \frac{\beta}{\mu_M} +\frac{\mu_1x_2}{\mu_1-\mu_M}\left(1-e^{-(\mu_1-\mu_M)\tau^{\kappa_2}(\mf x)}\right)\right)e^{-\mu_M \tau^{\kappa_2}(\mf x)}}{\left(x_1 + \frac{\beta}{\mu_M} +\frac{\mu_1x_2}{\mu_1-\mu_M}\left(1-e^{-(\mu_1-\mu_M)\tau^{\kappa_1}(\mf x)}\right)\right)e^{-\mu_M \tau^{\kappa_1}(\mf x)}}\\
    \nonumber\frac{\mu_M\kappa_2-\beta}{\mu_M\kappa_1-\beta} &\geq e^{\mu_M(\tau^{\kappa_1}(\mf x)-\tau^{\kappa_2}(\mf x))}\\
    \frac{1}{\mu_M}\log\left(\frac{\mu_M\kappa_2-\beta}{\mu_M\kappa_1-\beta}\right)&\geq \tau^{\kappa_1}(\mf x)-\tau^{\kappa_2}(\mf x).
\end{align}
Now, we can write
\begin{align*}
    f^{(1)}(\mf x) = \left\{\begin{array}{ll}
    \tau^{\kappa_2}(\mf x) +\int_{\tau^{\kappa_2}}^{\tau^{\kappa_1}} \phi\left(\frac{\beta}{\mu_M }- \left(x_1 + \frac{\beta}{\mu_M} -\frac{\mu_1x_2}{\mu_1-\mu_M}\left(1-e^{-(\mu_1-\mu_M)t}\right)\right)e^{-\mu_M t}\right)dt\\
    \qquad\qquad\qquad\qquad\qquad\qquad\qquad\qquad\qquad\qquad\qquad\qquad\qquad\mbox{ for } x_1<-\kappa_2\\
    \int_{0}^{\tau^{\kappa_1}} \phi\left(\frac{\beta}{\mu_M}- \left(x_1 + \frac{\beta}{\mu_M} -\frac{\mu_1x_2}{\mu_1-\mu_M}\left(1-e^{-(\mu_1-\mu_M)t}\right)\right)e^{-\mu_M t}\right)dt\\
    \qquad\qquad\qquad\qquad\qquad\qquad\qquad\qquad\qquad\qquad\qquad\qquad\qquad\mbox{ for } x_1\in[-\kappa_2,-\kappa_1]\\
    0 \qquad\qquad\qquad\qquad\qquad\qquad\qquad\qquad\qquad\qquad\qquad\qquad\quad\:\:\:\mbox{ for } x_1\in[-\kappa_1,0].
    \end{array}\right.
\end{align*}
We have $f_1^{(1)}(0,x_2)=f_2^{(1)}(0,x_2)=0$, and hence, $f^{(1)}(\mf x)$ satisfies the boundary condition. To check whether the suggested $f^{(1)}(\mf x)$ satisfies the PDE, we now calculate its derivatives. First, consider the derivative with respect to $x_1$ when $x_1<-\kappa_2$.
\begin{align*}
    f_1^{(1)}(\mf x) &= \tau^{\kappa_2}_1(\mf x)\\
    &\quad + \phi\left(\frac{\beta}{\mu_M }- \left(x_1 + \frac{\beta}{\mu_M} -\frac{\mu_1x_2}{\mu_1-\mu_M}\left(1-e^{-(\mu_1-\mu_M)\tau^{\kappa_1}}(x)\right)\right)e^{-\mu_M \tau_{\kappa_1}(\mf x)}\right)\tau_1^{\kappa_1}(\mf x)\\
    &\quad - \phi\left(\frac{\beta}{\mu_M }- \left(x_1 + \frac{\beta}{\mu_M} -\frac{\mu_1x_2}{\mu_1-\mu_M}\left(1-e^{-(\mu_1-\mu_M)\tau^{\kappa_2}}(x)\right)\right)e^{-\mu_M \tau_{\kappa_2}(\mf x)}\right)\tau_1^{\kappa_2}(\mf x)\\
    &\quad + \int_{\tau^{\kappa_2}}^{\tau^{\kappa_1}} e^{-\mu_Mt}\phi'\left(\frac{\beta}{\mu_M}- \left(x_1 + \frac{\beta}{\mu_M} -\frac{\mu_1x_2}{\mu_1-\mu_M}\left(1-e^{-(\mu_1-\mu_M)t}\right)\right)e^{-\mu_M t}\right)dt\\
    &= \int_{\tau^{\kappa_2}}^{\tau^{\kappa_1}} e^{-\mu_Mt}\phi'\left(\frac{\beta}{\mu_M}- \left(x_1 + \frac{\beta}{\mu_M} -\frac{\mu_1x_2}{\mu_1-\mu_M}\left(1-e^{-(\mu_1-\mu_M)t}\right)\right)e^{-\mu_M t}\right)dt.
\end{align*}
Note that the boundary terms reduce to $\phi(\kappa_1)=0$ and $\phi(\kappa_2)=1$. Now, considering the case $x_1\in[-\kappa_2,-\kappa_1]$, and using the same trick for the boundary term, we get 
\begin{align*}
    f_1^{(1)}(\mf x) &= \int_{0}^{\tau^{\kappa_1}} e^{-\mu_Mt}\phi'\left(\frac{\beta}{\mu_M}- \left(x_1 + \frac{\beta}{\mu_M} -\frac{\mu_1x_2}{\mu_1-\mu_M}\left(1-e^{-(\mu_1-\mu_M)t}\right)\right)e^{-\mu_M t}\right)dt.
\end{align*}
Combining, the results we have
\begin{align*}
    f_1^{(1)}(\mf x) &= \left\{
    \begin{array}{ll}
    \int_{\tau^{\kappa_2}}^{\tau^{\kappa_1}} e^{-\mu_Mt}\phi'\left(\frac{\beta}{\mu_M}- \left(x_1 + \frac{\beta}{\mu_M} -\frac{\mu_1x_2}{\mu_1-\mu_M}\left(1-e^{-(\mu_1-\mu_M)t}\right)\right)e^{-\mu_M t}\right)dt,\\\qquad\qquad\qquad\qquad\qquad\qquad\qquad\qquad\qquad\qquad\qquad\qquad\qquad\mbox{ for } x_1<-\kappa_2\\
     \int_{0}^{\tau^{\kappa_1}}  e^{-\mu_Mt}\phi'\left(\frac{\beta}{\mu_M}- \left(x_1 + \frac{\beta}{\mu_M} -\frac{\mu_1x_2}{\mu_1-\mu_M}\left(1-e^{-(\mu_1-\mu_M)t}\right)\right)e^{-\mu_M t}\right)dt,\\
     \qquad\qquad\qquad\qquad\qquad\qquad\qquad\qquad\qquad\qquad\qquad\qquad\qquad\mbox{ for } x_1\in[-\kappa_2,-\kappa_1]\\
    \end{array}\right..
\end{align*}
Using a similar argument, we get
\begin{align*}
    f_2^{(1)}(\mf x)&= \left\{\begin{array}{l}
    \int_{\tau^{\kappa_2}}^{\tau^{\kappa_1}} \frac{\mu_1}{\mu_1-\mu_M}\left(e^{-\mu_Mt}-e^{-\mu_1t}\right)\\
    \qquad\times\phi'\left(\frac{\beta}{\mu_M}- \left(x_1 + \frac{\beta}{\mu_M} -\frac{\mu_1}{\mu_1-\mu_M}\left(1-e^{-(\mu_1-\mu_M)t}\right)\right)e^{-\mu_M t}\right)dt,
    \qquad\;\mbox{ for }x_1<-\kappa_2\\
    \int_{0}^{\tau^{\kappa_1}} \frac{\mu_1x_2}{\mu_1-\mu_M}\left(e^{-\mu_Mt}-e^{-\mu_1t}\right)\\
    \qquad\times\phi'\left(\frac{\beta}{\mu_M }- \left(x_1 + \frac{\beta}{\mu_M} -\frac{\mu_1x_2}{\mu_1-\mu_M}\left(1-e^{-(\mu_1-\mu_M)t}\right)\right)e^{-\mu_M t}\right)dt,
    \qquad\mbox{ for } x_1\in[-\kappa_2,-\kappa_1]\\
    0,\qquad\qquad\qquad\qquad\qquad\qquad\qquad\qquad\qquad\qquad\qquad\qquad\qquad\qquad\quad\; \mbox{ for }x_1\in[-\kappa_1,0]
    \end{array}\right..
\end{align*}
Now, we can show that for $x_1<-\kappa_2$
\begin{align*}
    &\left(-\beta-\mu_M x_1 + \mu_1 x_2\right)f_1^{(1)}(\mf x) + \mu_1 x_2 f_2^{(1)}(\mf x)\\
    &\qquad\qquad = \int_{\tau^{\kappa_2}}^{\tau^{\kappa_1}}\left(\left(-\beta-\mu_M x_1 + \mu_1 x_2\right)e^{-\mu_Mt}-\mu_1x_2\frac{\mu_1}{\mu_1-\mu_M}\left(e^{-\mu_Mt}-e^{-\mu_1t}\right)\right)\\&\qquad\qquad\qquad\qquad \times\phi'\left(\frac{\beta}{\mu_M}- \left(x_1 + \frac{\beta}{\mu_M} -\frac{\mu_1}{\mu_1-\mu_M}\left(1-e^{-(\mu_1-\mu_M)t}\right)\right)e^{-\mu_M t}\right)dt\\
    &\qquad\qquad = \int_{\tau^{\kappa_2}}^{\tau^{\kappa_1}}\left(\left(-\beta-\mu_M x_1 \right)e^{-\mu_Mt}-\frac{\mu_1\mu_Mx_2}{\mu_1-\mu_M}e^{-\mu_Mt}+\frac{\mu_1^2x_2}{\mu_1-\mu_M}e^{-\mu_1t}\right)\\&\qquad\qquad\qquad\qquad \times\phi'\left(\frac{\beta}{\mu_M}- \left(x_1 + \frac{\beta}{\mu_M} -\frac{\mu_1}{\mu_1-\mu_M}\left(1-e^{-(\mu_1-\mu_M)t}\right)\right)e^{-\mu_M t}\right)dt\\
\end{align*}
Now, making the variable change
\begin{align*}
    u&=\frac{\beta}{\mu_M}- x_1e^{-\mu_Mt} + \frac{\beta}{\mu_M}e^{-\mu_Mt} -\frac{\mu_1}{\mu_1-\mu_M}e^{-\mu_Mt}-\frac{\mu_1}{\mu_1-\mu_M}e^{-\mu_1t},\\
    du&=\left(\beta+\mu_M x_1 \right)e^{-\mu_Mt}+\frac{\mu_1\mu_Mx_2}{\mu_1-\mu_M}e^{-\mu_Mt}-\frac{\mu_1^2x_2}{\mu_1-\mu_M}e^{-\mu_1t}dt,
\end{align*}
we have
\begin{align*}
    \left(-\beta-\mu_M x_1 + \mu_1 x_2\right)f_1^{(1)}(\mf x) + \mu_1 x_2 f_2^{(1)}(\mf x)&= \int_{\tau^{\kappa_2}}^{\tau^{\kappa_1}}-\phi'(u)du = -1 = -\phi(-x_1).
\end{align*}
A similar result holds for all other cases, hence $f^{(1)}(\mf x)$ satisfy the PDE along with the boundary conditions.
To prove a similar result for $f^{(2)}(\mf x)$ we will again follow the analysis in Braverman closely. First, we define four subdomains 
\begin{align*}
    S_0 &= \{\mf x\in \Omega: x_2\leq \kappa_1\}\\
    S_1 &= \{\mf x\in \Omega: x_2\geq \kappa_1, \mf x\leq \Gamma_{\kappa_1}\}\\
    S_2 &= \{\mf x\in \Omega: x_2\geq \kappa_1, \Gamma_{\kappa_1}\geq \mf x\geq \Gamma_{\kappa_2}\}\\
    S_3 &= \{\mf x\in \Omega: x_2\geq \kappa_1, \mf x\geq \Gamma_{\kappa_2}\}\\
\end{align*}
Then, we can define the function as
\begin{align*}
    f^{(2)}(\mf x)&=\left\{
    \begin{array}{ll}
    0, &\mbox{for }\mf x\in S_0\\
    \int_0^{\frac{\log(x_2/\kappa_1)}{\mu_1}}\phi(x_2e^{-\mu_1t}),  &\mbox{for }x_2\leq \kappa_2, \mf x\in S_1\\
    \frac{\log(x_2/\kappa_2)}{\mu_1} + \int_0^{\frac{\log(\kappa_2/\kappa_1)}{\mu_1}}\phi\left(\kappa_2 e^{-\mu_1t}\right)dt, &\mbox{for }x_2\geq \kappa_2, \mf x\in S_1\\
    \int_0^{\tau(\mf x)}\phi(x_2e^{-\mu_1t})dt + \frac{\sqrt{n}}{\beta}\int_{\kappa_1}^{x_2e^{-\mu_1\tau(\mf x)}}\phi(t)dt, &\mbox{for }x_2\leq \kappa_2, \mf x\in S_2\\
    \frac{\log(x_2/\kappa_2)}{\mu_1} + \int_{\frac{\log(x_2/\kappa_2)}{\mu_1}}^{\tau(\mf x)}\phi(x_2e^{-\mu_1t})dt
    + \frac{1}{\beta}\int_{\kappa_1}^{x_2e^{-\mu_1\tau(\mf x)}}\phi(t)dt, &\mbox{for }x_2\geq \kappa_2, \mf x\in S_2\\
    \tau(\mf x) + \frac{x_2e^{-\mu_1\tau(\mf x)}-\kappa_2}{\beta/} + \frac{1}{\beta}\int_{\kappa_1/\sqrt{n}}^{\kappa_2}\phi(t)dt, &\mbox{for }\mf x\in S_3.
    \end{array}\right.
\end{align*}

Now, we calculate the first derivatives with respect to $x_1$ and $x_2$
\begin{align*}
    f_1^{(2)}(\mf x) &=\left\{\begin{array}{ll}
    0, &\mbox{for } \mf x\in S_0\cup S_1\\
    \frac{1}{\beta}e^{-\mu_M\tau(\mf x)}\phi(x_2e^{-\mu_1\tau(\mf x) }), &\mbox{for }\mf x\in S_2\\
    \frac{1}{\beta}e^{-\mu_M\tau(\mf x)}, &\mbox{for }\mf x\in S_3
    \end{array}\right.\\
    f_2^{(2)}(\mf x) &=\left\{
    \begin{array}{ll}
    0, &\mbox{for }\mf x\in S_0\\
    \frac{1}{\mu_1x_2}\phi(x_2), &\mbox{for }\mf x\in S_1\\
    \frac{1}{\mu_1x_2}\left(\phi(x_2)-\phi(x_2e^{-\mu_1\tau(\mf x)}\right)+\frac{1}{\beta}\phi(x_2e^{-\mu_1\tau(\mf x)})\left(\frac{\mu_1e^{-\mu_M\tau(\mf x)}}{\mu_1-\mu_M}-\frac{\mu_Me^{-\mu_1\tau(\mf x)}}{\mu_1-\mu_M}\right), &\mbox{for }\mf x\in S_2\\
    \frac{1}{\beta}\frac{\mu_1e^{-\mu_M\tau(\mf x)}-\mu_Me^{-\mu_1\tau(\mf x)}}{\mu_1-\mu_M}, &\mbox{for }\mf x\in S_3
    \end{array}
    \right..
\end{align*}
Now, we prove that the suggested $f^{(2)}(\mf x)$ satisfies the PDE. 
\begin{enumerate}
    \item When $\mf x\in S_0$, $\cL f(\mf x) = 0=-\phi(0)$
    \item When $\mf x\in S_1$, 
    \begin{align*}
        &\left(-\beta-\mu_M x_1 + \mu_1 x_2\right)f_1^{(2)}(\mf x) + \mu_1 x_2 f_2^{(2)}(\mf x)=-\mu_1x_2f_2^{(2)}(\mf x)=-\phi(x_2)
    \end{align*}
    \item When $\mf x\in S_2$
    \begin{align*}
        &\left(-\beta-\mu_M x_1 + \mu_1 x_2\right)f_1^{(2)}(\mf x) + \mu_1 x_2 f_2^{(2)}(\mf x)\\
        &\qquad = -e^{-\mu_M\tau(\mf x)}\phi(x_2e^{-\mu_1\tau(\mf x)})-\frac{1}{\beta}\mu_Mx_1e^{-\mu_M\tau(\mf x)}\phi(x_2e^{-\mu_1\tau(\mf x)})\\
        &\qquad\quad + \frac{1}{\beta}\mu_1x_2e^{-\mu_M\tau(\mf x)}\phi(x_2e^{-\mu_1\tau(\mf x)}) -\phi(\mf x) +\phi(x_2e^{-\mu_1\tau(\mf x)}) \\
        &\qquad\quad- \frac{1}{\beta}\mu_1x_2\phi(x_2e^{-\mu_1\tau(\mf x)})\left(\frac{\mu_1e^{-\mu_M\tau(\mf x)}}{\mu_1-\mu_M}-\frac{\mu_Me^{-\mu_1\tau(\mf x)}}{\mu_1-\mu_M}\right)\\
        &\qquad =-\phi(\mf x) + \phi(x_2e^{-\mu_1\tau(\mf x)})\left(1 - e^{-\mu_M\tau(\mf x)}-\frac{1}{\beta}\mu_Mx_1e^{-\mu_M\tau(\mf x)} + \frac{1}{\beta}\frac{\mu_1^2x_2e^{-\mu_M\tau(\mf x)}}{\mu_1-\mu_M} \right.\\
        &\qquad\qquad\qquad\qquad\qquad\qquad\qquad\quad \left. - \frac{1}{\beta}\frac{\mu_1\mu_Mx_2e^{-\mu_M\tau(\mf x)}}{\mu_1-\mu_M}- \frac{1}{\beta}\frac{\mu_1^2x_2e^{-\mu_M\tau(\mf x)}}{\mu_1-\mu_M}\right.\\ &\qquad\qquad\qquad\qquad\qquad\qquad\qquad\qquad\qquad\qquad\qquad\qquad\qquad\quad \left.+\frac{1}{\beta}\frac{\mu_M\mu_1x_2e^{-\mu_1\tau(\mf x)}}{\mu_1-\mu_M}\right)\\
        &\qquad =-\phi(\mf x) + \phi(x_2e^{-\mu_1\tau(\mf x)})\left(1 - e^{-\mu_M\tau(\mf x)}-\frac{1}{\beta}\mu_Mx_1e^{-\mu_M\tau(\mf x)} \right.\\
        &\qquad\qquad\qquad\qquad\qquad\qquad\qquad\quad \left. - \frac{1}{\beta}\frac{\mu_1\mu_Mx_2e^{-\mu_M\tau(\mf x)}}{\mu_1-\mu_M}+\frac{1}{\beta}\frac{\mu_M\mu_1x_2e^{-\mu_1\tau(\mf x)}}{\mu_1-\mu_M}\right)\\
        &\qquad =-\phi(\mf x) + \phi(x_2e^{-\mu_1\tau(\mf x)})\left(1 - e^{-\mu_M\tau(\mf x)}\right.\\
        &\qquad\qquad\qquad\qquad\qquad\qquad\quad-\left(1-e^{-\mu_M\tau(\mf x)}-\frac{\mu_M}{\beta}\frac{\mu_1x_2e^{-\mu_M\tau(\mf x)}-\mu_1x_2e^{-\mu_1\tau(\mf x)}}{\mu_1-\mu_M}\right) \\
        &\qquad\qquad\qquad\qquad\qquad\qquad\quad \left. - \frac{1}{\beta}\frac{\mu_1\mu_Mx_2e^{-\mu_M\tau(\mf x)}}{\mu_1-\mu_M}+\frac{1}{\beta}\frac{\mu_M\mu_1x_2e^{-\mu_1\tau(\mf x)}}{\mu_1-\mu_M}\right)\\
        &\qquad =-\phi(\mf x).
    \end{align*}
    \item When $x\in S_3$
    \begin{align*}
        &\left(-\beta-\mu_M x_1 + \mu_1 x_2\right)f_1^{(2)}(\mf x) + \mu_1 x_2 f_2^{(2)}(\mf x)\\
        &\qquad = \left(-\beta-\mu_M x_1 + \mu_1 x_2\right)\frac{1}{\beta}e^{-\mu_M\tau(\mf x)} + \frac{1}{\beta}\frac{\mu_1^2 x_2e^{-\mu_M\tau(\mf x)}-\mu_1\mu_Mx_2e^{-\mu_1\tau(\mf x)}}{\mu_1-\mu_M}\\ 
        &\qquad= -e^{-\mu_M\tau(\mf x)}-\frac{1}{\beta}\mu_M\left(\frac{\beta}{\mu_M}-\frac{\beta}{\mu_M}e^{-\mu_M\tau(\mf x)}-\frac{\mu_1x_2e^{-\mu_M\tau(\mf x)}-\mu_1x_2e^{-\mu_1\tau(\mf x)}}{\mu_1-\mu_M}\right)\\
        &\qquad\quad + \frac{1}{\beta}\frac{\mu_1^2x_2e^{-\mu_1\tau(\mf x)}-\mu_1\mu_Mx_2e^{-\mu_1\tau(\mf x)}}{\mu_1-\mu_M}+ \frac{1}{\beta}\frac{\mu_1^2 x_2e^{-\mu_M\tau(\mf x)}-\mu_1\mu_Mx_2e^{-\mu_1\tau(\mf x)}}{\mu_1-\mu_M}\\ 
        &\qquad = -1 = -\phi(\mf x).
    \end{align*}
\end{enumerate}
Now, in a similar fashion we check the boundary conditions.
\begin{enumerate}
    \item When $\mf x\in S_0$, the equality is trivial as both derivatives are equal to 0.
    \item When $\mf x\in S_1$, $f_1^{(2)}(0,x_2)=0$. Note that the only point in $S_1$ with $x_1=0$ is $x_2=\kappa_1$ and $f_2^{(2)}(0,\kappa_1)=0$. Hence, the boundary condition is satisfied. 
    \item When $\mf x\in S_2$ and $x_1=0$, $\tau(\mf x)=0$, which implies $f_1^{(2)}(0, x_2)=f_2^{(2)}(0,x_2)=\frac{1}{\beta}\phi(x_2)$. 
    \item Similarly, when $\mf x\in S_3$, $f_1^{(2)}(0, x_2)=f_2^{(2)}(0,x_2)=\frac{1}{\beta}$ 
\end{enumerate}
Hence, the boundary condition is also satisfied for all values of $x_2$.

Now, we prove that the suggested $f^{(1)}(\mf x)$ satisfy the bounds in \eqref{eq:f1_bounds_stationary}-\eqref{eq:f2_der_bounds_stationary} as
\begin{align}\label{eq:f1_bound_kappa}
f^{(1)}(\mf x) &\leq \tau^{\kappa_1}(\mf x) \leq \tau^{\kappa_1}(\kappa_2)-\tau^{\kappa_2}(\kappa_2)\leq \frac{1}{\mu_M}\log\left(\frac{\mu_M\kappa_2-\beta}{\mu_M\kappa_1-\beta}\right) \mbox{ for all } x_1\leq \kappa_2 .
\end{align}
Similarly,  for $f^{(2)}(\mf x)$, we have
\begin{align}
    \label{eq:f2_bound_kappa1}&f^{(2)}(\mf x) \leq \frac{\log(\kappa_2/\kappa_1)}{\mu_1}, \mbox{ for }x_2\leq \kappa_2, \mf x\in S_1,\\
   \label{eq:f2_bound_kappa2} &f^{(2)}(\mf x) \leq \tau(\mf x) +\frac{x_2e^{-\mu_1\tau(\mf x)}-\kappa_1}{\beta}\leq \frac{\log(\kappa_2/\kappa_1)}{\mu_1} + \frac{\kappa_2-\kappa_1}{\beta}, \mbox{ for }x_2\leq \kappa_2, \mf x\in S_2,
\end{align}
where the second inequality in \ref{eq:f2_bound_kappa2} follows as
\begin{align*}
x_2e^{-\mu_1\tau(\mf x)}\geq \kappa_1 \mbox{ and }x_2\leq \kappa_2, \mbox{ when }x_2\in S_2.
\end{align*}
Now, we prove the bounds for $f^{(1)}(\mf x)$ and $f^{(2)}(\mf x)$. As shown in Braverman, we have \[|\phi'(\mf x)|\leq 4/(\kappa_2-\kappa_1)\] and hence,
\begin{align}
    \label{eq:f1_1_bound_kappa}f_1^{(1)}(\mf x) \leq \frac{4\tau^{\kappa_1}(\mf x)}{\kappa_2-\kappa_1}\leq \frac{4\log\left(\frac{\mu_M\kappa_2-\beta}{\mu_M\kappa_1-\beta}\right)}{\kappa_2-\kappa_1} \mbox{ for all } x_1\leq \kappa_2.
\end{align}
\begin{align}
    \label{eq:f2_1_bound_kappa}|f_1^{(2)}(x)|\leq \frac{1}{\beta} \mbox{ for all }\mf x\in \Omega.
\end{align}
Now, we will focus on $f_{11}^{(1)}(\mf x)$ and $f_{11}^{(2)}(\mf x)$.
\begin{align*}
    f_{11}^{(1)}(\mf x )&=\left\{\begin{array}{ll}
    \int_{\tau^{\kappa_2}}^{\tau^{\kappa_1}} e^{-2\mu_Mt}\phi''\left(\frac{\beta}{\mu_M}- \left(x_1 + \frac{\beta}{\mu_M} -\frac{\mu_1x_2}{\mu_1-\mu_M}\left(1-e^{-(\mu_1-\mu_M)t}\right)\right)e^{-\mu_M t}\right)dt,\\
    \qquad\qquad\qquad\qquad\qquad\qquad\qquad\qquad\qquad\qquad\qquad\qquad\qquad\mbox{ for } x_1<-\kappa_2\\
    \int_{0}^{\tau^{\kappa_1}} e^{-2\mu_Mt}\phi''\left(\frac{\beta}{\mu_M}- \left(x_1 + \frac{\beta}{\mu_M} -\frac{\mu_1x_2}{\mu_1-\mu_M}\left(1-e^{-(\mu_1-\mu_M)t}\right)\right)e^{-\mu_M t}\right)dt\\
    \qquad\qquad\qquad\qquad\qquad\qquad\qquad\qquad\qquad\qquad\qquad\qquad\qquad\mbox{ for } x_1\in[-\kappa_2,-\kappa_1]\\
    0 \qquad\qquad\qquad\qquad\qquad\qquad\qquad\qquad\qquad\qquad\qquad\qquad\quad\:\:\;\mbox{ for } x_1\in[-\kappa_1,0]
    \end{array}\right..
\end{align*}
\begin{align*}
    f_{11}^{(2)}(\mf x)&= \left\{\begin{array}{ll}
    0 &\mbox{for } \mf x\in S_0\cup S_1\\
    -\frac{1}{\beta}\mu_M\tau_1(\mf x)\phi(x_2e^{-\mu_1\tau(\mf x)}) - \frac{1}{\beta}e^{-(\mu_1+\mu_M)\tau(\mf x)}\phi'(x_2e^{-\mu_1\tau(\mf x)})x_2\mu_1\tau_1(\mf x) &\mbox{for }\mf x\in S_2\\
    -\frac{1}{\beta}\mu_M\tau_1(\mf x) &\mbox{for }\mf x\in S_3.
    \end{array}\right.
\end{align*}
Again, as shown in Braverman $|\phi''(\mf x)|\leq 12/(\kappa_2-\kappa_1)$, and hence, 
\begin{align}
   \label{eq:f1_11_bounds_kappa} |f_{11}^{(1)}(\mf x)|\leq \frac{12\log\left(\frac{\mu_M\kappa_2-\beta}{\mu_M\kappa_1-\beta}\right)}{\mu_M(\kappa_2-\kappa_1)^2}.
\end{align}
We also have
\begin{align}
    \label{eq:f2_11_bounds_kappa}|f_{11}^{(2)}(\mf x)| &\leq \frac{\mu_M}{\beta(\mu_1x_2-\beta)}+\frac{4}{\beta(\kappa_2-\kappa_1)} \mbox{ for all }x_2\geq \kappa_1
\end{align}
Choosing $\kappa_2=\frac{\beta}{\mu_M}+2\epsilon$ and $\kappa_1=\frac{\beta}{\mu_M}+\epsilon$ and realising $\mu_M\leq 1\leq \mu_1$, \eqref{eq:f1_bound_kappa}-\eqref{eq:f2_11_bounds_kappa} yields \eqref{eq:f1_bounds_stationary}-\eqref{eq:f2_der_bounds_stationary}. 
\endproof

{\em Proof of Theorem~\ref{thm:stationarity_diffusion}.}
Lemma~\ref{lem:stationary_pde_solution} implies that for any $\epsilon>0$, we can choose $f^{(1)}(\mf x)$ and $f^{(2)}(\mf x)$ solving \eqref{eq:pde_stationary1} and \eqref{eq:pde_stationary2}, respectively, such that
\begin{align*}
    \alpha\left(-1 + f_{11}^{(\Sigma)}(\mf x) + \alpha \left(f_1^{(\Sigma)}(\mf x)\right)^2\right)&\leq \alpha\left(-1 + \frac{12\log2}{\mu_M\epsilon^2} + \frac{\mu_M+4}{\beta\epsilon} + \alpha \left(\frac{4\log2}{\epsilon}+\frac{1}{\beta}\right)^2\right),
\end{align*}
and 
\begin{align*}
    \left(\alpha\left(f_{11}^{(\Sigma)}(\mf x)\right) + \alpha^2 \left(f_1^{(\Sigma)}(\mf x)\right)^2\right)&\leq \left(\alpha\left(\frac{12\log2}{\mu_M\epsilon^2} + \frac{\mu_M+4}{\beta\epsilon}\right) + \alpha^2 \left(\frac{4\log2}{\epsilon}+\frac{1}{\beta}\right)^2\right).
\end{align*}
Choosing $\alpha$ sufficiently small and $\epsilon$ sufficiently large, we can find constants $c,d>0$ to satisfy \eqref{eq:drift_bound} and \eqref{eq:drift_bound_conditional}, which proves the first part of the theorem. Then, the positive recurrence of the diffusion follows from Theorem 5.2 in \cite{down_etal95}. 
\qed
\endproof


\end{document}